\documentclass[]{IEEEtran}

\usepackage[utf8]{inputenc}
\usepackage[backend=biber, sorting=none]{biblatex}
\usepackage{microtype}
\usepackage{amsmath}
\usepackage{amsfonts}
\usepackage{amssymb}
\usepackage{mathtools}
\usepackage[version=4]{mhchem}
\usepackage{textcomp}
\usepackage{booktabs}
\usepackage{colortbl}

\usepackage{todonotes}
\usepackage{regexpatch}
\makeatletter
\xpatchcmd{\@todo}{\setkeys{todonotes}{#1}}{\setkeys{todonotes}{inline,#1}}{}{}
\makeatother

\makeatletter
\let\MYcaption\@makecaption
\makeatother

\usepackage[font=footnotesize]{subcaption}

\makeatletter
\let\@makecaption\MYcaption
\makeatother

\usepackage{graphicx}
\usepackage{siunitx}
\usepackage{tabularx}
\usepackage{hyperref}
\usepackage[capitalize]{cleveref}

\title{Cold-Start Modeling and On-Line Optimal Control of the Three-Way Catalyst}

\author{Jonathan Lock, Kristoffer Clasén, Jonas Sjöblom, Tomas McKelvey}
\begin{document}

\maketitle

\begin{abstract}
We present a three-way catalyst (TWC) cold-start model, calibrate the model based on experimental data from multiple operating points, and use the model to generate a Pareto-optimal cold-start controller suitable for implementation in standard engine control unit hardware. The TWC model is an extension of a previously presented physics-based model that predicts carbon monoxide, hydrocarbon, and nitrogen oxides tailpipe emissions. The model axially and radially resolves the temperatures in the monolith using very few state variables, thus allowing for use with control-policy based optimal control methods. In this paper we extend the model to allow for variable axial discretization lengths, include the heat of reaction from hydrogen gas generated from the combustion engine, and reformulate the model parameters to be expressed in conventional units. We experimentally measured the temperature and emission evolution for cold-starts with ten different engine load points, which was subsequently used to tune the model parameters (e.g.~chemical reaction rates, specific heats, and thermal resistances). The simulated cumulative tailpipe emission modeling error was found to be typically -20\% to +80\% of the measured emissions. We have constructed and simulated the performance of a Pareto-optimal controller using this model that balances fuel efficiency and the cumulative emissions of each individual species. A benchmark of the optimal controller with a conventional cold-start strategy shows the potential for reducing the cold-start emissions.
\end{abstract}

\nocite{lock2021}

\section{Introduction}

The Three-Way Catalyst (TWC) is used in nearly all conventional vehicles with spark-ignited (SI) engines to reduce the level of harmful emissions generated by the combustion engine that would otherwise exit the tailpipe. The toxic emissions generated by the combustion engine, broadly categorized as nitrogen oxides (\ce{NO_x}), carbon monoxide (\ce{CO}), and residual hydrocarbons (\ce{THC}), are in the TWC converted to primarily form non-toxic nitrogen gas (\ce{N_2}), carbon dioxide (\ce{CO_2}), and water (\ce{H_2O}) \cite{heywood1988internal, hedinger_optimal_2017}. Modern TWCs are very effective at removing emissions, with conversion efficiencies of over 95\% being commonplace \cite{hedinger_optimal_2017} and in some cases significantly higher, as was found in the experimental results in this paper. However, the TWC must be sufficiently hot to function, which is in ordinary operation maintained by virtue of the hot exhaust gases passing through it from the combustion engine. However, when a vehicle is cold-started (i.e.~started after the TWC has had sufficient time to cool to the ambient temperature) the tailpipe emissions are much larger until the TWC is heated to its ordinary operation temperature, an interval typically taking on the order of 40--100 seconds\cite{hedinger_optimal_2017}. These cold-start emissions are very significant, and for many regulatory test procedures are responsible for 60--80\% of the emissions generated from an entire test (which are for reference on the order of 30 minutes).

Several methods for reducing cold-start emissions have been studied from a multiple perspectives. These range from methods of constructing the TWC that reduces the cold-start time \cite{ramanathan_optimal_2004}, methods for preheating the TWC before starting the combustion engine \cite{gao_review_2019}, and control schemes that focus on controlling the combustion engine's operation to limit the emissions generated during the cold-start \cite{hedinger_optimal_2017, schori_optimal_2014, pannag_coldstart_2009, fiengo_control_2002, keynejad_suboptimal_2013}.

In this paper our goal is to develop a model-based optimal TWC cold-start controller that can be feasibly be implemented in a standard engine control unit (ECU). We will consider a conventional SI engine and TWC, where the engine's load point can be freely controlled during the cold-start, making the controller suitable for e.g.~hybrid vehicles. More specifically, we view the TWC cold-start problem as determining the optimal engine speed, load, and spark timing to apply over time while balancing the conflicting goals of maximizing fuel efficiency and minimizing the cumulative emissions. Generating a Pareto-optimal controller will therefore both require a dynamic thermal model of the TWC as well as a suitable optimal control design method. This paper can naturally be divided into two parts, one where we develop a TWC cold-start model, and one where we evaluate the performance of an optimal controller generated using said model.

The model presented in this paper extends on a TWC model previously developed by the authors \cite{lock2021}. In the previous work, we derived a thermal model of the TWC with very few state variables suited for fast off-line simulation or on-line control systems. The model resolved both axial and radial temperature variations from a first-principles perspective. In this paper we extend the previous model by allowing the axial discretization to vary along the TWC's length, model the heat generation caused by oxidation of hydrogen in the exhaust gas, and reformulate the tuning parameters (heat capacity, thermal conductivity, etc) to be expressed in well-known units (\si{\joule \per \kelvin \per \kilogram}, \si{\watt \per \meter \per \kelvin}, etc). Furthermore, we expand on our previous work by here considering a TWC consisting of two separate monoliths and use separate training and validation datasets for tuning and evaluation. 

In this paper we will first briefly discuss categories of existing models and their strengths, weaknesses, and relevant applications. Following this we will introduce our extended model and the experimental setup used to calibrate the model. Finally, we will study the experimental results and evaluate the simulated performance of the optimal control scheme using the calibrated model.  A listing of all abbreviations is shown in \cref{tab:abbrev}, while all model parameters and their units are presented in \cref{tab:twc-params}.

\begin{table}
	\centering
	\caption{Table of used abbreviations.} \label{tab:abbrev}
	\begin{tabular}{c|c}
		\hline
		CAbTDC & Crank Angle before Top Dead Center\\
		\hline
		\ce{CO} & Carbon Monoxide \\
		\hline
		ECU & Engine Control Unit \\
		\hline
		SI & Spark Ignited \\
		\hline
		\ce{NO_x} & Nitrogen Oxides \\
		\hline
		SA & Spark Angle \\
		\hline
		\ce{THC} & Total Hydrocarbon \\
		\hline
		TWC & Three Way Catalyst \\
		\hline
		TWC1 & The first monolith in the TWC\\
		\hline
		TWC2 & The second monolith in the TWC\\
		\hline
	\end{tabular}
\end{table}

\begin{table}
	\centering
	\caption{TWC parameters.\label{tab:twc-params}}
	\begin{tabularx}{\columnwidth}{l|l|X}
		\footnotesize{Parameter} & \footnotesize{Unit} & \footnotesize{Description} \\
		\hline
		\hline
		$L$ & \si{\meter} & Total TWC length \\ 
		\hline 
		$R$ & \si{\meter} & TWC radius \\ 
		\hline
		$L_n$ & \si{\meter} & Length of axial slice $n$ \\ 
		\hline 
		$L_{n-1,n}$ & \si{\meter} & Center-to-center distance between axial slice $n-1$ and $n$ \\ 
		\hline 
		$t_w$ & \si{\meter} & TWC wall thickness \\ 
		\hline 
		$l_c$ & \si{\meter} & TWC channel width \\ 
		\hline
		$\mathrm{OFA}$ & - & TWC open frontal area\\
		\hline
		$m_\mathrm{TWC}$ & \si{\kilogram} & Monolith mass \\
		\hline 
		$c_p$ & \si{\joule \per \kelvin \per \kilogram} & Net monolith specific heat\\
		\hline 
		$\dot{m}_\mathrm{exh}$ & \si{\kilogram \per \second} & Exhaust massflow \\
		\hline 
		$c_{p,\mathrm{exh}}$ & \si{\joule \per \kelvin} & Exhaust gas specific heat\\
		\hline
		$T_\mathrm{amb}$ & \si{\degreeCelsius} & Ambient temperature\\
		\hline
		$T_\mathrm{exh}$ & \si{\degreeCelsius} & Temperature of gas feeding TWC\\
		\hline
		$k_\mathrm{ax}$ & \si{\watt \per \meter \per \kelvin} & Effective axial thermal conductivity\\
		\hline
		$k_\mathrm{rad}$ & \si{\watt \per \meter \per \kelvin} & Effective radial thermal conductivity\\
		\hline 
		$k_\mathrm{amb}$ & \si{\watt \per \meter \per \kelvin} & Effective thermal conductivity to ambient\\
		\hline  
		$t_\mathrm{amb}$ & \si{\meter} & Thickness of insulation to ambient\\
		\hline 
		\hline
		$N$ & - & Number of axial TWC segments \\ 
		\hline
		$M$ & - & Number of resolved radial TWC channels \\ 
		\hline 
		$T$ & \si{\degreeCelsius} & Temperature state vector\\		
		\hline 
		$T_n$ & \si{\degreeCelsius} & Central temperature in axial slice $n$\\ 
		\hline 
		$\Delta_T$ & \si{\degreeCelsius} & Center/periphery temperature difference\\ 
		\hline 
		\hline
		$k_{n,m}^s$ & \si{\mole \per \second} & Reaction rate, species $s$ in cell $n,m$ \\ 
		\hline 
		$A^s$ & - & Pre-exponential factor, species $s$ \\ 
		\hline 
		$E_a^s$ & \si{\joule \per \mole} & Activation energy of emission species $s$ \\ 
		\hline 
		$y_{n,m}^s$ & - & Mole fraction, species $s$ in cell $n,m$ \\ 
		\hline 
		$\dot{m}_{n,m}^{\mathrm{s,conv}}$ & \si{\kilogram \per \second} & Mass conversion, species $s$, cell $n,m$\\
		\hline
		$\dot{m}_{n,m}^{\mathrm{s,in}}$ & \si{\kilogram \per \second} & Incoming mass, species $s$, cell $n,m$\\
		\hline
		$\dot{m}_{n,m}^{\mathrm{s,conv}}$ & \si{\kilogram \per \second} & Outgoing mass, species $s$, cell $n,m$\\
		\hline
		$\dot{m}^s_\mathrm{tp}$ & \si{\kilogram \per \second} & Tailpipe massflow, species $s$\\
		\hline
		$t_{r,n}$ & \si{\second} & Gas residence time in slice $n$ \\
		\hline 
		$p_\mathrm{TWC}$ & \si{\pascal} & Absolute pressure in TWC\\
		\hline 
		\hline
		$P_\mathrm{ctr}$ & \si{\watt} & Heating power, radial center\\
		\hline 
		$P_\mathrm{per}$ & \si{\watt} & Heating power, radial periphery\\
		\hline 
		$P_\mathrm{ax}$ & \si{\watt} & Heating power, axial conduction\\
		\hline 
		$P_\mathrm{rad}$ & \si{\watt} & Heating power, radial conduction\\
		\hline 
		$P_\mathrm{con,ctr}$ & \si{\watt} & Heating power, convection, center \\
		\hline 
		$P_\mathrm{con,per}$ & \si{\watt} & Heating power, convection, periphery\\
		\hline 
		$P_{n,m}$ & \si{\watt} & Exothermic power generated in cell $n,n$\\
		\hline 
		$P_\mathrm{exo,ctr}$ & \si{\watt} & Weighted central exothermic power\\
		\hline 
		$P_\mathrm{exo,per}$ & \si{\watt} & Weighted peripheral exothermic power\\
		\hline 
		$P_\mathrm{amb}$ & \si{\watt} & Heating power, loss to ambient\\
		\hline 
		$q_\mathrm{ax}$ & \si{\watt \per \meter \squared} & Axial heat flux\\
		\hline 
		$q_\mathrm{rad}$ & \si{\watt \per \meter \squared} & Radial heat flux\\
		\hline 
		$q_\mathrm{amb}$ & \si{\watt \per \meter \squared} & Heat flux to ambient\\
		\hline
	\end{tabularx} 
\end{table}

\subsection{Literature survey}

The topic of modeling the TWC for cold-start purposes has been considered from a wide range of perspectives. Some authors derive fully three-dimensional or two-dimensional models that capture many of the fundamentally complex chemical kinetics, transport dynamics, and temperature dynamics \cite{braun2002three, chen1988three, zygourakis1989transient}. Though accurate, models of this caliber are very computationally demanding and primarily suited for in-depth analysis.

One commonly used method of reducing the computational demand is to exploit the characteristic structure of modern TWC's, which consist of a large number of parallel channels. Assuming each channel's construction and composition is identical it is sufficient to study the behavior of a single channel and afterwards scale the result by the number of channels in the TWC. This approximation is typically referred to as the single channel approximation or single channel model, and results in models with significantly reduced computational demands.

Several authors\cite{pontikakis2004three, baba1996numerical, oh1982transients, real_modelling_2021, yan_modeling_2019, ramanathan_kinetic_2011, brandt_dynamic_2000} model the TWC with the single channel approximation, and depending on the number of modeled chemical phenomenon and the level of spatial resolution some models can approach realtime simulation speeds on powerful PC's (i.e.~where generating one second's worth of data requires one second of computation time). Beyond modeling the temperature dependence of the TWC, many of these models also include terms to capture a TWC phenomenon where oxygen is absorbed and released in the TWC. The stored oxygen greatly influences the TWC's capacity to convert emissions, where \ce{CO} and \ce{THC} are more effectively oxidized when the stored level of oxygen is large, while \ce{NO_x} is more effectively reduced when the level is low.

These models allow for simulating either the whole TWC or a representative channel to a fairly high degree of accuracy. However, they are primarily of use for analyzing the performance of the TWC in an off-line manner, for instance in a TWC design process. In this paper we instead focus on on-line TWC control, and in particular consider the cold-start problem. This requires a model that is significantly simpler computationally, both as optimal control methods place specific demands on the complexity and structure of the system models, and as the ECU has a very limited computational capacity. There are several classes of control-oriented models, ranging from models that approximate the spatially varying temperature distribution as a scalar temperature \cite{shaw_simplified_2002, schori_optimal_2014, pannag_coldstart_2009, fiengo_control_2002}, to more complex models that axially resolve the TWC temperature and/or include oxygen-storage terms \cite{michel_optimizing_2017, azad_determining_2012, zhu_development_2019}. The simpler models are fairly well suited for direct use with on-line optimal control methods, while most of the more complex models are used primarily as a starting point for creating a suboptimal controller.

The model presented in this paper (an extension of \cite{lock2021}) has been constructed for the specific purpose of subsequently being implemented for optimal control in conventional ECUs. More specifically, it is well-suited to control-policy based optimal control methods, where the optimal control signal (e.g.~engine speed, load, spark angle, and so on) is precomputed in an offline phase and stored in a table for a discrete set of TWC temperatures. A subsequent realtime controller can ultimately determine the optimal control signal by consulting the table of stored temperatures and associated optimal control signals \cite{Bertsekas2017, bryson1975applied}. Controllers of this class are very powerful, as they allow for nearly arbitrarily nonlinear model dynamics, costs, and constraints, but are limited in that their memory demand scales exponentially with the number of state variables and require the states variables to either be measured or estimated. In an effort to limit the number of state variables we have chosen to not dynamically model the stored oxygen in the TWC as this is not as significant as the temperature dynamics during a cold-start \cite{shaw_simplified_2002, lock2021}.

\section{TWC Model}
We will in this section introduce the TWC model. This model is based on and extends a model previously presented by the authors \cite{lock2021}. The model is extended by allowing the size of the axially discretized slices to vary over the length of the catalyst, reformulated so all parameters are based on easily determined physical parameters, and we consider the case where two separate TWC monoliths are placed in series.

The model can naturally be divided into three distinct subsections; one modeling the chemical kinetics, one modeling the temperature dynamics, and one that interpolates the low-dimensional state variables to a higher-dimensional temperature distribution using a physics-based method. We will initially consider a single TWC monolith, and later return to the case where two are placed in series.

\begin{figure}
	\centering
	\begin{subfigure}[b]{\columnwidth}
		\includegraphics[width=\textwidth]{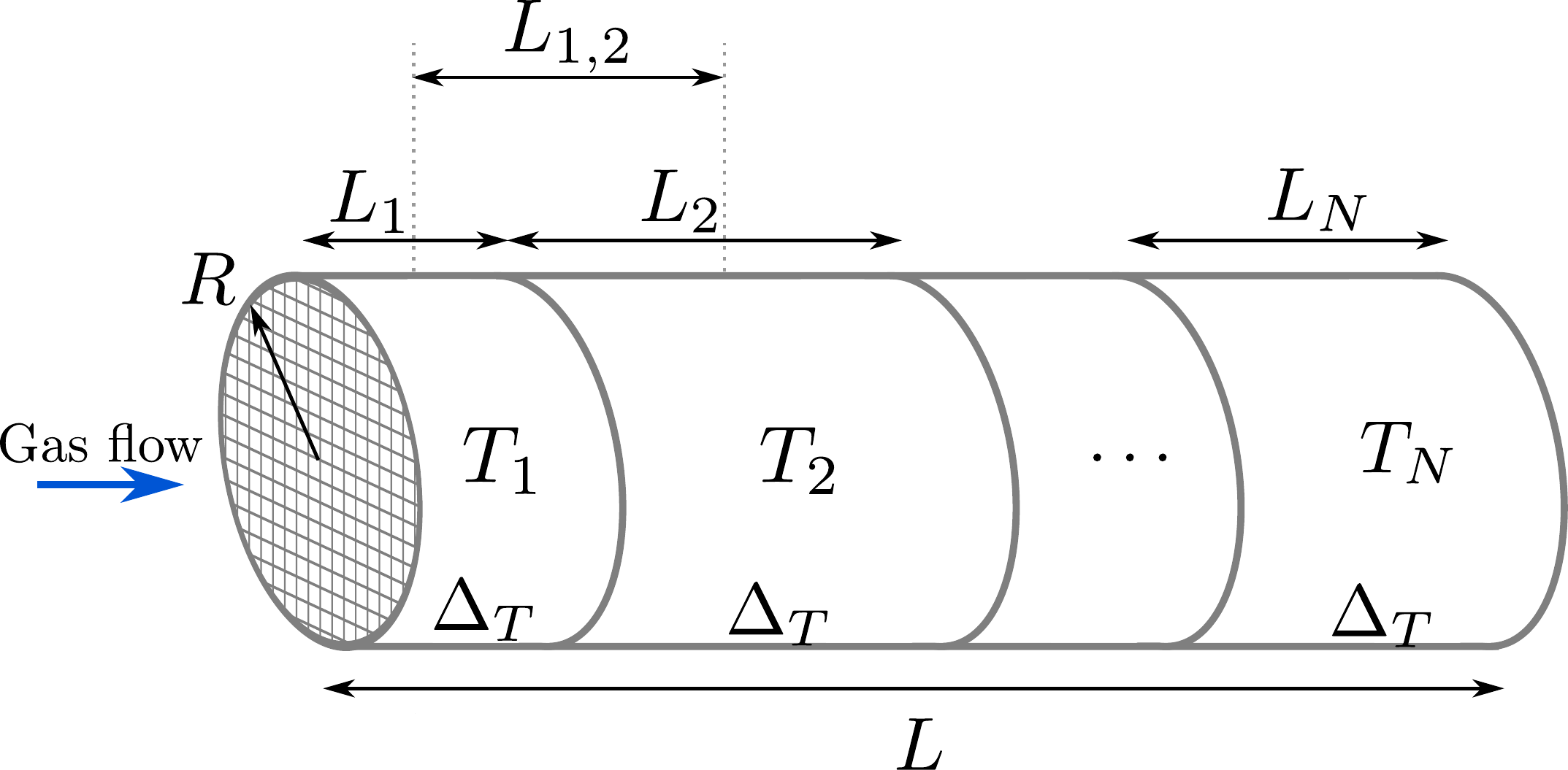}
		\caption{TWC body. Engine exhaust travels from left to right. Figure derived from \cite{lock2021}. \label{fig:twc-geometry-body}}
	\end{subfigure}\\
	\vspace{\baselineskip}
	\begin{subfigure}[b]{0.35\columnwidth}
		\includegraphics[width=\textwidth]{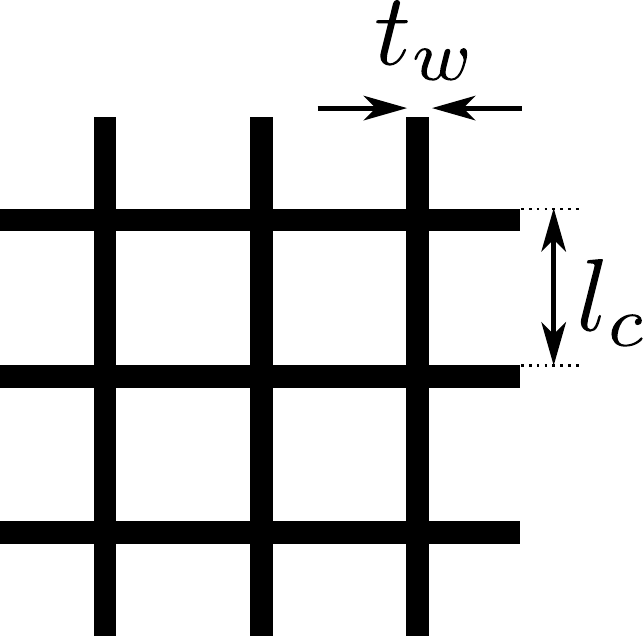}
		\caption{Detail of TWC channels.\label{fig:twc-geometry-channels}}
	\end{subfigure}
	\caption{Illustration of geometry and measurement definitions for a single TWC monolith.\label{fig:twc-geometry}}
\end{figure}

An illustration of the assumed TWC geometry is shown in \cref{fig:twc-geometry-body,fig:twc-geometry-channels}. Specifically, we assume that the TWC is cylindrical with radius $R$ and length $L$. We make the modeling choice of dividing the TWC into $N$ different axial slices, and extend the previously presented model \cite{lock2021} that assumed equally-sized slices by allowing the associated lengths $L_1, L_2, \dots, L_N$ (where $\sum L_n = L$) to be different for each slice. We also define the lengths $L_{1,2}, L_{2,3}, \dots, L_{N-1,N}$ as the axial distances between the midpoints of neighboring slices. Finally, we assume the TWC has a monolithic structure with square, axially traversing channels of wall thickness $t_w$ and channel length $l_c$, as illustrated in \cref{fig:twc-geometry-channels}.

Importantly, we make the approximation of axially discretizing the TWC temperature. More specifically, we model the TWC temperature in the radial center of slice $n$ as $T_n$. This implies that the temperature in the radially central channel is modeled as $N$ segments of constant temperature. Furthermore, we model the difference in temperature between the radial center and radial periphery of each slice as as $\Delta_T$. Note that $\Delta_T$ is not axially resolved, i.e.~we assume that $\Delta_T$ is identical for all axial slices, i.e.~the periphery temperature of slice 1 is $T_1 + \Delta_T$, slice 2 is $T_2 + \Delta_T$, and so on. Finally, as is described in more detail in \cref{subsec:rad-temp-interp}, we use an interpolation scheme to approximate the temperature at $M$ different radial locations ranging from the radial center to the periphery. Ultimately, a single TWC monolith is at any given instance in time characterized by the state variable vector
\begin{align}
	T &= \begin{bmatrix}
		T_1\\
		T_2\\
		\dots\\
		T_N\\
		\Delta_T
	\end{bmatrix}\,.
\end{align}

In the following subsections we will detail the individual parts of the model. First, we define how the state variable is used to generate a dense representation of the TWC temperature that is axially and radially resolved. This is followed by the chemical kinetics model that determines the conversion efficiency of the TWC as well as the heating power generated by the exothermic reactions. Finally, we introduce the thermal model that generates the state variable derivative and define the interface between the two separate TWC monoliths.

\subsection{Radial temperature interpolation}\label{subsec:rad-temp-interp}
The model presented in this work uses the single-channel approximation for fundamental TWC material properties, while resolving the radial temperature profile by simulating several parallel channels corresponding to different radial positions with different associated temperatures. This allows for capturing the experimentally observed behavior where the periphery of the TWC is significantly colder than the radially central sections (as we will be discussed in \cref{subsec:cold-start}).

\begin{figure}
	\centering
	\includegraphics[width=\columnwidth]{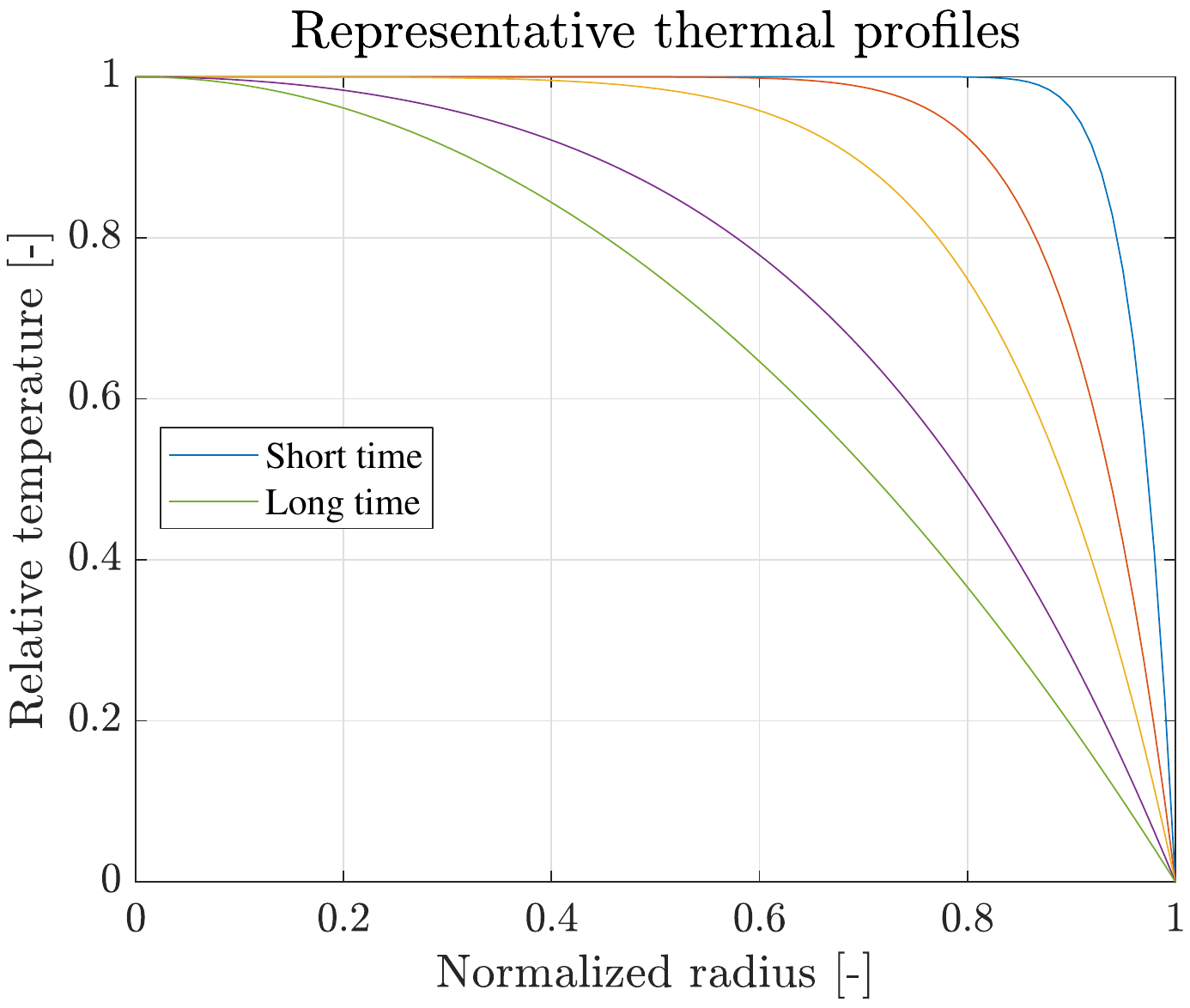}
	\caption{Representative solutions to the transient heat equation in a flat disc with a homogeneous power term and a Dirichlet boundary condition $\hat{T}(R) = 1$. The displayed temperatures and radii have been normalized to the range $[0,1]$. Figure reused from \cite{lock2021}. \label{fig:heat-eq-sol}}
\end{figure}

More specifically, we model the radial temperature profile $\hat{T}(t,r)$, at time $t$ and radius $r$, as a solution to the transient heat equation in a flat circular disc with radius $R$ and an initial temperature of zero, i.e.~$\hat{T}(0,r) = 0, r=[0,R]$. Note that $\hat{T}$ has no relation to the state variable $T$ or the axial slice temperatures $T_n$ despite the similar notation. Furthermore, we assume a Dirichlet boundary condition, i.e.~$\hat{T}(t,R) = 0$, and assume that the plate develops a constant homogeneous power. This power is intended to be analogous to the power delivered to a slice in the TWC by convection, axial conduction, and exothermic heat generation. Solving the time-evolution of $\hat{T}(t,r)$ is a textbook problem (e.g. \cite[p. 148]{logan1998applied}) with a solution that can be expressed as a Fourier-Bessel series. Solving this numerically over time and radii can also be done easily (for instance with MATLAB's \texttt{pdepe} function), generating the radially-varying time evolution of the plate's temperature. Normalized solutions (where $\hat{T}(t,r)$ is scaled to the range $[0,1]$, and $R=1$) are shown in \cref{fig:heat-eq-sol}.

\begin{figure}
	\centering
	\includegraphics[width=\columnwidth]{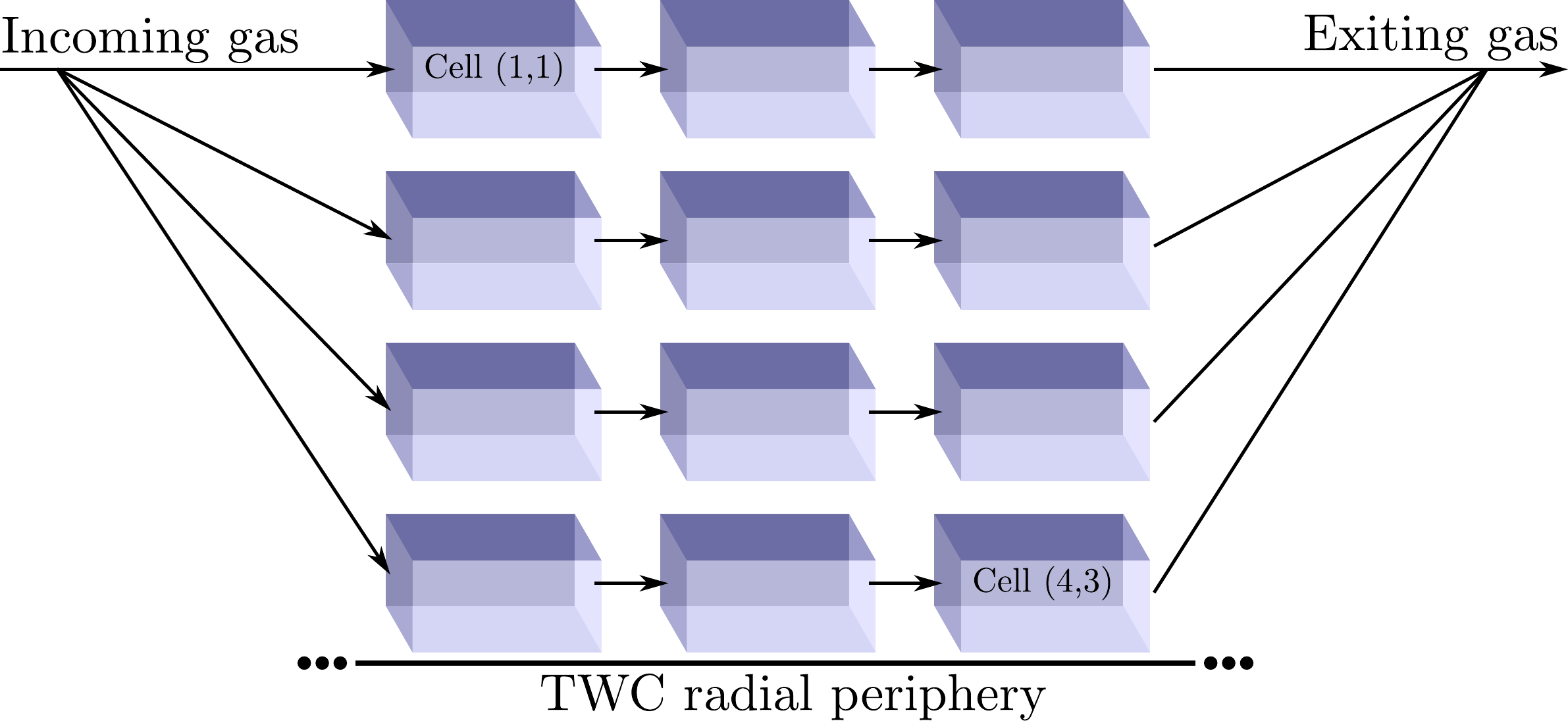}
	\caption{Fully-resolved TWC, here shown for $N=3$ and $M=4$. Figure reused from \cite{lock2021}. \label{fig:fully-resolved-temps}}
\end{figure}

In this paper we interpolate the radial TWC temperature profile by using precomputed solutions to the above flat-plate problem. More specifically, we assume that the radial temperature profile in slice $n$ is given by $A_1 \hat{T}(t',r) + A_2$ for a given time instant $t'$.  Here, $A_1$ and $A_2$ are selected so that $\hat{T}(t',0) = T_n$ and $\hat{T}(t',R) = T_n + \Delta_T$ (i.e.~$\hat{T}$ is scaled and offset to match the known radial center and periphery temperatures), and $t'$ is selected to give a radial temperature profile that matches the experimentally measured temperature profile for a given engine operating point (more on this in \cref{subsec:steady-state}). In summary, we interpolate the radial temperature distribution using the known radially central and peripheral temperatures and the instantaneous engine operating point. Letting $M$ denote the number of independent single-channel models we wish to resolve, this interpolation scheme allows us to convert the $N+1$ state variables to a representation with $M$ single-channel models of $N$ segments. This is illustrated in \cref{fig:fully-resolved-temps}, where cells $([1 \dots M],1)$ are fed with the incoming gas, cell $(m,n+1)$ is fed with the output of cell $(m,n)$, and the output from cells $([1 \dots M], N)$ are combined to form the total exhaust from the TWC. 

As the TWC is assumed to be circular, regions near the periphery have a larger associated area than regions near the axial center. This is taken into account in our model by approximating the massflow in the physical TWC as equal at all locations, and scaling the proportion of gas passing through each channel to match. Letting $\dot{m}_\mathrm{exh}$ be the total exhaust massflow from the engine and scaling by the relative area of an annular ring with a major radius of $m/M$ and a minor radius of $(m-1)/M$ gives the massflow into a given cell as
\begin{subequations}
	\begin{align}
	\dot{m}_{m,1} &= \frac{\pi (m^2 - (m - 1)^2)}{\pi M^2} \dot{m}_\mathrm{exh} \nonumber \\
	&= \frac{(2m - 1)}{M^2} \dot{m}_\mathrm{exh} \label{eq:weighting}\\
	\dot{m}_{m,n+1} &= \dot{m}_{m,n}\,.
	\end{align}
\end{subequations}

\subsection{Chemical kinetics}
The total range of chemical reactions occurring in the TWC are very complex and involve a wide range of compounds. However, there are fewer that contribute to the legislated emissions or significant heat generation. We will therefore limit our scope to net reactions (i.e.~without considering intermediary steps). This is done both for simplicity, and as a detailed approach would require the addition of numerous state variables that track the concentration of the emission species and their intermediaries in the TWC. \cite{brandt_dynamic_2000, pattas_transient_1994, lock2021} give the most significant reactions (apart from \cref{eq:h2}) as
\begin{subequations}
	\begin{align}
	\ce{2CO + O_{2}} & \ce{ -> 2CO_{2}}\label{eq:co}\\
	\ce{2H_{2} + O{2}} & \ce{ -> 2H_{2}O} \label{eq:h2} \\
	\ce{2NO + 2CO} & \ce{-> N_{2} + 2CO_{2}}\label{eq:no}\\
	\ce{2NO_{2}} & \ce{-> N_{2} + 2O_{2}}\label{eq:no2}\\
	\ce{2C_{3}H_{6} + 9O_{2}} & \ce{-> 6CO_{2} + 6H_{2}O}\label{eq:propene}\\
	\ce{C_{3}H_{8} + 5O_{2}} & \ce{-> 3CO_{2} + 4H_{2}O}\,.\label{eq:propane}
	\end{align}\label{eq:reactions}
\end{subequations}
We will assume all reactions are balanced, i.e.~for a TWC with 100\% conversion efficiency all \ce{CO}, \ce{H_2}, \ce{C_3H_6}, and \ce{C_3H_8} emitted from the engine are fully oxidized and all \ce{NO} and \ce{NO_2} are fully reduced.

In this paper we include the heat generated by the oxidation of hydrogen gas, i.e.~\eqref{eq:h2}, which is generated in the engine by the water-gas shift reaction. As \ce{H_2} is not typically experimentally measured we will instead estimate its mole-fraction from the measured mole-fractions of \ce{CO} and \ce{CO_2}, as given by \cite[eq. (4.68)]{heywood1988internal}
\begin{subequations}
	\begin{align}
		y_{\ce{H_2}} &= \frac{y_{\ce{H_{2}O}} y_{\ce{CO}}}{K y_{\ce{CO_{2}}}}\\
		\shortintertext{where}
		y_{\ce{H_2O}} &= \frac{m'}{2n'} \frac{y_{\ce{CO}} + y_{\ce{CO_{2}}}}{1 + y_{\ce{CO}}/(K y_{\ce{CO_{2}}}) + (m'/2n')(y_{\ce{CO}} + y_{\ce{CO_{2}}})}\,.
	\end{align} \label{eq:h2-est}
\end{subequations}
Here, $y$ corresponds to the mole-fraction of each corresponding compound, $K$ is a constant value set to 3.8 \cite[eq (4.63)]{heywood1988internal}, and $n'$ and $m'$ correspond to the number of carbon and hydrogen atoms respectively\footnote{\cite{heywood1988internal} indicates in (4.68) that $n'$ and $m'$ correspond to hydrogen and carbon respectively. Consulting the previous derivations instead indicates that $n'$ and $m'$ correspond to carbon and hydrogen.} in each molecule of the fuel. With the RON95 E10 fuel studied in this paper we used the supplier-specified value of $m'/2n' = 0.258$.

Typically, nitrogen oxides (\ce{NO} and \ce{NO_2}) and hydrocarbon (\ce{C_3H_6} and \ce{C_3H_8}) emissions are lumped together and denoted as \ce{NO_x} and \ce{THC} respectively \cite[pp. 572--597]{heywood1988internal}. We will in this paper assume a constant ratio of 99:1 for \ce{NO} to \ce{NO_2} as indicated by \cite[p. 578]{heywood1988internal}, and by \cite{tischer_three-way-catalyst_2007, svraka_model_2019} a constant ratio of 3:1 for \ce{C_3H_6} to \ce{C_3H_8}.

We model the reaction rate $k_{n,m}^s$ of an emissions species $s$ in any given cell $n,m$ using an Arrhenius expression of form
\begin{align}
	k_{n,m}^s = A^se^\frac{-E_a^s}{RT_{n,m}}\,,\label{eq:arrhenius}
\end{align}
where $R$ is the ideal gas constant, $T_{n,m}$ is the temperature of cell $n,m$, $E_a^s$ is the activation energy of emission species $s$, and $A^s$ is the apparent pre-exponential factor for species $s$. Letting $y_{n,m}^s$ indicate the mole fraction of emission species $s$ in cell $n,m$, we model the evolution of the mole fraction as
\begin{align}
	\frac{\mathrm{d}y_{n,m}^s}{\mathrm{d}t} &= -k_{n,m}^s y_{n,m}^s\,. \label{eq:kin-diff-eq}
\end{align}
Note that we do not include an inhibition factor in \cref{eq:kin-diff-eq} in order to limit the complexity of the model and the model tuning process. However, including an inhibition term (e.g.~as in \cite{ramanathan_kinetic_2011}) is viable and would not interfere with the optimal control method that will be described in \cref{sec:opt-ctrl}. Furthermore, note that \cref{eq:kin-diff-eq} does not include an oxygen concentration term. As the engine is operated stoichiometrically the \ce{O_2} concentration is fairly constant, implying that it can be lumped into $k_{n,m}^s$. This is beneficial as we avoid the need to explicitly measure or model the \ce{O_2} concentration.

Though \cref{eq:h2-est} allows for generating an estimate of the hydrogen gas concentration for a given \ce{CO} and \ce{CO_2} concentration (quantities which are easily measured with conventional emissions-measurement equipment), as $y^{\ce{H_{2}}}$ is not typically measured it is difficult to tune the associated reaction rate parameters. By \cite{ramanathan_kinetic_2011, svraka_model_2019} we have chosen to instead model the reaction rate of \cref{eq:h2} as identical to that of \cref{eq:co}, i.e.~we assume $E_a^{\ce{H_2}} = E_a^{\ce{CO}}$ and $A^{\ce{H_2}} = A^{\ce{CO}}$. Finally, as \ce{H_2} is not typically viewed as a problematic emission species we will in the remainder of this paper only consider \cref{eq:h2} from the perspective of determining the heat of reaction, in contrast to \ce{CO}, \ce{THC}, and \ce{NO_x} emissions which both contribute with their associated heat of reaction and whose tailpipe emissions are important to track.

As described in more detail in \cite{lock2021}, the gas residence time in each axial slice is short enough for the monolith temperature to be close-to constant. Using this constant-temperature approximation we can explicitly solve \cref{eq:kin-diff-eq} as
\begin{align}
	y_{n,m}^s(t_{r,n}) = y_{n,m}^s(0)e^{-k_{n,m}^s t_{r,n}}\label{eq:conc-exp}
\end{align}
for a residence time in slice $n$ of $t_{r,n}$. We will approximate the residence time by assuming a plug-flow reactor model (i.e.~assuming there is no axial dispersion), giving
\begin{align}
	t_{r,n} &= \frac{V_{\mathrm{slice},n}}{\nu} \label{eq:t_r-base}
\end{align}
where $V_{\mathrm{slice},n}$ is the gas volume of slice $n$ and $\nu$ is the volumetric flow-rate of the exhaust gases. Using the geometry of the TWC as defined in \cref{fig:twc-geometry}, and using the ideal gas law we can approximate \cref{eq:t_r-base} as
\begin{align}
	t_{r,n} &= \frac{\mathrm{OFA} \cdot L_n \pi R^2}{\dot{m}_\mathrm{exh}R_\mathrm{specific}T_{n,m}p_\mathrm{TWC}^{-1}}\,,
\end{align}
where $\mathrm{OFA}$ is the open frontal area of the TWC, defined by
\begin{align}
	\mathrm{OFA} &= (l_c-t_w)^2l_c^{-2}\,.\label{eq:ofa}
\end{align}
In this paper we extend \cite{lock2021} to use physically meaningful SI units for all parameters. Here $\dot{m}_\mathrm{exh}$ is the exhaust massflow (\si{\kilogram \per \second}), $P$ is the absolute pressure (\si{\pascal}) in the TWC, which is typically close to the ambient pressure, and $R_\mathrm{specific}$ is the specific gas constant (\si{\joule \per \kelvin \per \kilogram}) for tailpipe ratio of \ce{N_2}, \ce{CO_2}, \ce{O_2}, and \ce{H_2O} that was experimentally measured for a hot TWC. This specific gas constant is used as it is easily determined and the remaining gases only marginally contribute to $R_\mathrm{specific}$.

Ultimately, \crefrange{eq:arrhenius}{eq:ofa} give a simple physics-based model of the most significant reactions that occur in the TWC that takes temperature, gas composition, and residence time into account.

Note that we have implicitly assumed that the incoming gas composition is time-invariant, as this significantly reduces the number of required state variables. (Explicitly modeling a time-varying incoming gas concentration can require an additional $3N$ state variables, one for each emission species concentration in each slice.) This implies that the model is suited for quasi-static combustion engine operation, where the engine-out emission species and massflow varies slowly with respect to the residence time in the TWC. Fortunately, as the residence time in the entire TWC is fairly short (on the order of 0.05 -- 0.1\,\si{\second}) \cite[p. 64]{pannag_coldstart_2009} we hypothesize that moderately-varying dynamic operation with transitions on the order of 0.5 -- 1\,\si{\second} will show accuracy similar to that of constant engine operation.

By \cref{eq:conc-exp}, we can compute the massflow emitted from cell $n,m$ as
\begin{align}
	\dot{m}_{n,m}^{s,\mathrm{out}} &= \dot{m}_{n,m}^{s,\mathrm{in}} e^{-k_{n,m}^s t_{r,n}}\,,
\end{align}
and, by the conservation of mass, the converted massflow is trivially
\begin{align}
	\dot{m}_{n,m}^{s,\mathrm{conv}} &= \dot{m}_{n,m}^{s,\mathrm{in}} - \dot{m}_{n,m}^{s,\mathrm{out}}\,.
\end{align}

This lets us model the tailpipe emissions of emission species $s$ as the sum of the outputs from each individual cell in the last axial segment, i.e.
\begin{align}
	\dot{m}_\mathrm{tp}^s &= \sum_{m=1}^M \dot{m}_{m,N}^s\,.
\end{align}
For convenience, we also define the conversion efficiency of the entire TWC for a given emission species as
\begin{align}
	\eta^s  &= 1 - \frac{\dot{m}_\mathrm{tp}^s}{\dot{m}^s_\mathrm{exh}}\,,
\end{align}
which we can view as the proportion of emissions converted in the TWC.

We generate an estimate of the exothermic reaction power generated by the above reactions by computing the (temperature-dependent) heat of reaction for each mole of reactant species as
\begin{subequations}
	\begin{align}
	\mathrm{d}H_{\ce{CO}} &= H_{0,\ce{CO_{2}}}-H_{0,\ce{CO}}-1/2H_{0,\ce{O_{2}}}\label{eq:dHCO}\\
	\mathrm{d}H_{\ce{H_{2}}} &= H_{0, \ce{H_{2}O}} - H_{0,\ce{H_{2}}} - 1/2H_{0,\ce{O_{2}}}\label{eq:dHH2}\\
	\mathrm{d}H_{\ce{NO}} &= 1/2H_{0,\ce{N_{2}}}+H_{0,\ce{CO_{2}}}-H_{0,\ce{NO}}-\mathrm{d}H_{\ce{CO}}\label{eq:dHNO}\\
	\mathrm{d}H_{\ce{NO}_{2}} &= 1/2H_{0,\ce{N_{2}}}+H_{0,\ce{O_{2}}}-H_{0,\ce{NO_{2}}}\label{eq:dHNO2}\\
	\mathrm{d}H_{\ce{C_{3}H_{6}}} &= 3H_{0,\ce{CO_{2}}}+3H_{0,\ce{H_{2}O}}-H_{0,\ce{C3H6}}-9/2H_{0,\ce{O_{2}}}\label{eq:dHC3H6}\\
	\mathrm{d}H_{\ce{C_{3}H_{8}}} &= 3H_{0,\ce{CO_{2}}}+4H_{0,\ce{H_{2}O}}-H_{0,\ce{C3H8}}-5H_{0,\ce{O_{2}}}\,.\label{eq:dHC3H8}
	\end{align}
\end{subequations}
For brevity, we have not explicitly stated the temperature dependence of the above terms but include their temperature dependence in the numerical model. We use the Shomate equation and reference constants given by the NIST (available at \url{https://webbook.nist.gov}) to compute the numerical values of the above terms. Using the previous concentration ratios for the lumped terms gives the effective reaction power
\begin{subequations}
	\begin{align}
	\mathrm{d}H_{\ce{NO_{x}}} & =(99\mathrm{d}H_{\ce{NO}}+\mathrm{d}H_{\ce{NO_{2}}})/100\label{eq:dHNOx}\\
	\mathrm{d}H_{\ce{THC}} & =(3\mathrm{d}H_{\ce{C_{3}H_{6}}}+\mathrm{d}H_{\ce{C_{3}H_{8}}})/4\,.\label{eq:dHTHC}
	\end{align}
\end{subequations}

By \cref{eq:dHCO,eq:dHH2,eq:dHNOx,eq:dHTHC} the total temperature-dependent heat of reaction generated in each cell is thus
\begin{align}
	P_{n,m} &= \dot{m}_{n,m}^{\ce{CO},\text{conv}}\cdot\mathrm{d}H_{\ce{CO}} +
	\dot{m}_{n,m}^{\ce{H_{2}},\text{conv}} \cdot \mathrm{d}H_{\ce{H_{2}}}\nonumber \\
	&+ \dot{m}_{n,m}^{\ce{NO_{x}},\text{conv}}\cdot\mathrm{d}H_{\ce{NO_{x}}}
	+\dot{m}_{n,m}^{\ce{THC},\text{conv}}\cdot\mathrm{d}H_{\ce{THC}}\,.\label{eq:net_molar_exothermic}
\end{align}

\subsection{Temperature dynamics}
We model the temperature dynamics using a heat balance ODE. Introducing the relative length-weighting matrix
\begin{align}
	\mathbf{W}_L &= \begin{bmatrix}
	L_1/L & 0 & \dots &  0\\
	0 & L_2/L & \dots &  0\\
	\vdots & & \ddots & \vdots\\
	0 & 0 & \dots & L_N/L\\
	\end{bmatrix},
\end{align}
we can define the heat-balance ODEs as
\begin{subequations}
	\begin{align}
	m_\mathrm{TWC} \mathbf{W}_L c_p  \frac{\mathrm{d}T_\mathrm{ctr}}{\mathrm{d}t} &= P_\mathrm{ctr} \\
	m_\mathrm{TWC} \mathbf{W}_L c_p  \frac{\mathrm{d}T_\mathrm{per}}{\mathrm{d}t} &= P_\mathrm{per}\,,
	\end{align}
\end{subequations}
where
\begin{subequations}
	\begin{align}
	P_\mathrm{ctr} &= P_\mathrm{ax} - P_\mathrm{rad} + P_\mathrm{con,ctr} + P_\mathrm{exo,ctr} \label{eq:p_ctr} \\
	P_\mathrm{per} &= P_\mathrm{ax} + P_\mathrm{rad} + P_\mathrm{con,per} + P_\mathrm{exo,per} - P_\mathrm{amb}\,. \label{eq:p_per}
	\end{align}
\end{subequations}
Here, $m_\mathrm{TWC}$ is the mass of the TWC, $c_p$ is its specific heat, $P_\mathrm{ctr}$ and $P_\mathrm{per}$ are $N \times 1$ vectors corresponding to the total power developed in the radial center and periphery respectively, and $\frac{\mathrm{d}T_\mathrm{ctr}}{\mathrm{d}t}$ and $\frac{\mathrm{d}T_\mathrm{per}}{\mathrm{d}t}$ are $N \times 1$ vectors representing the temperature derivative in the radial center and periphery. With these terms, we can construct the total state vector ODE as
\begin{align}
	\frac{\mathrm{d}T}{\mathrm{d}t} &= \begin{bmatrix}
		\frac{\mathrm{d}T_\mathrm{ctr}}{\mathrm{d}t} \\
		N \mathrm{mean}(\mathbf{W}_L (\frac{\mathrm{d}T_\mathrm{per}}{\mathrm{d}t} - \frac{\mathrm{d}T_\mathrm{ctr}}{\mathrm{d}t}))
		\end{bmatrix}\,,\label{eq:diff-eq}
\end{align}
where $\mathrm{mean}([x_1, x_2, \dots, x_n])$ corresponds to the arithmetic mean of the elements in $x$, i.e. $1/n \sum_{i=1}^n x_i$. Note that we can interpret $N \mathrm{mean}(\mathbf{W}_L (\frac{\mathrm{d}T_\mathrm{per}}{\mathrm{d}t} - \frac{\mathrm{d}T_\mathrm{ctr}}{\mathrm{d}t}))$ as corresponding to the average difference between radially central and peripheral powers, weighted by the relative length of each slice. 

The power terms in the right hand side of \cref{eq:p_ctr,eq:p_per} are separated into axial, radial, convection, exothermic, and ambient loss terms respectively, which we will define below. This extends on our previous work \cite{lock2021}, which used lumped-element parameters without an explicit power-balance formulation.

\subsubsection{Axial conduction}
\begin{figure}
	\centering
	\includegraphics[width=0.5\columnwidth]{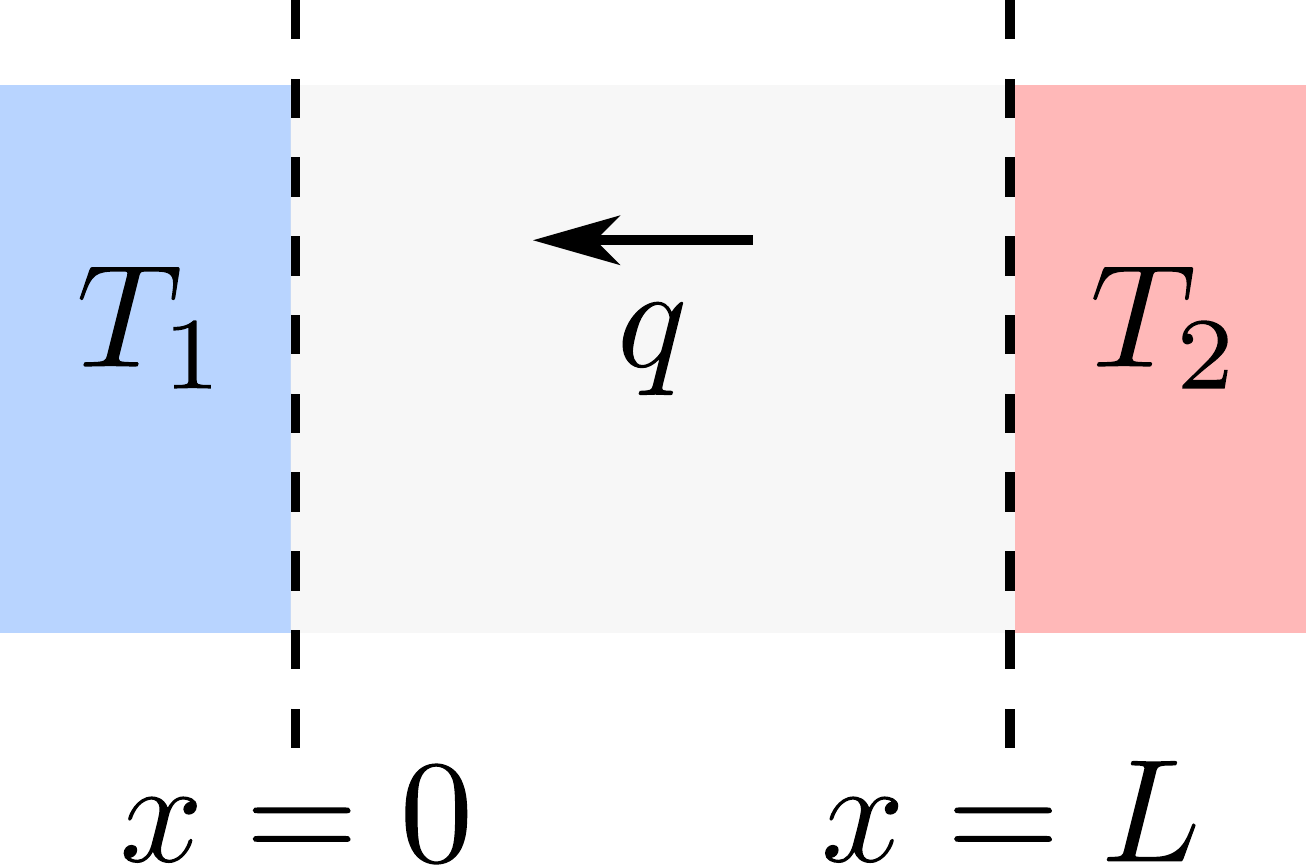}
	\caption{Heat flux between two materials of known temperature. \label{fig:heat-transfer-ref}}
\end{figure}

The axial heat conduction power $P_\mathrm{ax}$ is modeled by Fourier's heat law \cite{logan1998applied}. Using conventional notation, the heat flux between two materials of constant temperature separated by a material of thickness $l$ is in general
\begin{align}
	q &= k \frac{T_2 - T_1}{l},
\end{align}
where $q$ is the heat flux (\si{\watt \per \meter^2}), $T_1$ and $T_2$ are two known temperatures (\si{\kelvin}), and $k$ is the thermal conductivity of the material (\si{\watt \per \meter \per \kelvin}), as illustrated in \cref{fig:heat-transfer-ref}.

In this paper, we extend our previous model \cite{lock2021} by modeling the axial flux between successive axial slices as
\begin{align}
	q_\mathrm{ax} &= k_\mathrm{ax} 
		\begin{bmatrix}
			\dfrac{T_2 - T_1}{L_{1,2}} \\[2ex]
			\dfrac{T_3 - T_2}{L_{2,3}} \\
			\vdots \\
			\dfrac{T_N - T_{N-1}}{L_{N-1,N}}
		\end{bmatrix}
\end{align}
where $k_\mathrm{ax}$ is the axial thermal conductivity and $L_{n,n+1}$ is the distance between the center of axial slice $n$ and $n+1$, as illustrated in \cref{fig:twc-geometry-body}. Note that for $N$ axial slices we thus have $N-1$ axial fluxes between slices.

We model the power associated with each flux term by scaling by the surface area of the solid mass of the TWC, i.e.
$q_\mathrm{ax} (1-\mathrm{OFA}) R^2 \pi$. We can then model the total developed power in each axial slice due to conduction as the difference in incoming and outgoing power fluxes, i.e.
\begin{align}
	P_\mathrm{ax} = \begin{bmatrix}
		0\\
		q_\mathrm{ax} (1-\mathrm{OFA}) R^2 \pi 
	\end{bmatrix} - \begin{bmatrix}
		q_\mathrm{ax} (1-\mathrm{OFA}) R^2 \pi \\
		0
	\end{bmatrix}.
\end{align}

\subsubsection{Radial conduction}
The radial heat conduction is modeled in a manner similar to the axial heat conduction. The radial flux is modeled as
\begin{align}
	q_\mathrm{rad} &= k_\mathrm{rad} \dfrac{\Delta_T}{R/2}\,,
\end{align}
i.e.~a temperature difference of $\Delta_T$ and separation of $R/2$. Approximating the surface area conducting heat as that of a cylinder with half the radius of the TWC and length equal to the TWC's length gives a developed radial conduction power of
\begin{align}
	P_\mathrm{rad} &= q_\mathrm{rad} \pi R L \mathbf{1}\,,
\end{align}
where $\mathbf{1}$ is the ones vector of size $N \times 1$.

\subsubsection{Convection}
The convection heat powers $P_\mathrm{con,ctr}$ and $P_\mathrm{con,per}$ are modeled under the assumption that each cell is sufficiently long and narrow for the gas temperature to reach the cell temperature, i.e.~the gas travels slowly enough to reach thermal equilibrium with the TWC walls. This was considered in our previous model \cite{lock2021}, where we found that five axial slices was a suitable upper limit. As we now also allow for the slice lengths to vary, it is thus reasonable to require that $L_n \leq L/5\,\forall n$.

This gives the convection powers as
\begin{align}
	P_\mathrm{con,ctr} &= \dot{m}_\mathrm{exh} c_{p,\mathrm{exh}} \begin{bmatrix}
		T_\mathrm{exh} - T_1 \\
		T_1 - T_2 \\
		\vdots \\
		T_{N-1} - T_N
	\end{bmatrix} \\
	P_\mathrm{con,per} &= \dot{m}_\mathrm{exh} c_{p,\mathrm{exh}} \begin{bmatrix}
		T_\mathrm{exh} - (T_1 + \Delta_T) \\
		(T_1 + \Delta_T) - (T_2 + \Delta_T) \\
		\vdots \\
		(T_{N-1} + \Delta_T) - (T_N + \Delta_T)
	\end{bmatrix}
\end{align}
where $\dot{m}_\mathrm{exh}$ is the exhaust massflow (\si{\kg \per \second}), $c_{p,\mathrm{exh}}$ is the constant-pressure specific heat of the exhaust gases (\si{\joule \per \kilogram \per \kelvin}), and $T_\mathrm{exh}$ is the temperature of the exhaust gas fed into the TWC (\si{\kelvin}).

\subsubsection{Exothermic power}
The exothermic power terms $P_\mathrm{exo,ctr}$ and $P_\mathrm{exo,per}$ are modeled by weighting the densely-resolved single-channel exothermic power into an effective central and peripheral powers. Here we use a linear weighting scheme as a first approximation, given as
\begin{align}
	P_\mathrm{exotherm, ctr} &= \sum_{m=1}^{M} P_{n,m} (1-\frac{m-1}{M-1}) \\
	P_\mathrm{exotherm, per} &= \sum_{m=1}^{M} P_{n,m} \frac{m-1}{M-1}\,.
\end{align}
Note that the term $\frac{m-1}{M-1}$ varies from 0 to 1 as $m$ varies from 1 to $M$.

\subsubsection{Ambient losses}
The heat losses to the ambient environment are modeled as conductive, with a total flux of

\begin{align}
	q_\mathrm{amb} &= k_\mathrm{amb} \begin{bmatrix}
		(T_1 + \Delta_T - T_\mathrm{amb}) t_\mathrm{amb}^{-1} \\
		(T_2 + \Delta_T - T_\mathrm{amb}) t_\mathrm{amb}^{-1} \\
		\vdots \\
		(T_N + \Delta_T - T_\mathrm{amb}) t_\mathrm{amb}^{-1}\,. \\
	\end{bmatrix}
\end{align}
Here, $k_\mathrm{amb}$ is the effective thermal conductivity of the insulating material and $t_\mathrm{amb}$ is its associated thickness. Modeling the exposed surface area as a cylinder with radius and length equal to the whole TWC (and thus neglecting heat loss through the circular ends of the cylinder) gives the power loss to the ambient environment as
\begin{align}
	P_\mathrm{amb} &= q_\mathrm{amb} 2 \pi R L\,.\label{eq:last-temp-dyn}
\end{align}

\subsection{Two-monolith structure}

\begin{figure}
	\centering
	\includegraphics[width=\columnwidth]{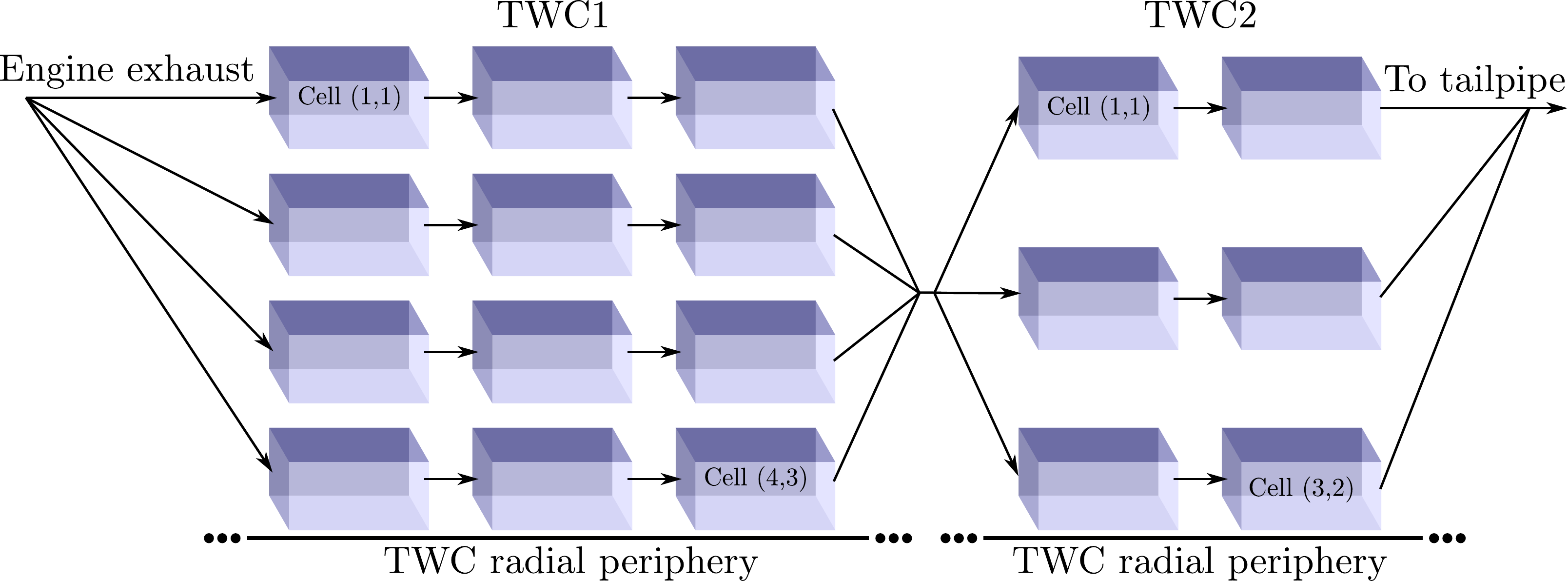}
	\caption{The series-connected TWC is modeled as two independent TWC's with the gas exiting the first being mixed and fed into the second.\label{fig:two-twc}}
\end{figure}

In this paper we extend the model previously presented in \cite{lock2021} by studying a TWC with series-coupled monoliths. More specifically, we consider two physically separated TWC's, where the gas leaving the first is assumed to be completely mixed and then fed into the second, as illustrated in \cref{fig:two-twc}. We model this by assuming two completely independent sets of TWC parameters (as listed in \cref{tab:twc-params}), and let the exhaust gas first travel from the engine through the first TWC. As a first approximation, the gas leaving the first TWC is assumed to be perfectly mixed (both with respect to temperature and emission species concentrations) and then fed through the second TWC, which then finally exits to the tailpipe.

Ideal mixing implies that the temperature of the gas feeding the second TWC is given as
\begin{align}
	T_\mathrm{exh,TWC2} &= \sum_{m=1}^{M_\mathrm{TWC1}} \frac{m^2 - (m-1)^2}{M_\mathrm{TWC1}^2} T_{N,m,\mathrm{TWC1}}\,,
\end{align}
where 
i.e.~the gas leaving the first TWC (TWC1) is combined and scaled by the its relative flow rate. Trivially, we also have that the emission species concentration entering TWC2 is
\begin{align}
	\dot{m}_{1,m,\mathrm{TWC2}}^{s,\mathrm{in}} &= \frac{m^2 - (m-1)^2}{M_\mathrm{TWC2}^2} \dot{m}_\mathrm{tp,TWC1}^s\,,
\end{align}
i.e.~we use the same weighting scheme previously defined in \cref{eq:weighting}.

\section{Experimental setup}\label{sec:exp-setup}

\begin{table}
	\centering{}
	\caption{Engine parameters. \label{tab:exp-engine-setup}}
	\begin{tabularx}{\columnwidth}{l|l}
		\hline 
		Engine type & VEA Gen I, VEP4 MP\\
		\hline 
		Number of cylinders & Four, in-line\\
		\hline 
		Displaced volume & 1969\,cc\\
		\hline 
		Bore/Stroke & 82\,mm/93.2\,mm\\
		\hline 
		Compression ratio & 10.8:1\\
		\hline 
		Valve train & DOHC, 16 valves\\
		\hline 
		Intake camshaft & Variable 0-48\textdegree CA advance\\
		\hline 
		Exhaust camshaft & Variable 0-30\textdegree CA retard\\
		\hline 
		Ignition system & DCI, standard J-gap spark plugs\\
		\hline 
		Fuel system/Injection pressure & DI/200 bar\\
		\hline 
		Fuel & Gasoline RON95 E10\\
		\hline 
		Start of injection & 308-340 CAbTDCf\\
		\hline 
		Boosting system & Turbocharger\\
		\hline 
		Rated power/Rated torque & 187\,kW/350\,Nm\\
		\hline 
		Stoichiometric air/fuel ratio & 14.01:1\\
		\hline 
	\end{tabularx}
\end{table}

As in our previous study in \cite{lock2021}, the experimental setup consisted of a production Volvo Cars two liter in-line four-cylinder direct injected spark ignited turbocharged engine rated for 187\,kW and 350\,Nm, as listed in \cref{tab:exp-engine-setup}.  The engine was connected to an electrical dynamometer that regulated the engine speed and measured the generated torque. A prototyping ECU was used to sample and change engine parameters. The TWC was close-coupled to the turbocharger outlet.

\begin{figure}
	\centering
	\includegraphics[width=\columnwidth]{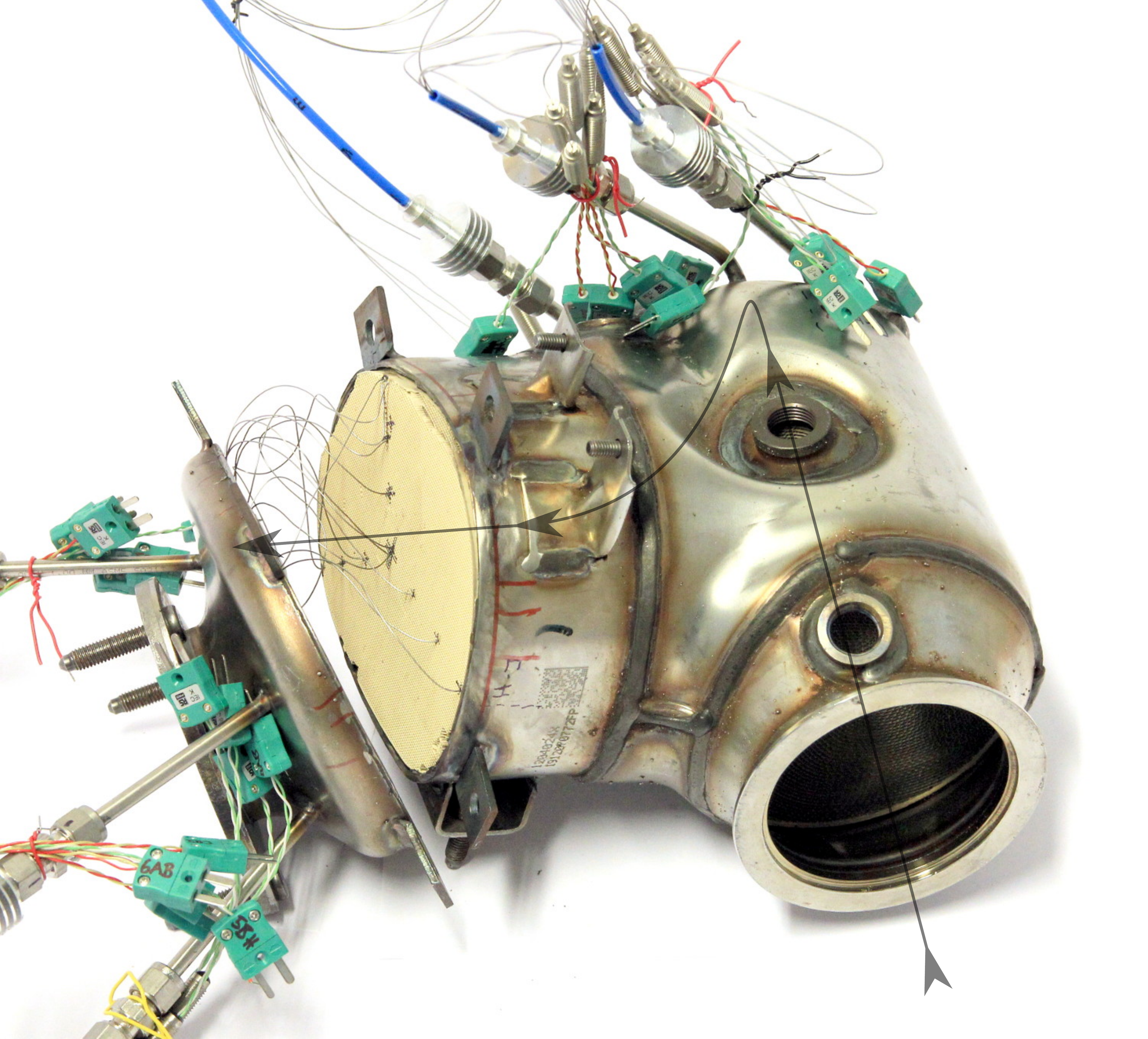}
	\caption{The instrumented TWC, with exposed ceramic (left) and metallic (right) sections. The exhaust gas flow is shown with the highlighted arrow. The cover to the left is open here for illustrative purposes and tightly connected to the main body during operation. \label{fig:exp-setup-instrumented-twc}}
\end{figure}

\begin{figure}
	\centering
	\includegraphics[width=\columnwidth]{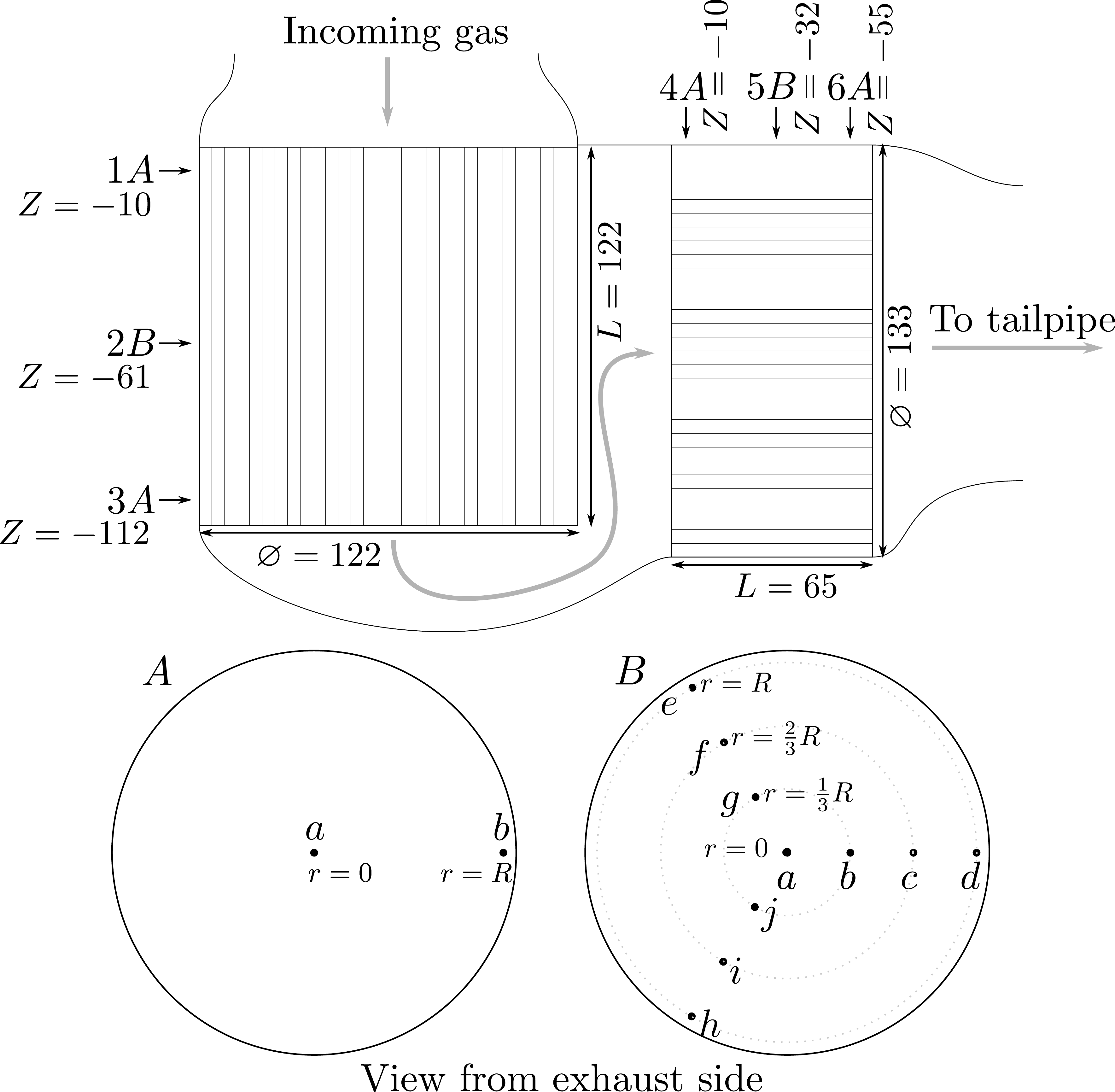}
	\caption{Detailed cross-section of the TWC structure and thermocouple locations. A total of six axial positions were measured (three per monolith) as indicated. Each axial position was sampled in one of two configurations, A and B, as indicated. Each specific sensor is referenced as TS\texttt{XYZ}, where \texttt{X} is the slice number (1-6), \texttt{Y} is the slice type (A/B), and \texttt{Z} is the sensor position (A/B or A-J). \label{fig:exp-setup-sensor-locs}}
\end{figure}

The TWC was instrumented with 28 thermocouples (14 in each monolith) and three exhaust gas sampling locations. The thermocouples, 0.5\,\si{\milli \meter} type-K with a grounded hot junction\footnote{Manufacturer: RS PRO, model number: 847-1110}, were inserted into a TWC channel and held by friction. A close-up of the instrumented TWC is shown in \cref{fig:exp-setup-instrumented-twc}, which also shows the direction of gas flow through the two monoliths. A  more detailed drawing of the TWC construction and the thermocouple locations is shown in \cref{fig:exp-setup-sensor-locs}.

Several thermocouples failed during the experimental campaign. We believe this to be due to the combination of fairly sharply bending the thermocouples in order to reach the required monolith channels and a high level of vibration in an initial TWC mounting fixture. Fortunately, the most critical sensors (in slices 1, 3, 4, and 6) were fully functional and only sensors in slices 2 and 5 were damaged. We excluded data from the damaged sensors in our analysis.

\begin{figure}
	\centering
	\includegraphics[width=\columnwidth]{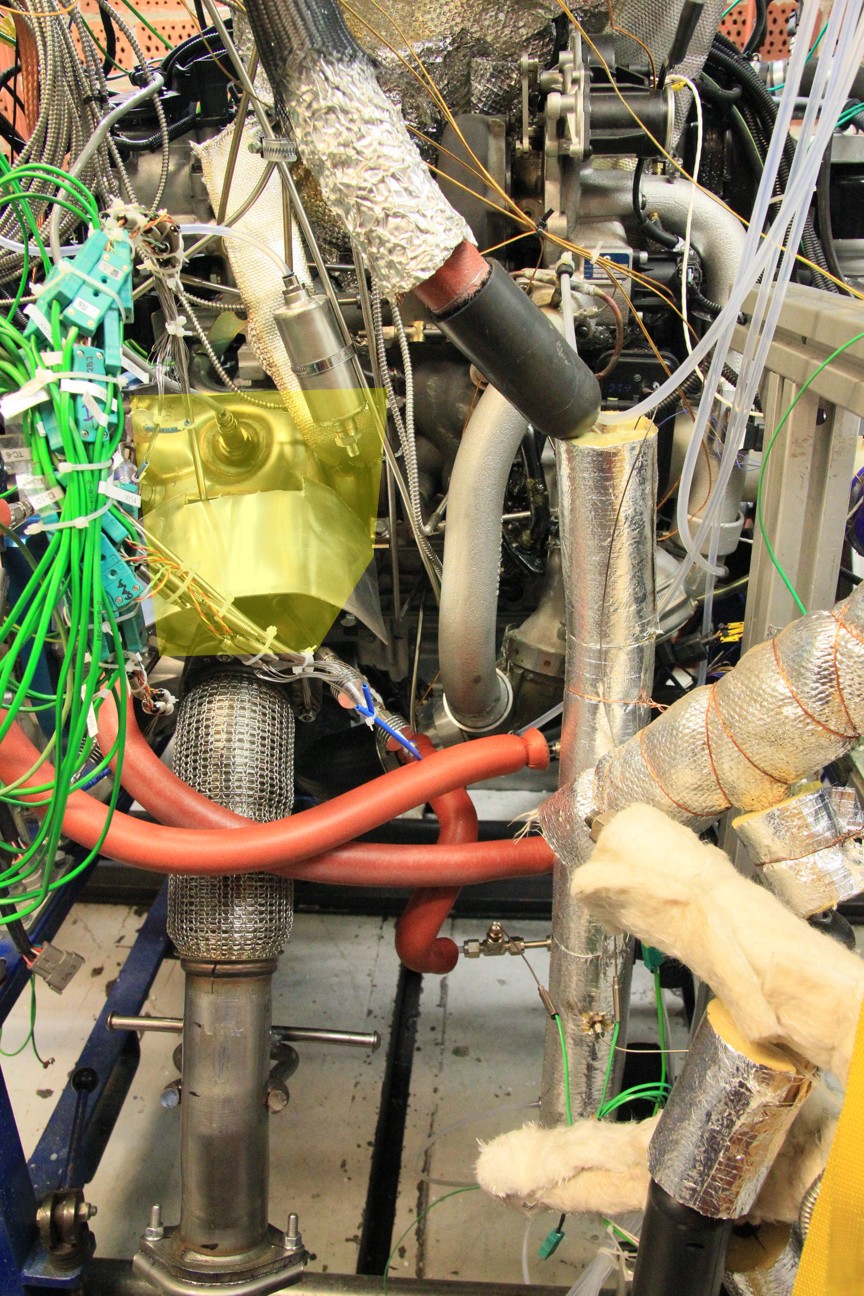}
	\caption{The experimental setup, with the TWC housing and heat shield highlighted. The turbocharger is just visible to the right of the TWC. The exhaust from the TWC is fed down through the visible ducting. \label{fig:exp-setup-photo}}
\end{figure}

\begin{figure}
	\centering
	\includegraphics[width=\columnwidth]{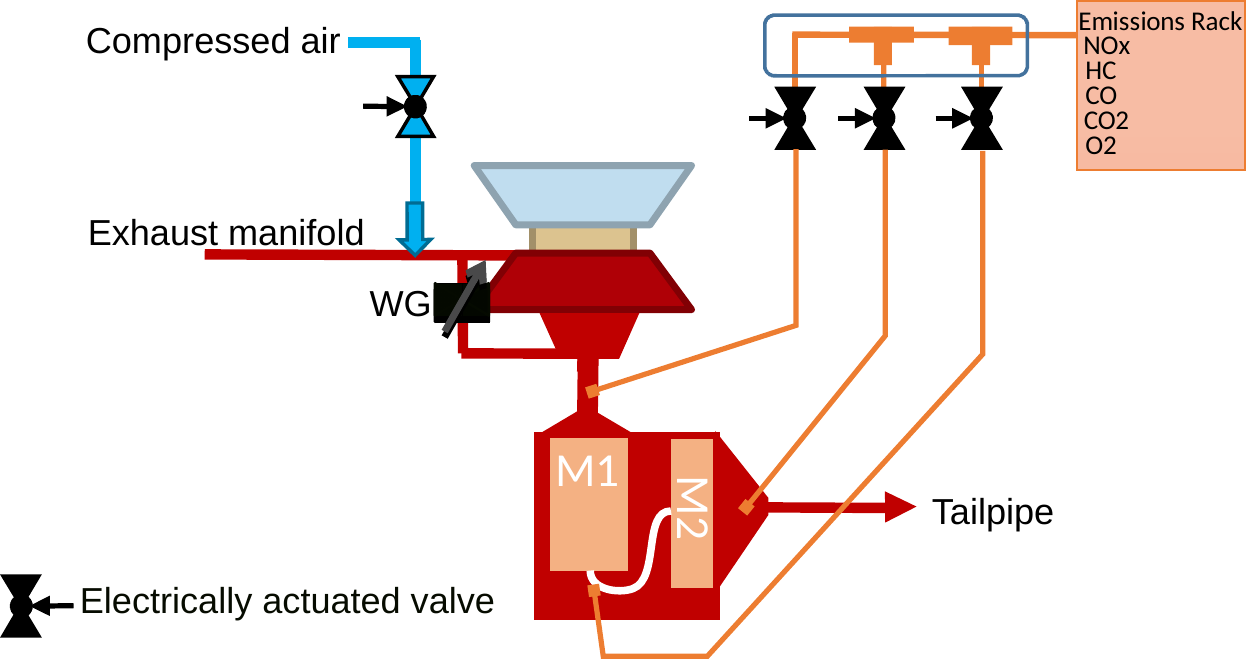}
	\caption{Schematic representation of the experimental setup, with exhaust gas passing through the turbocharger, through monolith 1 (M1), monolith 2 (M2), and finally exiting to the tailpipe. Exhaust gases were sampled at three locations; directly after the turbocharger, between the monoliths, and after monolith 2. \label{fig:exp-setup-schem}}
\end{figure}

An auxiliary air feed was added to the exhaust manifold, which allowed for flushing the entire exhaust subsystem with room-temperature air. By running the engine in fuel-cut mode (i.e.~disabling fuel injection and motoring the engine with the dynamometer) and injecting auxiliary air into the exhaust manifold the exhaust aftertreatment system could be cooled to under 100\textdegree C in approximately 5 minutes. The auxiliary airflow was set to 1000\,\si{\liter \per \minute} STP, which was the maximum flow-rate supported by the mass flow controller. The auxiliary airflow was completely disabled during normal operation (i.e.~when fuel injection was enabled). A photograph of the experimental setup is shown in \cref{fig:exp-setup-photo}, where the engine is visible and the TWC is highlighted. A schematic representation of the experimental set-up and gas flows is shown in \cref{fig:exp-setup-schem}, which also highlights the auxiliary air feed, exhaust gas flow, and gas sampling locations.

The emission sampling points after TWC1 and after TWC2 measured the emissions exiting the radially central channel of each respective TWC. As both TWCs at times displayed a large radial temperature differential, this implies that the average emissions leaving each TWC can be significantly different from the emissions measured at the radial center.

\subsection{Data Acquisition}
Emissions signals from instruments, fuel consumption, and dynamometer readings were sampled with a National Instruments DAQ and an associated LabVIEW program. Engine temperatures, pressures, and the air-fuel ratio was sampled using acquisition units over a CAN ETAS module. All thermocouples were of type K. Fuel massflow was measured with a Coriolis meter. All parameters were sampled at a 10\,\si{\hertz} rate.

Exhaust gases were sampled from three different locations (as illustrated in \cref{fig:exp-setup-schem}). All sampled gases were extracted with a heated hose (180\,\si{\degreeCelsius}), followed by a heated conditioning unit (190\,\si{\degreeCelsius}) with a heated filter and pump. Emissions concentrations were measured with separate instruments. \ce{THC} emissions were measured using a flame ionization detector, \ce{NO_x} using a chemiluminescence analyzer, and \ce{CO} using a non-dispersive infrared detector. The propagation delay and axial dispersion in hoses and instruments was identified by recording the measured engine-out emissions during the transition from fuel-cut operation to normal operation (as will be described in \cref{subsec:measurement-procedure}). With this data we compensated for the propagation delay and applied a first-order high-pass filter to mitigate some of the axial dispersion. This compensation was applied to the remaining two sampling locations, allowing for studying transient emission concentration changes moderately well using an instrument rack primarily intended for steady-state analysis.

As our experimental set-up only allowed for measuring the emissions at one location at any given time it was crucial for the engine-out emissions to be consistent between different runs. Due to this we chose to run the combustion engine in stationary operation, with the goal of maximizing the exhaust gas composition repeatability. We hypothesize that using hardware that measures the emission species at every sample point simultaneously would allow for non-stationary engine operating during cold-start tests.

\subsection{Measurement procedure}\label{subsec:measurement-procedure}
The emission measurement equipment was calibrated before measurements using calibration gases and the engine was heated to its working temperature by operating it at a moderate load until the coolant reached its working temperature. The engine was kept warm during the entire test procedure, implying that the cold-starts studied in this paper refer to the case where the TWC is initially cold while the engine is at operating temperature. Furthermore, the TWC was instrumented with heated lambda sensors, and the engine operated with the conventional closed-loop lambda control scheme during TWC cold-start tests.

\subsubsection{Steady-state analysis}
\begin{table}
	\caption{Tested steady-state engine operating points. The listed speeds, loads, and spark angles (SA) are setpoint values. \label{tab:steady-state-operation-points}}
	\begin{tabularx}{\columnwidth}{c|c|X}
		Speed [RPM] & BMEP [bar] & SA [CAbTDC]\\
		\hline
		\hline
		1000 & 2 & [12, 14, 16, 18, 20, 22, 24]\\
		\hline
		1000 & 5 & [6, 8, 10, 12, 14, 16, 18]\\
		\hline
		1000 & 8 & [-2, 2, 4]\\
		\hline
		1500 & 2 & [12, 14, 16, 18, 20, 21, 22, 24]\\
		\hline
		1500 & 5 & [6, 8, 10, 12, 14, 16, 18]\\
		\hline
		1500 & 8 & [-2, 2, 4, 6, 8, 10]\\
		\hline
		2000 & 2 & [16, 18, 20, 22, 24, 26, 28]\\
		\hline
		2000 & 5 & [8, 10, 12, 14, 16, 18, 20]\\
		\hline
		2000 & 8 & [2, 4, 6, 8, 10, 12]\\
		\hline
		2000 & 10 & [-4, 4]\\
		\hline
		2000 & 12 & [-2, 2, 4]\\
		\hline
		2000 & 14 & [-2, 2]\\
		\hline
		2500 & 2 & [14, 16, 18, 20, 22, 24, 26]\\
		\hline
		2500 & 5 & [12, 14, 16, 18, 20, 22, 24]\\
		\hline
		2500 & 8 & [8, 10, 12, 14, 16]\\
		\hline
		2500 & 13 & [4, 6]\\
		\hline
		3000 & 8 & [8, 10, 12, 14, 16, 18]\\
		\hline
	\end{tabularx}
\end{table}

The goal of this test was to identify the steady-state engine-out emissions and the associated steady-state radial temperature distribution in the TWC. This was performed by statically running the engine at a given speed and BMEP and sweeping the spark angle from the default value and retarding it to the edge of combustion stability. \Cref{tab:steady-state-operation-points} lists the tested speeds, BMEPs, and spark angles tested.

\subsubsection{TWC cold-start characterization}

\begin{table}
	\caption{TWC cold-start load points. SA set to ECU default value.\label{tab:cold-start-op-pts}}
	\centering
	\begin{tabular}{c|c|c|c}
		Index & Speed [rpm] & BMEP [bar] & SA [CAbTDC]\\
		\hline
		\hline
		1 & 1000 & 2 & 24\\
		\hline
		2 & 1000 & 5 & 18\\
		\hline
		3 & 1500 & 5 & 18\\
		\hline
		4 & 1500 & 2 & 24\\
		\hline
		5 & 2000 & 2 & 28\\
		\hline
		6 & 2000 & 5 & 20\\
		\hline
		7 & 3000 & 8 & 18\\
		\hline
		8 & 1000 & 8 & 4\\
		\hline
		9 & 1500 & 8 & 10\\
		\hline
		10 & 2000 & 8 & 12\\
		\hline
	\end{tabular}
\end{table}

The goal of this test was to characterize the cold-start parameters of the two TWC's. The combustion engine was kept at a warm and constant temperature throughout these tests, i.e.~we evaluated the behavior of a cold TWC and warm engine. We performed this experiment by
\begin{itemize}
	\item disabling fuel injection (i.e.~motoring the engine with the dynamometer) and opening the auxiliary air valve until all the TWC thermocouples reported a temperature of under 100\,\si{\degreeCelsius},
	\item first closing the auxiliary air valve, and then immediately enabling ordinary fuel injection until the TWC reached near-equilibrium temperature and emissions.
\end{itemize}
This procedure was repeated for each emission sample point for each of the load points listed in \cref{tab:cold-start-op-pts}. These test points were chosen so that some generated a heating profile that gave a long time to light-off (primarily load points 1--5), while others reached light-off more quickly (load points 6--10). The low-load points gave a longer data stream and reduced the relative error due to axial dispersion in the emission sampling lines. The remaining load points (points 6--10) were more representative of a conventional heating strategy, where light-off is reached more quickly. Furthermore, the load points were characterized by BMEP rather than IMEP due to limitations in the measurement equipment. The engine was kept at a warm and constant temperature to ensure that the exhaust gas composition was not influenced by changes in friction from one load point to the next. It is plausible that regulating for a given IMEP would give an exhaust gas profile that is less sensitive to engine temperature.

\section{Experimental Results}
\subsection{Steady-state} \label{subsec:steady-state}
\begin{table*}
	\centering
	\caption{TWC1 steady-state data for representative load points (speed, BMEP, and spark advance). Figures shown with measured values.\label{tab:steady-state-data}}
\begin{tabular}{c|c|c|c|c|c|c|c|c|c|c|c}
	Speed & BMEP & SA & $\dot{m}_\mathrm{exh}$ & BSFC & \ce{CO} & \ce{NO_x} & \ce{THC} & $TT_{0}$ & $TT_{R/3}$ & $TT_{2R/3}$ & $TT_{R}$ \\
	{[RPM]} & [bar] & [CAbTDC] & [g/s] & [g/kWh] & [ppm] & [ppm] & [ppm] & [\si{\degreeCelsius}] & [\si{\degreeCelsius}] & [\si{\degreeCelsius}] & [\si{\degreeCelsius}]\\
	\hline
	\hline
	1010 & 6.9 & 10 & 13 & 257 & 3580 & 2500 & 427 & 568 & 567 & 566 & 537 \\
	\hline                                                                   
	3000 & 7.9 & 18 & 41 & 259 & 7620 & 2550 & 284 & 839 & 839 & 832 & 813 \\
	\hline                                                                   
	1000 & 1.8 & 12 & 5 & 410 & 6090 & 391 & 548 & 620 & 615 & 609 & 569 \\  
	\hline                                                                   
	1000 & 2.0 & 24 & 5 & 384 & 6650 & 1150 & 684 & 498 & 497 & 491 & 442 \\ 
	\hline                                                                   
	1000 & 4.7 & 6 & 9 & 299 & 5950 & 1200 & 459 & 525 & 526 & 519 & 472 \\  
	\hline                                                                   
	1000 & 5.0 & 18 & 9 & 280 & 6550 & 2070 & 533 & 552 & 552 & 545 & 503 \\ 
	\hline                                                                   
	1570 & 4.8 & 6 & 16 & 313 & 7720 & 533 & 345 & 726 & 725 & 719 & 684 \\  
	\hline                                                                   
	1570 & 5.1 & 19 & 15 & 279 & 9070 & 1350 & 459 & 668 & 667 & 661 & 627 \\
	\hline                                                                   
	2000 & 5.0 & 8 & 21 & 316 & 6730 & 279 & 325 & 715 & 716 & 710 & 674 \\  
	\hline                                                                   
	2000 & 5.2 & 21 & 19 & 274 & 9150 & 686 & 442 & 696 & 697 & 693 & 664 \\ 
	\hline                                                                   
	3000 & 8.1 & 8 & 46 & 277 & 7060 & 1610 & 195 & 893 & 893 & 884 & 864 \\ 
	\hline                                                                   
	3010 & 8.2 & 19 & 43 & 258 & 7560 & 2590 & 328 & 803 & 805 & 797 & 773 \\
	\hline
\end{tabular} 
\end{table*}

\begin{figure}
	\centering
	\begin{subfigure}[b]{\columnwidth}
		\includegraphics[width=\columnwidth]{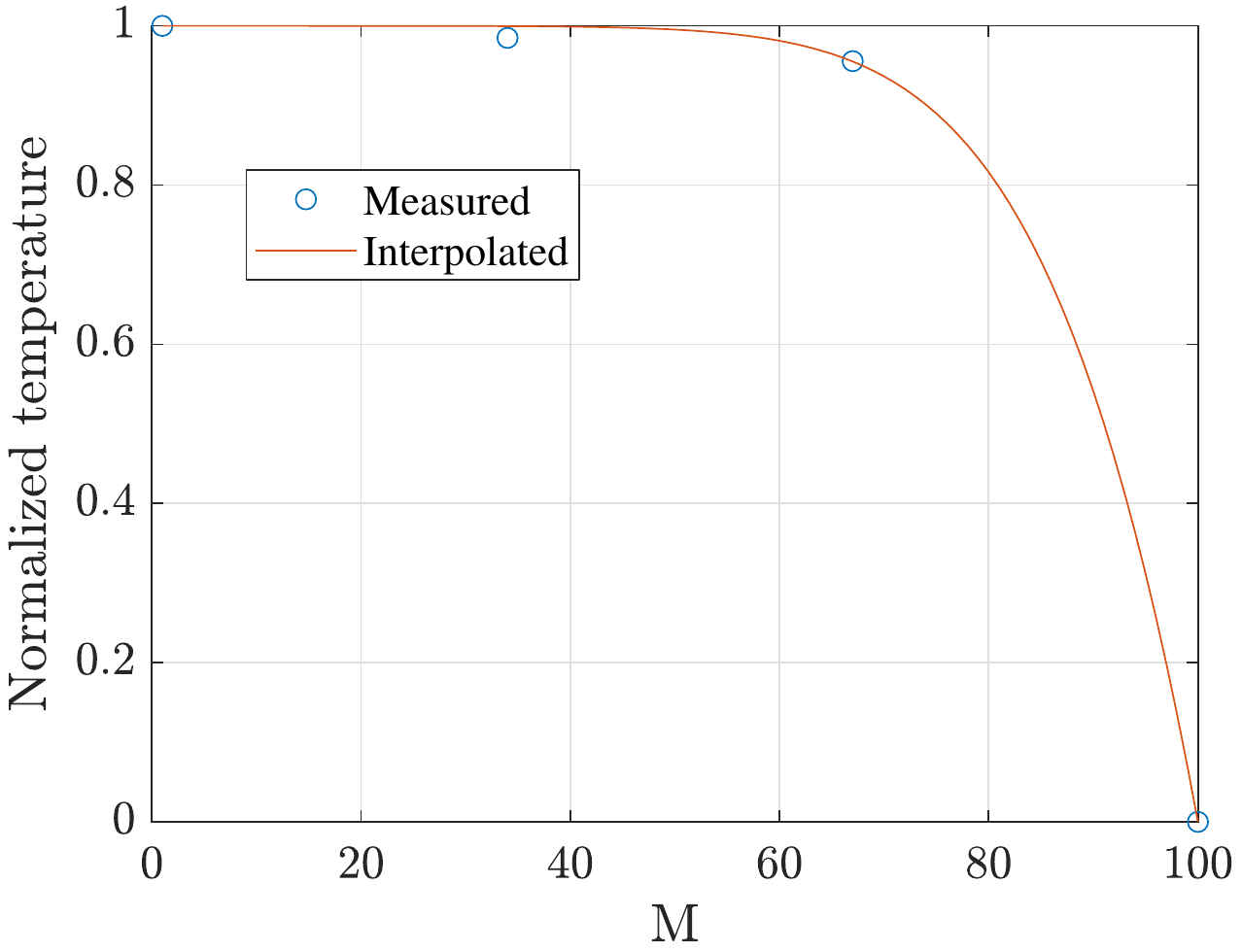}
		\caption{1000\,RPM, 8\,bar, SA 10 CAbTDC}
	\end{subfigure}
	\begin{subfigure}[b]{\columnwidth}
		\includegraphics[width=\columnwidth]{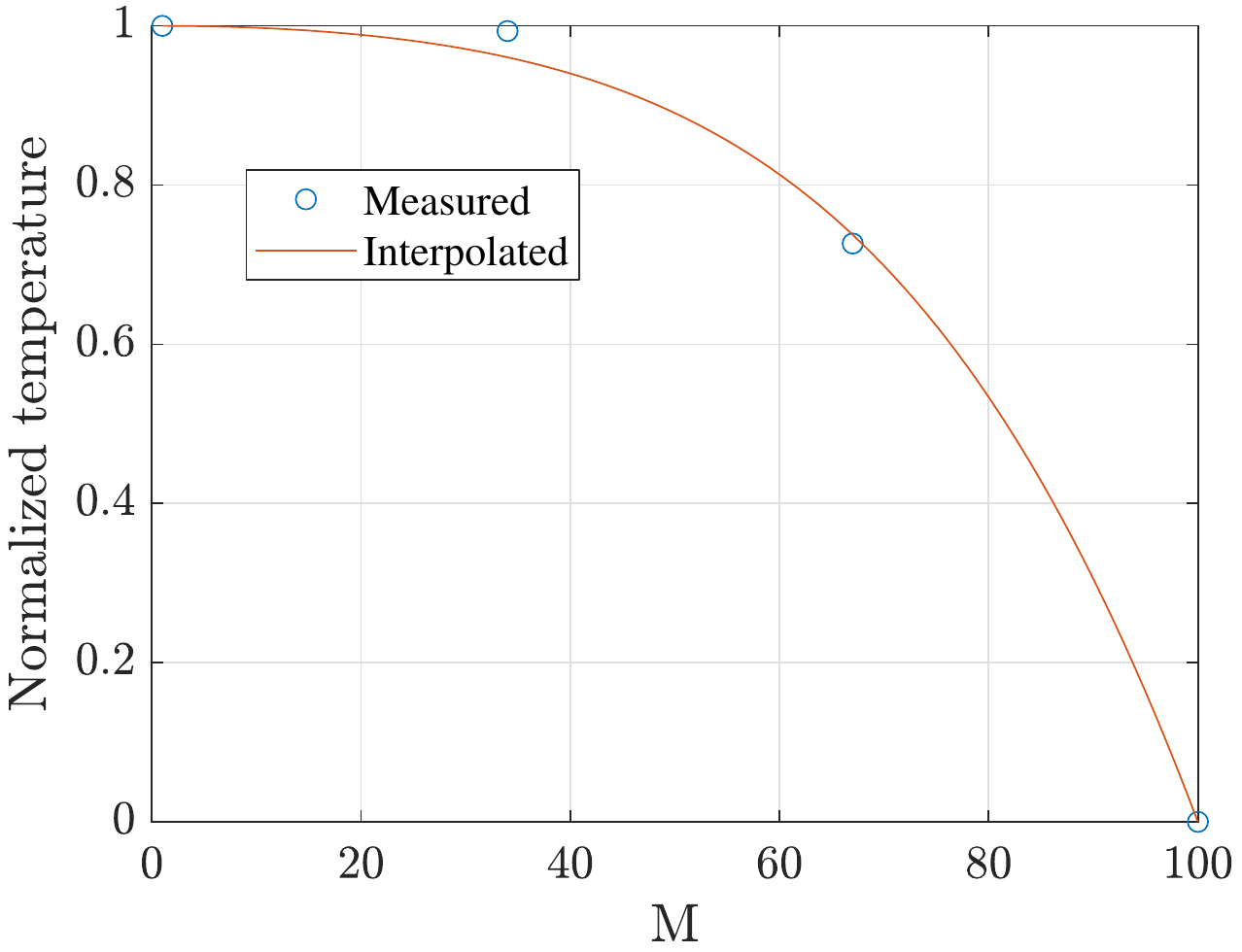}
		\caption{3000\,RPM, 8\,bar, SA 18 CAbTDC}
	\end{subfigure}
	\caption{Measured normalized radial temperature distribution and least-squares interpolation for $M=100$ at two representative load points.\label{fig:rad-profiles}}
\end{figure}

The steady-state experimental results were used to generate a table of the mean equilibrium TWC temperatures, engine-out emissions, exhaust massflow, and engine BSFC for each of the load points listed in \cref{tab:steady-state-operation-points}. Representative data is shown in \cref{tab:steady-state-data} for TWC1. The parameters $TT_0$, $TT_{R/3}$, $TT_{2R/3}$, and $TT_R$ correspond to the mean thermocouple temperature for thermocouples in slice 2 at radius $r=0$, $r=R/3$, $r=2R/3$, $r=R$ respectively (i.e. $TT_{R/3}$ is the mean of TS2BB, TS2BG, TS2BJ). Importantly, we also also use this experimental data to determine the radial interpolation profile outlined in \cref{subsec:rad-temp-interp}, i.e.~for each load point we generate an associated interpolated radial temperature profile. More specifically, for each load point we assume an interpolation function of form
\begin{subequations}
	\begin{align}
	f_\mathrm{interp} &= \hat{T}(t',r) \\
	\shortintertext{where}
	t' &= \underset{t}{\mathrm{argmin}} \sum_{n=0}^3 |\hat{T}(t,nR/3) - TT_{nR/3}|\,,
	\end{align}
\end{subequations}
i.e.~we let $f_\mathrm{interp}$ be the optimal solution in the one-norm sense that minimizes the deviation between the measured temperatures and the solution to the heat equation over all time $t'$. The one-norm is consistently used in this paper in an effort to reduce the effect of outliers. Measured temperatures and the associated interpolation function is shown \cref{fig:rad-profiles} for two representative load points.

\begin{figure*}
	\centering
	\includegraphics[width=0.9\textwidth]{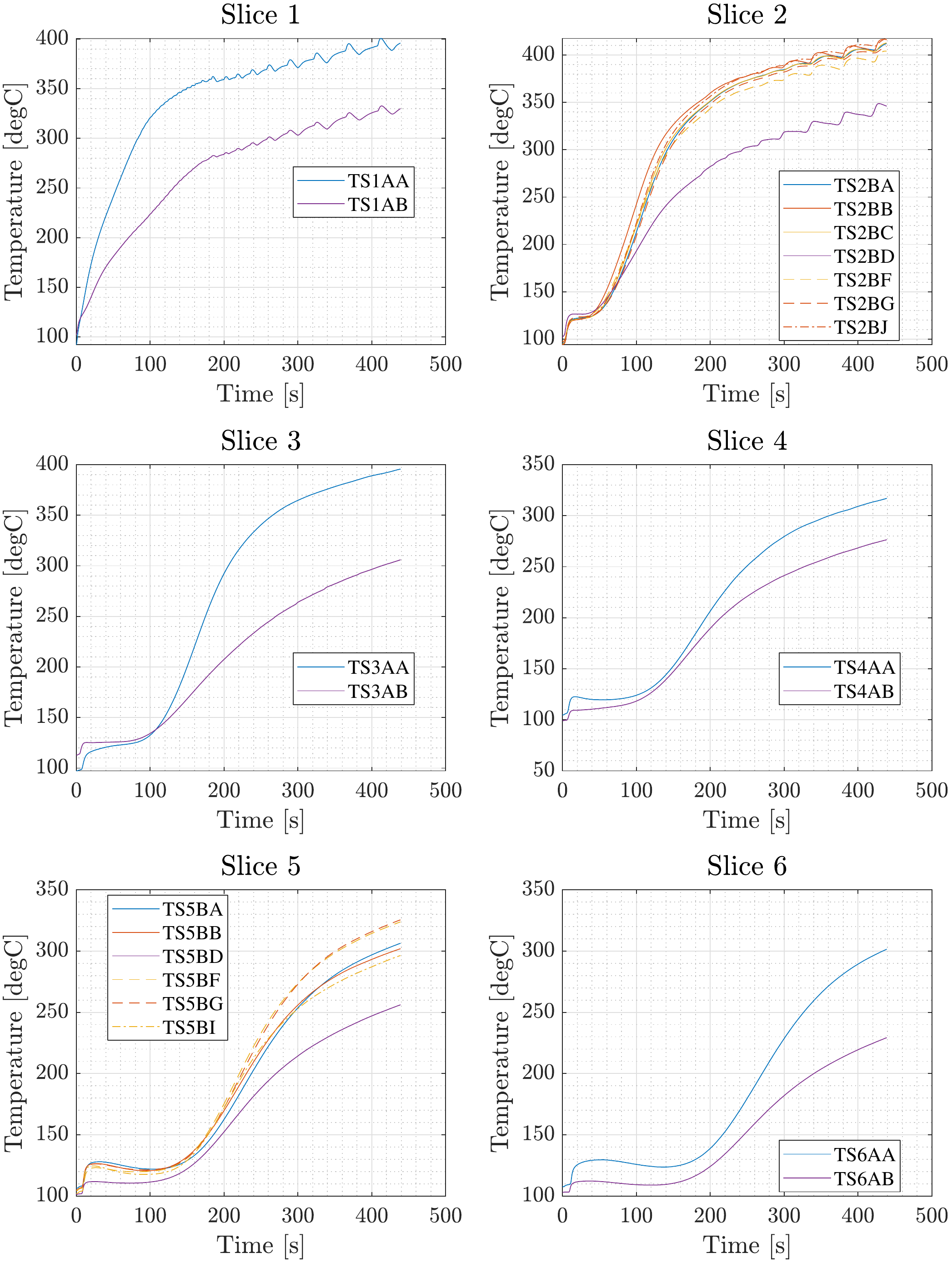}
	\caption{Representative cold-start temperature evolution. Here shown for load point 1 in \cref{tab:cold-start-op-pts} (1000\,RPM, 2\,bar BMEP). Thermocouple locations listed in \cref{fig:exp-setup-sensor-locs}. Damaged thermocouples (TS2BE, TS2BH, TS2BI, TS5BC, TS5BE, TS5BH, TS5BJ) excluded from plots. \label{fig:temp-ev-res}}
\end{figure*}

\subsection{Cold-start} \label{subsec:cold-start}
\begin{figure}
	\centering
	\includegraphics[width=\columnwidth]{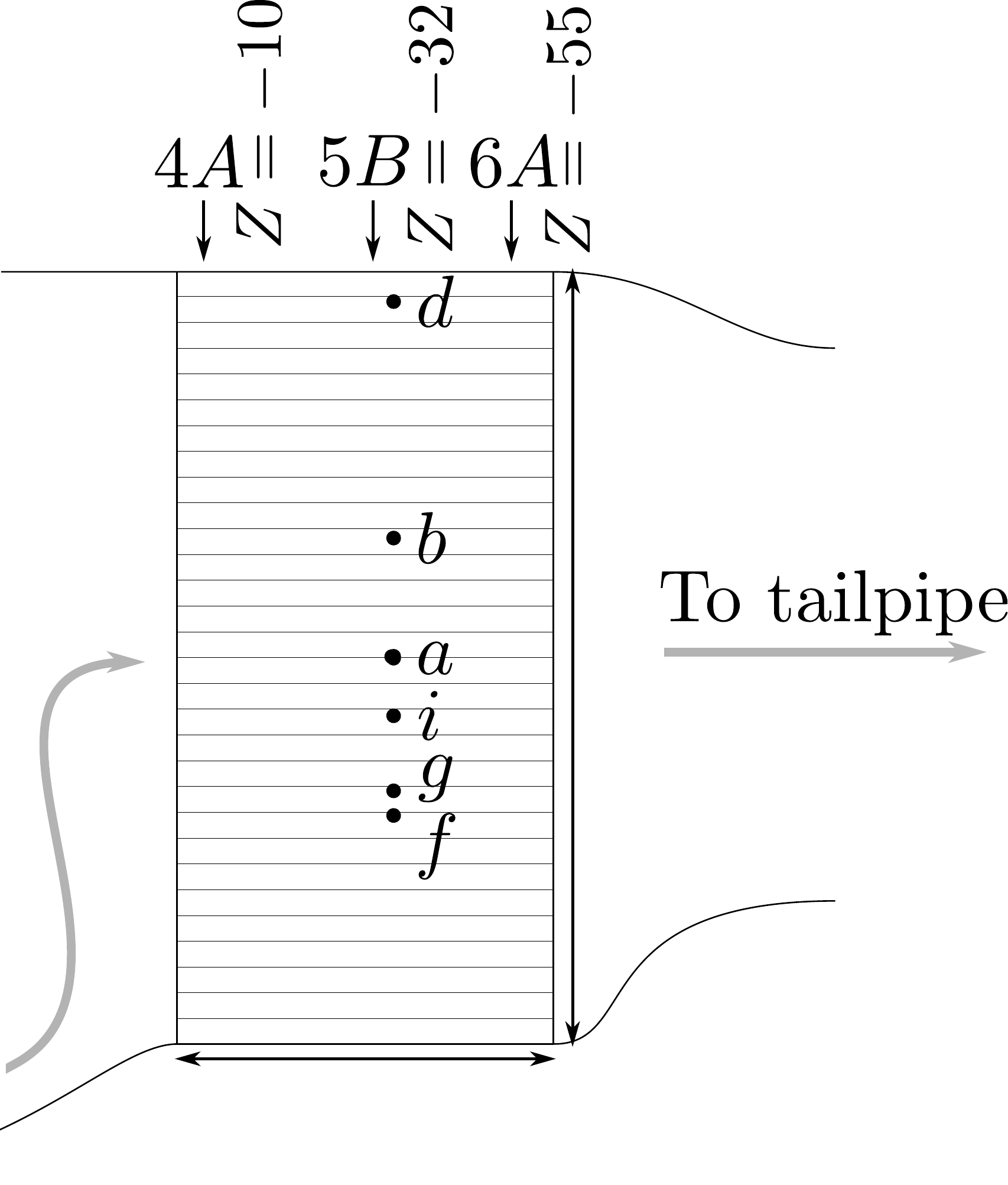}
	\caption{Projection of TS5Bx thermocouples onto the plane. Approximate gas flow shown by gray arrows. The thermocouple closest to the bottom (TS5BF) reaches the highest temperature in \cref{fig:fully-resolved-temps}, while the thermocouple closest to the top (TS5BD) reaches the lowest.\label{fig:sensor-locations-twc2}}
\end{figure}

\Cref{fig:temp-ev-res} shows a representative cold-start temperature evolution, here for 1000\,RPM and 2\,bar BMEP.  We show this specific load point as it gives the longest system dynamics. We can draw several useful conclusions from this test;
\begin{itemize}
	\item The radial temperature distribution is significant throughout both TWC's, with a temperature difference between the radial center and periphery of up to 100\,\si{\degreeCelsius} near light-off.
	\item The first TWC shows no major azimuth temperature variation.
	\item The second TWC does show variations along the azimuth, with the hottest regions nearer the bottom section of the cross-section (see \cref{fig:sensor-locations-twc2}). We hypothesize that this is due to increased massflow near the lower sections, as the sharp bend in the TWC housing causes an uneven pressure distribution across the inlet to the second TWC.
	\item TWC2 shows less pronounced axial temperature variations when compared to TWC1. This could plausibly be due to the length of TWC2, which is only half of TWC1.
\end{itemize}
These results are consistent for the other load points, which display similar results.

Based on these results we have chosen to model TWC1 as consisting of three axial slices, while TWC2 is modeled with a single axial slice. This gives a total of 3+1 state variables for the TWC1 and 1+1 state variables for TWC2, i.e.~a total of six state variables. Though additional slices would have the benefit of improving accuracy, the dynamic-programming based optimal control method we use in \cref{sec:opt-ctrl} is suited for no more than 4--6 state variables. We have chosen to allocate more slices to the first TWC as it displays the most significant axial temperature variations.

\subsection{Model tuning}
The primary goal of the experimental work is to generate measurement data that is used to tune the TWC model. The tuning process was divided into two distinct sections where we first tuned the reaction rate parameters, and afterwards tuned the temperature dynamics parameters. The problem was divided into two sections to reduce the number of degrees of freedom in each optimization step. Furthermore, as the temperature dynamics depends on the exothermic power, it is prudent to first determine the reaction rate parameters.

Of the 10 cold-start simulations listed in \cref{tab:cold-start-op-pts}, we designated operating points $i_\mathrm{train} = [1,3,5,7,9]$ as a training set, and $i_\mathrm{valid} = [2,4,6,8,10]$ as a validation set. All model tuning was done solely using $i_\mathrm{train}$, allowing us to later study the model's accuracy by studying results of applying the model to the validation set's operating points.

\subsection{Reaction rate parameters}\label{subsec:reaction-rate-params}
Here, we consider tuning the per-species reaction-rate parameters $A^s$ and $E_a^s$, giving a total of six parameters to tune per TWC. With the experimental setup as described in \cref{sec:exp-setup}, the gas composition entering and leaving TWC1's radially central channel is well-measured. However, the gas composition entering TWC2 is not as well characterized, as the gas composition leaving TWC1 is inhomogeneous (due to the large radial temperature gradient) and partially mixed before entering TWC2. Due to this, we have chosen to first tune the reaction rate parameters for TWC1, and afterwards make use of the identical precious metal composition of TWC1 and TWC2 (which differ only in their washcoat thickness and loading). This allows us to estimate TWC1's reaction rate parameters using experimental data and then compute the equivalent parameters for TWC2.
 
With respect to TWC1, for a given set of reaction-rate parameters we used the (measured) temperature evolution of each axial slice to simulate the outgoing emission concentration. More specifically, we let the measured state evolution for each operating point be
\begin{subequations}
	\begin{align}
	T_\mathrm{meas, TWC1} &= [T_1, T_2, T_3, \Delta_T]\\
\shortintertext{where}
	T_1 &= \mathrm{TS1AA}\\
	T_2 &= \mathrm{TS2BA}\\
	T_3 &= \mathrm{TS3AA}\\
	\Delta_T &= \mathrm{mean}([\mathrm{TS1AB} - \mathrm{TS1AA}, \mathrm{TS2BD} - \mathrm{TS2BA},\nonumber \\ 
	& \mathrm{TS3AB} - \mathrm{TS3AA}])\,,
	\end{align}\label{eq:meas-state-def-twc1}
\end{subequations}
i.e.~the measured radially central temperatures and the mean difference between the radial center and periphery respectively. This gives a state vector time-evolution $T(k)_\mathrm{meas, TWC1}$, where $k=0,1,2,\dots$ indicates the time-sample of the state vector, sampled at a rate of 10\,\si{\hertz}. In this and later stages we simulated 100 radial channels. We chose to simulate a relatively large number of channels as the kinetics submodel was implemented in a semi-parallel manner that did not require significantly longer to evaluate than for instance 10 channels\footnote{This was in part due to the nature of our model implementation in MATLAB. As the chemical kinetics submodel consists of only basic arithmetic operations we could implement the model in a manner that effectively makes use of MATLAB's numerically efficient matrix operations. This gave a model implementation whose execution time was dominated by the the number of function calls, i.e.~the sample rate and simulation time, rather than the number of radial channels.}.

We started the tuning process by first determining an initial guess, where we selected $E_a^s$ as determined by \cite{ramanathan_kinetic_2011} and $A^s$ was selected to give a light-off temperature near the experimentally measured behavior. We then used MATLAB's \texttt{patternsearch} utility (a zero'th order gradient-descent optimization method) to minimize the 1-norm penalty
\begin{multline}
 J^* = \min_{E_{a,\mathrm{TWC1}}^s, A_\mathrm{TWC1}^s} \sum_{i \in i_\mathrm{train}} \sum_{k=1}^K|\dot{m}(k)_{3,1,\mathrm{TWC1},i}^{s,\mathrm{out}} \\
  - \dot{m}(k)_\mathrm{meas,\mathrm{TWC1 out},i}^s|
\end{multline}
using the \texttt{MADSPositiveBasis2N} polling method and with the \texttt{UseCompletePoll} flag set. Here, $\dot{m}(k)_{3,1,\mathrm{TWC1},i}^{s,\mathrm{out}}$ is the simulated concentration of emission species $s$ emitted by the third (i.e.~last) axial slice at the radial center at sample $k$ and operating point $i$, $\dot{m}(k)_\mathrm{meas,\mathrm{TWC1 out},i}^s$ is the measured emissions at the same time instance and operating point, and $i$ iterates over the operating points in the training dataset.

With the kinetics parameters for TWC1 determined, we estimate the parameters for TWC2 by assuming an identical activation energy and with the pre-exponential term scaled by the amount of catalytically active material. This gives the estimates
\begin{subequations}
	\begin{align}
		E_{a,\mathrm{TWC2}}^s &= E_{a,\mathrm{TWC1}}^s \\
		A_{\mathrm{TWC2}}^s &= \frac{t_\mathrm{wash,TWC2} \cdot w_\mathrm{TWC2}}{t_\mathrm{wash,TWC1} \cdot w_\mathrm{TWC1}} A_{\mathrm{TWC1}}^s\,,
	\end{align}
\end{subequations}
where $t_\mathrm{wash}$ and $w$ are the washcoat thickness and loading respectively for each TWC.

\begin{table}
	\centering
	\caption{Tuned reaction rate parameters.\label{tab:tuned-reaction-params}}
	\begin{tabular}{c|c|c}
		Parameter & Value & Unit\\
		\hline
		\hline
		$E_{a,\mathrm{TWC1,2}}^{\ce{CO}}$ & $84.0\cdot10^3$ & \si{\joule \per \mole} \\
		\hline
		$E_{a,\mathrm{TWC1,2}}^{\ce{NO_x}}$ & $82.1\cdot10^3$ & \si{\joule \per \mole} \\
		\hline
		$E_{a,\mathrm{TWC1,2}}^{\ce{THC}}$ & $51.0\cdot10^3$ & \si{\joule \per \mole} \\
		\hline
		$A_{\mathrm{TWC1}}^{\ce{CO}}$ & $59.5\cdot10^9$ & - \\
		\hline
		$A_{\mathrm{TWC1}}^{\ce{NO_x}}$ & $27.6\cdot10^6$ & - \\
		\hline
		$A_{\mathrm{TWC1}}^{\ce{THC}}$ & $11.2\cdot10^9$ & - \\
		\hline
		$A_{\mathrm{TWC2}}^{\ce{CO}}$ & $23.8\cdot10^9$ & - \\
		\hline
		$A_{\mathrm{TWC2}}^{\ce{NO_x}}$ & $16.1\cdot10^6$ & - \\
		\hline
		$A_{\mathrm{TWC2}}^{\ce{THC}}$ & $5.63\cdot10^9$ & - \\
		\hline
	\end{tabular}
\end{table}

\Cref{tab:tuned-reaction-params} lists the identified parameters for both TWCs, which are on the same order of magnitude as literature suggests \cite{ramanathan_kinetic_2011, kang_detailed_2014}. Note that the \texttt{patternsearch} method is similar to gradient-descent methods, and as the problem is non-convex is therefore not guaranteed to return a globally optimal solution.

\begin{figure}
	\centering
	\includegraphics[width=\columnwidth]{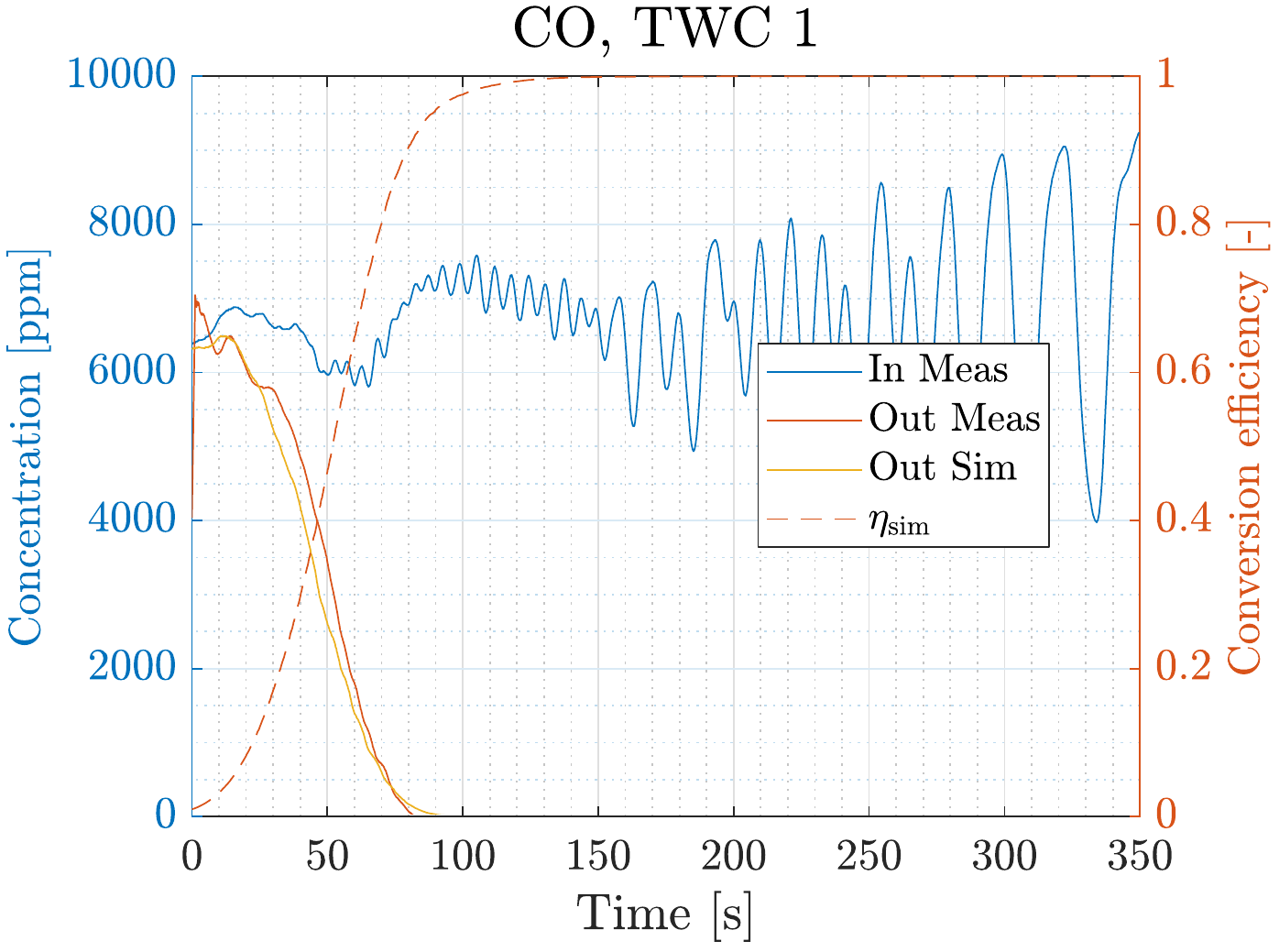}
	\includegraphics[width=\columnwidth]{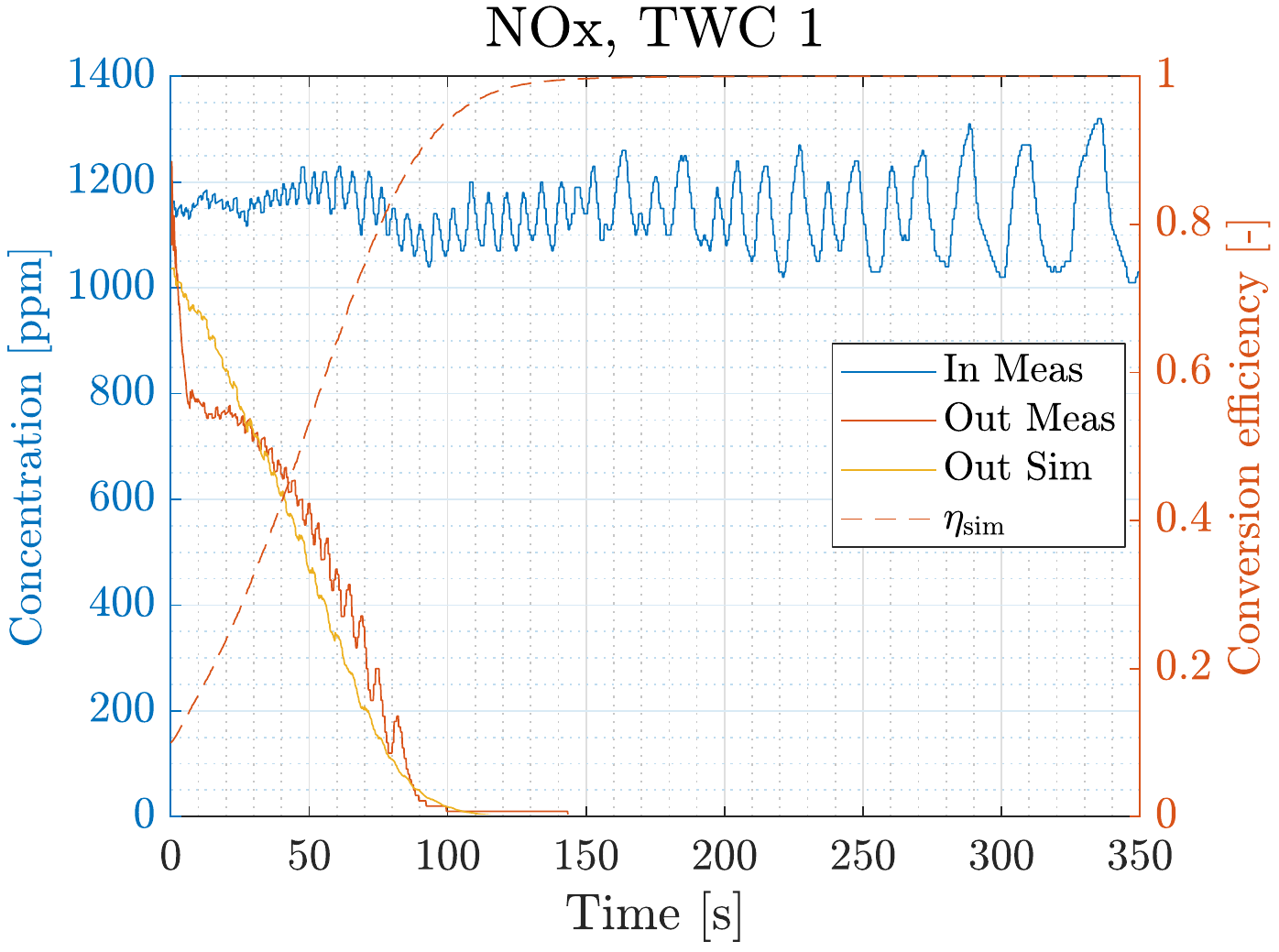}
	\includegraphics[width=\columnwidth]{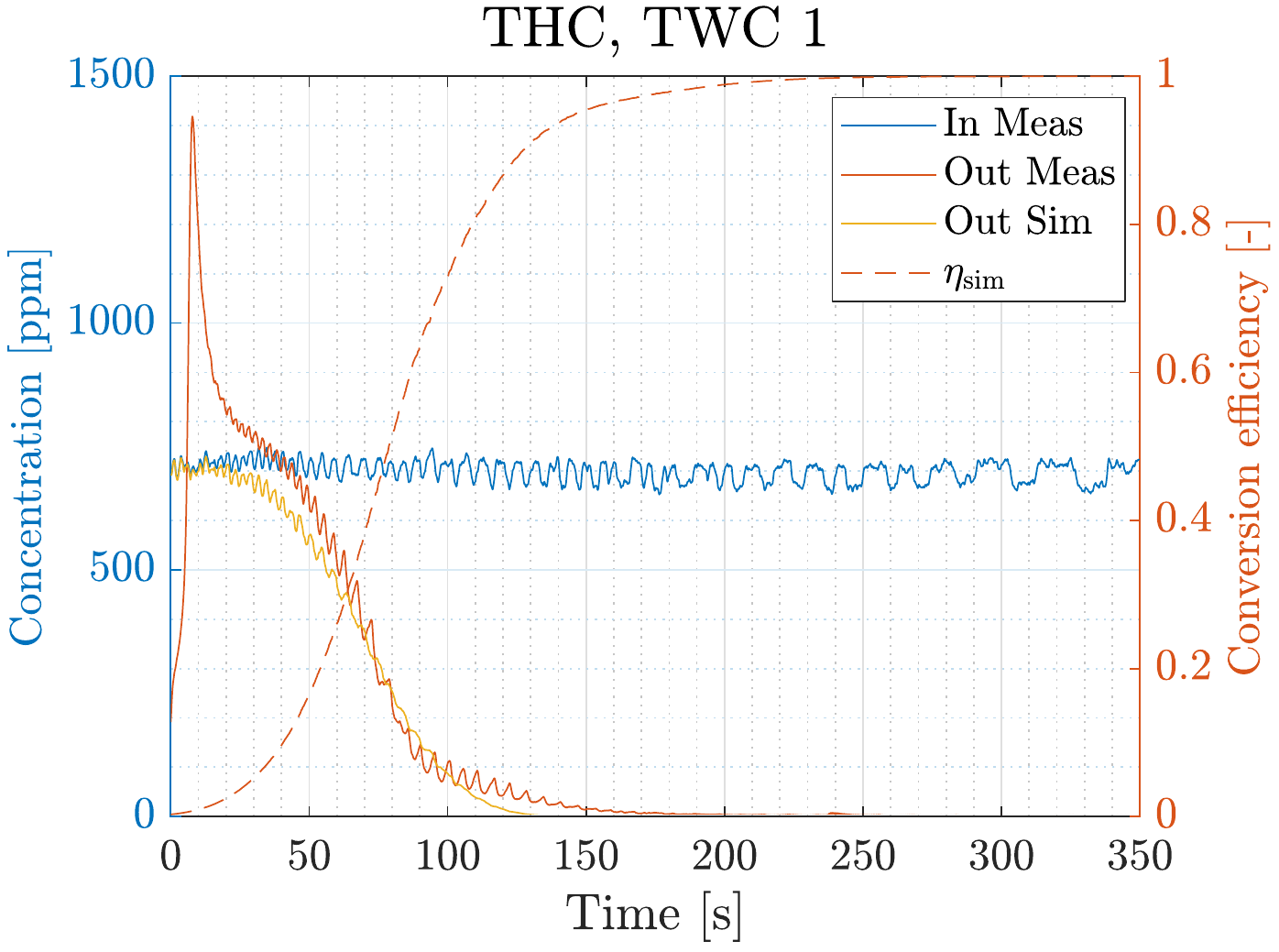}
	\caption{Emissions profiles for 1000\,rpm, 2\,bar BMEP load point. See \cref{fig:temp-ev-res} for the associated measured temperature evolution and \cref{fig:temp-ev-twc1} for the associated state evolution. \label{fig:emissions-sim-vs-modeled}}
\end{figure}

Time-resolved plots of the measured and simulated emissions profiles and the simulated conversion efficiencies are shown in \cref{fig:emissions-sim-vs-modeled} for the lowest-load operating point. The simulated outgoing emissions are shown for the simulated radially central channel, which corresponds to the location of the measured emissions. The most significant deviations are seen at $t\in[0, 30]$ for \ce{NO_x} and at $t \in [0, 50]$ for \ce{THC}. We hypothesize that the former is due to adsorption and the latter due to the poor transient response of the measurement equipment.

\subsection{Temperature dynamics parameters}
With the reaction rate parameters determined we consider the temperature dynamics parameters $k_\mathrm{ax}$, $k_\mathrm{ra}$, $k_\mathrm{amb}$, $c_p$, and $L_n$ separately for each TWC. Note that as TWC2 is modeled with only a single axial slice, there is no modeled axial conduction and $k_\mathrm{ax}$ and $L_n$ therefore have no meaning. We tuned each TWC to the measured data by assigning the initial state to the measured temperature at the start of the cold-start test and then applying a 1-norm penalty to the deviation between the simulated and measured states, i.e.
\begin{align}
	J^* &= \min_{k_\mathrm{ax}, k_\mathrm{ra}, k_\mathrm{amb}, c_p, L_n} \sum_{i \in i_\mathrm{train}} \sum_{k=1}^K |T(k)_{\mathrm{meas},i} - T(k)_{\mathrm{sim},i}|
\end{align}
where $T_\mathrm{sim}$ is the simulated state evolution generated by solving \cref{eq:diff-eq} using an explicit fourth-order Runge-Kutta solver with a fixed time-step of 0.1\,\si{\second}, $T(0)_\mathrm{sim}$ is initialized as $T(0)_\mathrm{sim} = T(0)_\mathrm{meas}$, $T_\mathrm{meas}$ is defined by \cref{eq:meas-state-def-twc1} for TWC1, and for TWC2 defined as
\begin{align}
	T_\mathrm{meas} &= [T_1,\Delta_T]\\
\shortintertext{where}
	T_1 &= \mathrm{mean}([\mathrm{TS4AA}, \mathrm{TS5BA}, \mathrm{TS6AA}])\\
	\Delta_T &= \mathrm{mean}([\mathrm{TS4AB} - \mathrm{TS4AA}, \nonumber \\
	& \mathrm{TS5BD} - \mathrm{TS5BA}, \mathrm{TS6AB} - \mathrm{TS6AA}])\,.
\end{align}
As in the reaction rate parameters, we used the \texttt{patternsearch} method to determine the optimal parameters. We supplied the initial guess for $k_\mathrm{ax}$, $k_\mathrm{ra}$, $k_\mathrm{amb}$, and $c_p$ by setting them to the values specified by the TWC manufacturer, and $L_n$ to geometrically ideal values, where we assume the thermocouples are placed in the center of each slice. Referencing \cref{fig:exp-setup-sensor-locs} gives the initial guess $L_1 = 20 \cdot 10^{-3}$\,\si{\meter}, $L_2 = 102 \cdot 10^{-3}$\,\si{\meter}, and $L_3 = 20 \cdot 10^{-3}$\,\si{\meter}. \Cref{tab:tuned-temp-dyn-params} lists the parameter values found after tuning, as well as the known (i.e.~assigned) fixed model parameters.

\begin{table}
	\centering
	\caption{Tuned and fixed temperature dynamics parameters.\label{tab:tuned-temp-dyn-params}}
	\begin{tabular}{c|c|c}
		Tuned Parameter & Value & Unit\\
		\hline
		\hline
		$k_\mathrm{ax, TWC1}$ & 319 & \si{\watt \per \meter \per \kelvin} \\
		\hline
		$k_\mathrm{ra, TWC1}$ & 46.6 & \si{\watt \per \meter \per \kelvin} \\
		\hline
		$k_\mathrm{amb, TWC1}$ & 0.421 & \si{\watt \per \meter \per \kelvin} \\
		\hline
		$c_{p,\mathrm{TWC1}}$ & 2318 & \si{\joule \per \kelvin \per \kilogram} \\
		\hline
		$L_{1,\mathrm{TWC1}}$ & $32.1 \cdot 10^{-3}$ & \si{\meter} \\
		\hline
		$L_{2,\mathrm{TWC1}}$ & $48.2 \cdot 10^{-3}$ & \si{\meter} \\
		\hline
		$L_{3,\mathrm{TWC1}}$ & $41.8 \cdot 10^{-3}$ & \si{\meter} \\
		\hline
		$k_\mathrm{ra, TWC2}$ & 4.53 & \si{\watt \per \meter \per \kelvin} \\
		\hline
		$k_\mathrm{amb, TWC2}$ & 0.602 & \si{\watt \per \meter \per \kelvin} \\
		\hline
		$c_{p,\mathrm{TWC2}}$ & 2360 & \si{\joule \per \kelvin \per \kilogram} \\
		\hline
		\multicolumn{1}{c}{} & \multicolumn{1}{c}{} & \multicolumn{1}{c}{} \\
		Fixed Parameter & Value & Unit \\
		\hline
		\hline
		$\mathrm{OFA}_\mathrm{TWC1}$ & 0.935 & - \\
		\hline
		$\mathrm{OFA}_\mathrm{TWC2}$ & 0.846 & - \\
		\hline
		$m_\mathrm{TWC1}$ & 0.418 & \si{\kilogram}\\
		\hline
		$m_\mathrm{TWC2}$ & 0.248 & \si{\kilogram} \\
		\hline
		$c_{p,\mathrm{exh}}$ & 1050 & \si{\joule \per \kelvin \per \kilogram} \\
		\hline
		$t_\mathrm{amb, TWC1}$ & $10 \cdot 10^{-3}$ & \si{\meter}\\
		\hline
		$t_\mathrm{amb, TWC2}$ & $10 \cdot 10^{-3}$ & \si{\meter}\\
		\hline		
		$T_\mathrm{amb}$ & 25 & \si{\degreeCelsius} \\
		\hline
	\end{tabular}
\end{table}

\begin{figure}
	\begin{subfigure}{\columnwidth}
		\includegraphics[width=\textwidth]{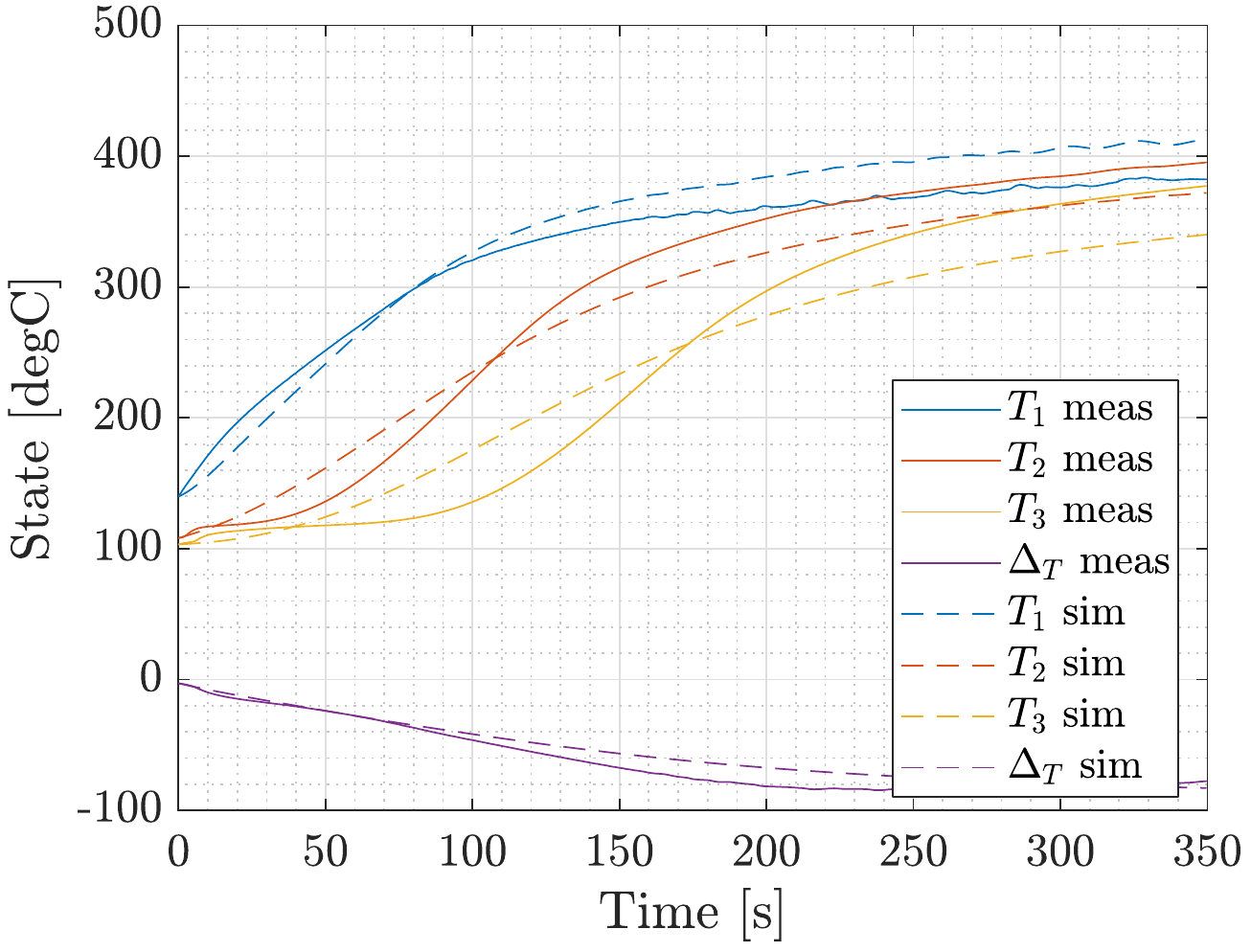}
		\caption{Measured and simulated temperature evolution for TWC1 with 3 axial slices.}\label{fig:temp-ev-twc1}
	\end{subfigure}\\
	\begin{subfigure}{\columnwidth}
		\includegraphics[width=\textwidth]{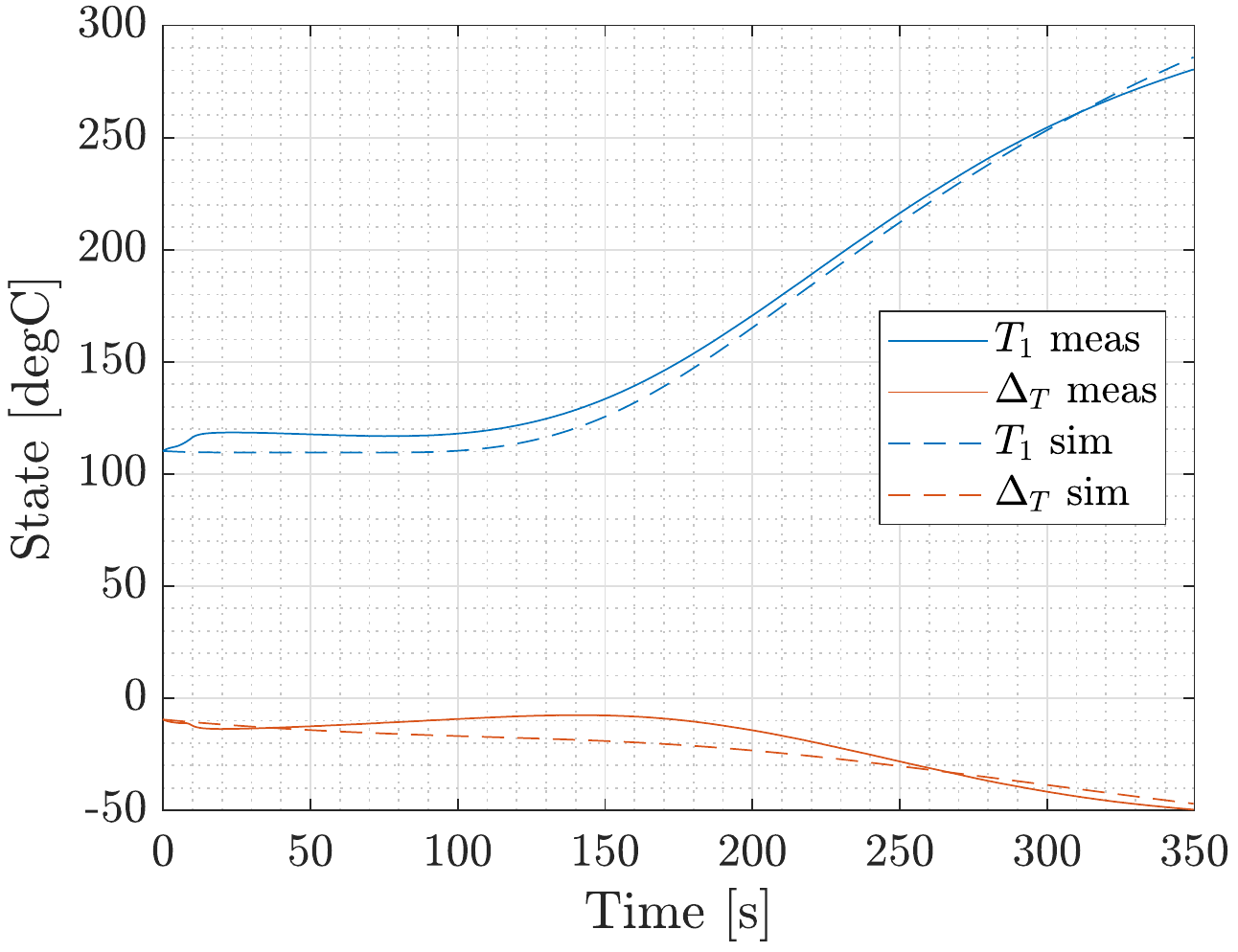}
		\caption{Measured and simulated temperature evolution for TWC2.}\label{fig:temp-ev-twc2}
	\end{subfigure}\\
	\caption{Temperature evolution for 1000\,rpm, 2\,bar BMEP operating point.}\label{fig:temp-ev}
\end{figure}

\begin{figure}
	\begin{subfigure}{\columnwidth}
		\includegraphics[width=\textwidth]{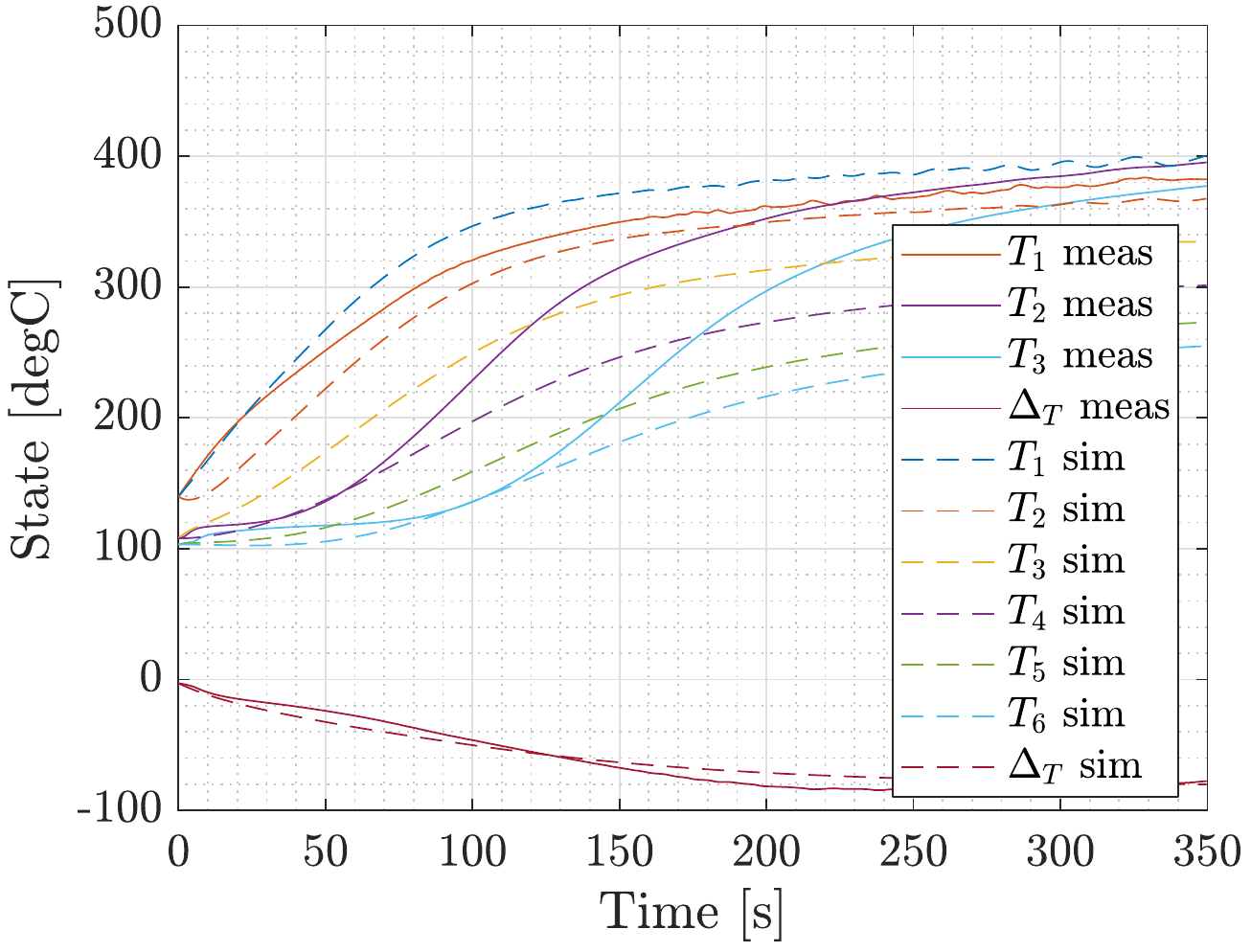}
		\caption{Measured and simulated temperature evolution at the 1000\,rpm, 2\,bar BMEP operating point, TWC1 model modified to use 6 axial slices. Other model parameters unchanged.}\label{fig:temp-ev-twc1-slices}
	\end{subfigure}
	\begin{subfigure}{\columnwidth}
		\includegraphics[width=\textwidth]{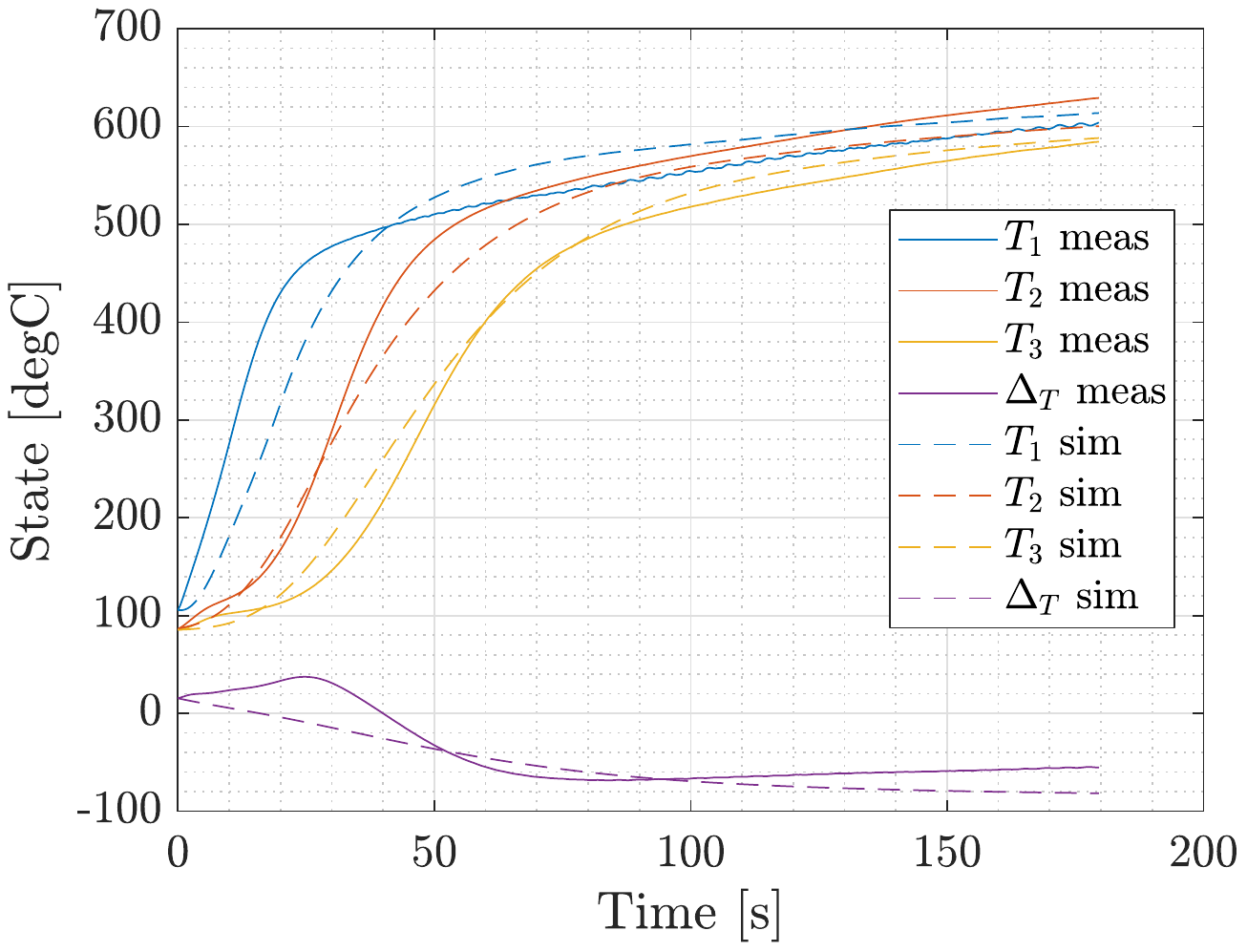}
		\caption{Measured and simulated temperature evolution for TWC1 at the 1500\,rpm, 8\,bar BMEP load point.}\label{fig:temp-ev-twc1-highload}
	\end{subfigure}
	\caption{Increasing the number of resolved slices (\cref{fig:temp-ev-twc1-slices}) and/or heating the TWC more quickly (\cref{fig:temp-ev-twc1-highload}) reduces the modeling error caused by the limited ability to represent a delay.}\label{fig:temp-ev-twc1-delay}
\end{figure}

An illustration of a representative temperature evolution is shown in \cref{fig:temp-ev-twc1} and \cref{fig:temp-ev-twc2}. Though the first slice of TWC1 and TWC2 capture the measured temperature evolution well, the second and third slices of TWC1 do not capture the characteristic delay shown in the measured data. We hypothesize that this is independent of the chosen tuning parameters and an inherent limitation of our modeling assumption of a small number of axial slices. More specifically, as a general discrete-time delay of $n$ samples (trivially) requires storing the values of the $n$ samples, our shown model can only represent a true delay of three samples, i.e. an insignificant~0.3\,seconds, before the last axial segment starts displaying a positive temperature derivative. This can be alleviated somewhat either by increasing the number of axial segments or by considering cold-starts with a less prominent delay, as is shown in \cref{fig:temp-ev-twc1-delay}.

\subsection{Cumulative emissions accuracy}
With the model tuned, we will now turn to quantitatively evaluating the TWC model's accuracy. Here, we consider the relative difference between the cumulative measured and simulated emissions (i.e.~cold-start ``bag emissions'') for each emission species and TWC. Using the notation where $\Delta_i^s$ corresponds to the $i$'th TWC for emission species $s$ gives
\begin{subequations}
	\begin{align}
		\Delta_1^s &= \frac{\sum_{k=0}^K \dot{m}(k)^{s,\mathrm{out}}_{3,1,\mathrm{TWC1}}}{\sum_{k=0}^K\dot{m}(k)^s_\mathrm{meas,TWC1 out}} - 1 \\
		\Delta_2^s &= \frac{\sum_{k=0}^K \dot{m}(k)^{s,\mathrm{out}}_{1,1,\mathrm{TWC2}}}{\sum_{k=0}^K\dot{m}(k)^s_\mathrm{meas,TWC2 out}} - 1\,.
	\end{align}
\end{subequations}

\begin{figure}
	\centering
	\includegraphics[width=\columnwidth]{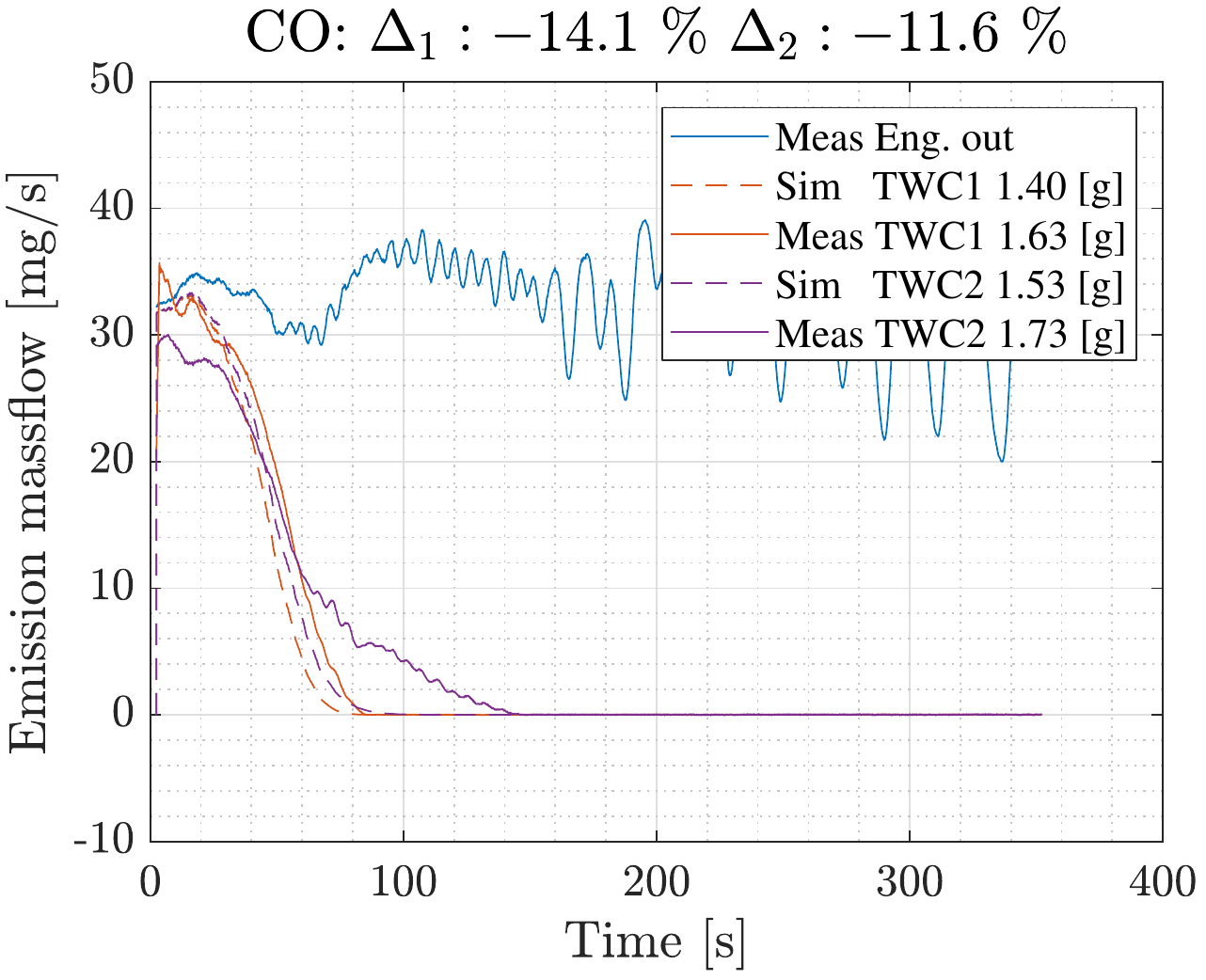}\\
	\includegraphics[width=\columnwidth]{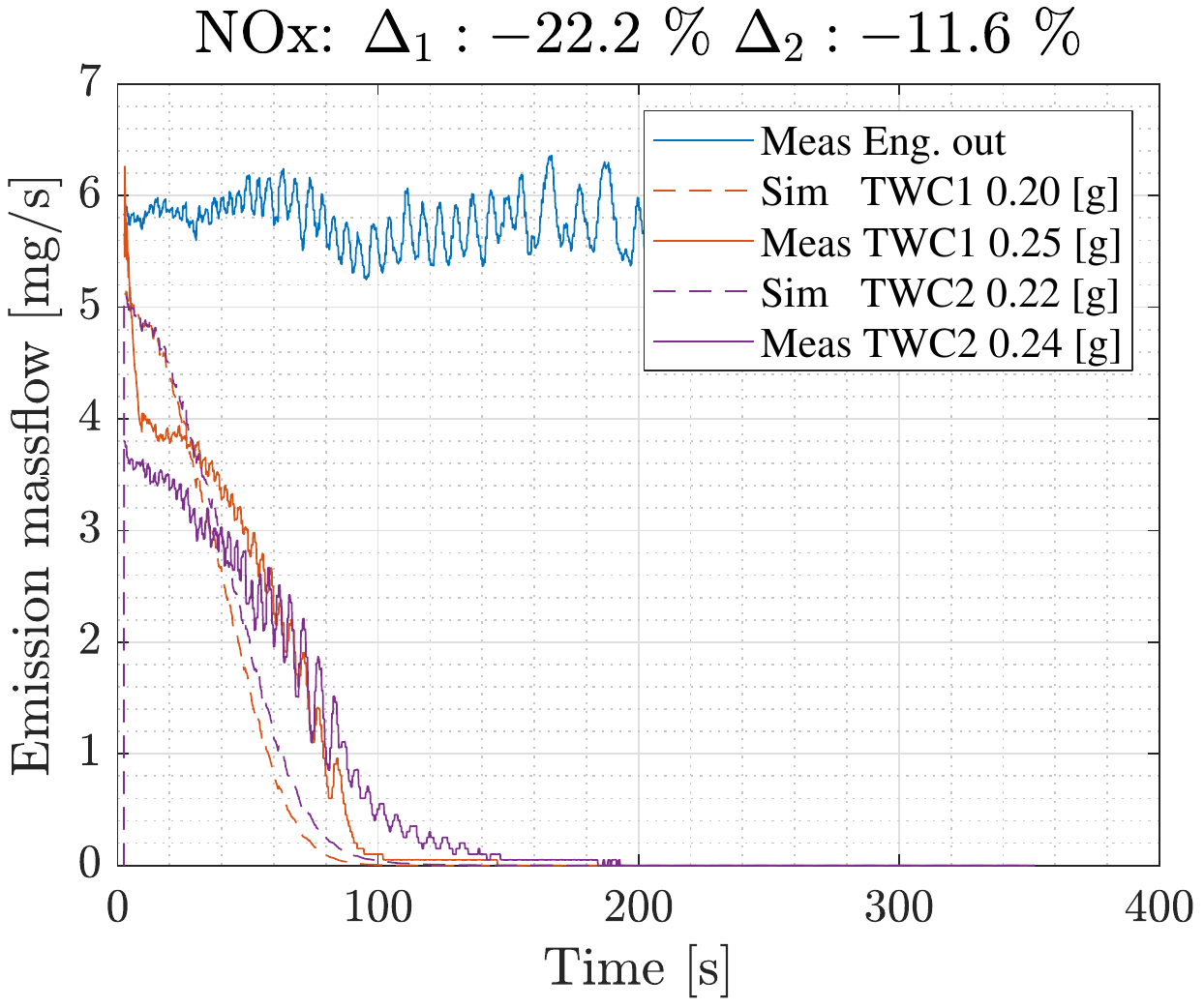}\\
	\includegraphics[width=\columnwidth]{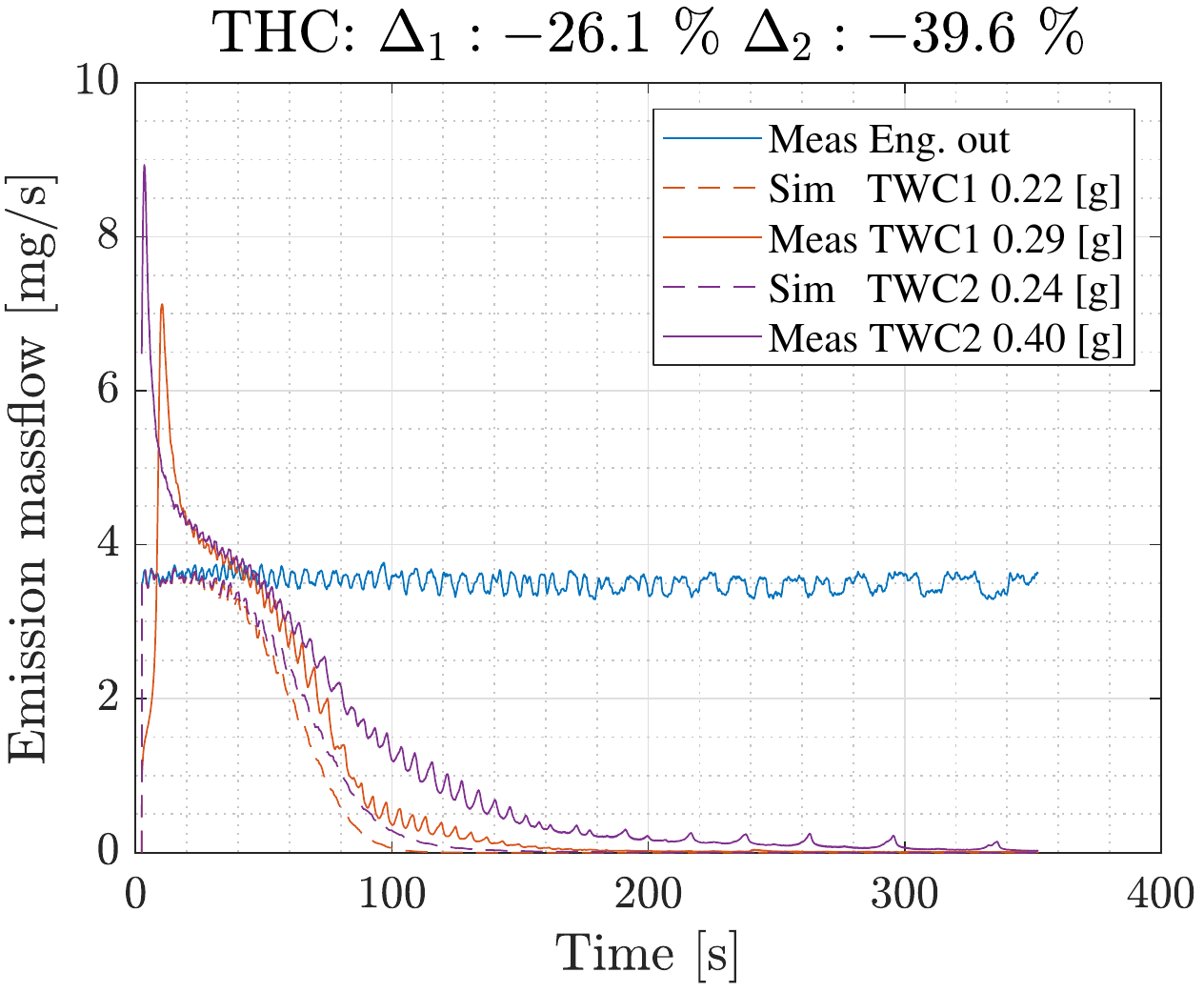}
	\caption{Measured and simulated emissions for each emission species at a 1000\,rpm, 2\,bar BMEP load point.\label{fig:delta-plot}}
\end{figure}

An illustration of the simulated and measured cold-start emissions is shown in \cref{fig:delta-plot} for the 1000\,rpm, 2\,bar BMEP load point, along with the associated cumulative simulation error. The figure indicates that one significant contribution to the cumulative error is due to inaccuracies in the measurement equipment (which is primarily designed for analyzing stationary operation). This is most clearly seen during the first 20 seconds of operation for \ce{THC}, where the measured emissions are significantly larger than than the engine-out emissions. It is possible that the \ce{NO_x} emissions are also incorrectly measured, as the measured emissions are only half of the engine-out emissions after 3-4 seconds (while light-off occurs after approximately 60 seconds at this load point). However, it is also plausible that the low \ce{NO_x} emissions are correctly measured and this anomaly is instead due to unmodeled adsorption in the TWC.

\begin{figure}
	\centering
	\includegraphics[width=\columnwidth]{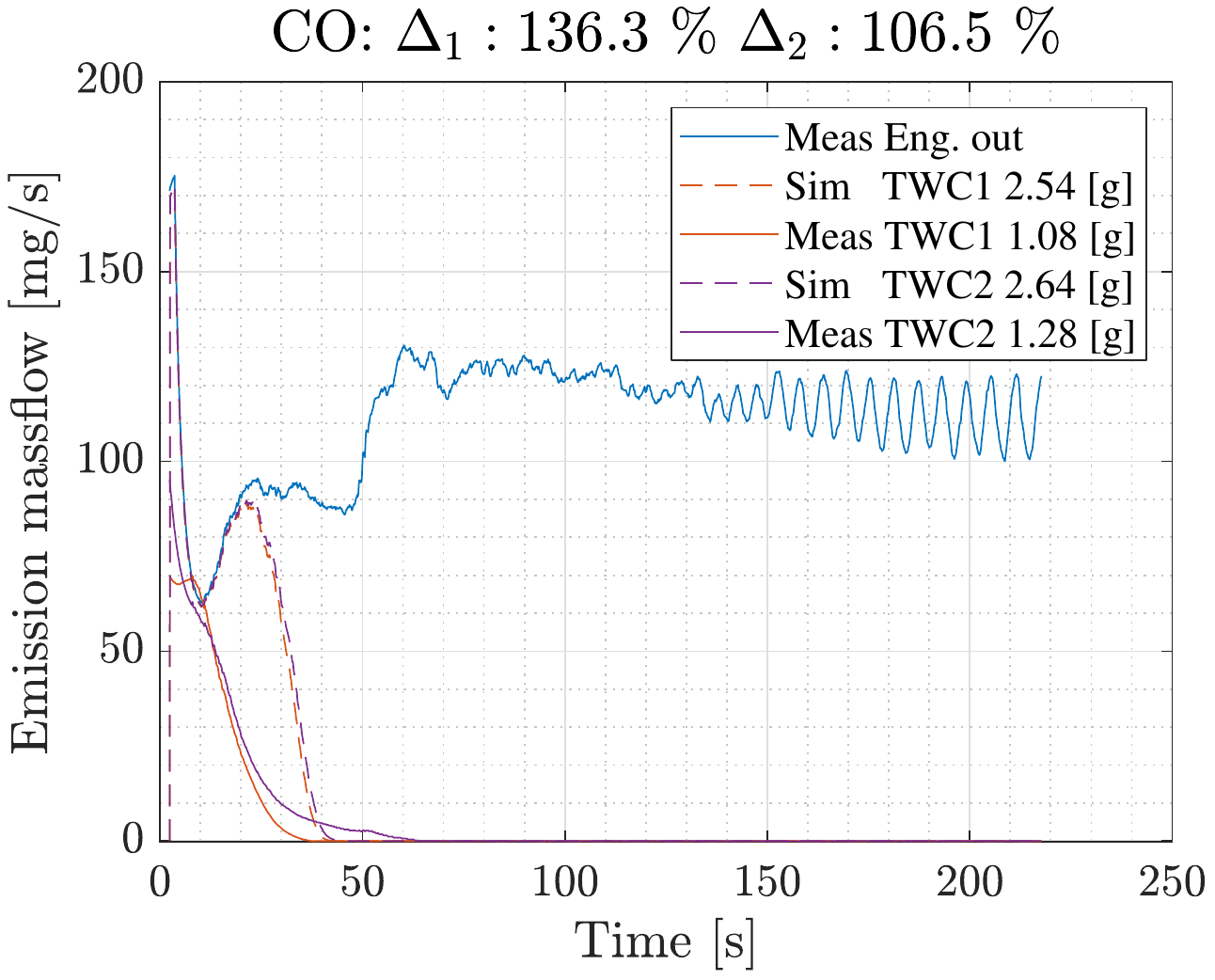}
	\caption{Measured and simulated emissions for \ce{CO} emissions at the 1500\,rpm, 5\,bar BMEP load point. Note the large increase in engine emissions after 10\,seconds, giving rise to an increase in simulated emissions without an associated increase in measured emissions. \label{fig:delta-plot-bad} }
\end{figure}

A table listing the relative cumulative tailpipe error is shown in \cref{tab:delta-accuracy} for each load point. The validation dataset (the lower half) displays an accuracy comparable to the training dataset (the upper half), indicating that the model is not overfitted. Furthermore, though there is a significant degree of variability between the measured and predicted cumulative emissions this is only somewhat worse than a significantly more complex model \cite{ramanathan_kinetic_2011} which displays a typical cumulative error on the order of $\pm20\%$ to $\pm50\%$. Furthermore, we hypothesize that the most significant outliers (e.g.~the 1500\,rpm, 2\,bar load point) are to some extent due to process variability and/or measurement error. Consulting the time evolution for this load point (\cref{fig:delta-plot-bad}, here shown for \ce{CO} emissions) indicates that this can be a factor, as the majority of the modeling error arises after 10\,seconds when the engine-out emissions significantly increase but a similar increase is not seen in the measured emissions.
 
\begin{table}
	\centering
	\caption{Relative model accuracy for training (upper half) and validation (lower half) load points.\label{tab:delta-accuracy}}
	\begin{tabular}{c|c|c|c|c}      
		Speed & BMEP & $\Delta_2^{\ce{CO}}$ & $\Delta_2^{\ce{NO_x}}$ & $\Delta_2^{\ce{THC}}$ \\
		{[RPM]} & [bar] & [-] & [-] & [-]\\
		\hline
		\hline
		997 & 1.99 & -11.6\% & -11.6\% & -39.6\% \\  
		\hline                                 
		1500 & 4.84 & +106.5\% & -20.8\% & +66.0\% \\
		\hline                                 
		2000 & 2.07 & +44.8\% & +34.2\% & -7.8\% \\  
		\hline                                 
		3000 & 8.17 & +39.5\% & +38.9\% & -18.2\% \\ 
		\hline                                 
		1500 & 8.09 & +72.7\% & +56.4\% & -10.7\% \\ 
		\hline
		\hline
		998 & 4.93 & +2.0\% & -19.5\% & -23.4\% \\  
		\hline                                
		1500 & 2.07 & +38.6\% & +41.7\% & -10.8\% \\
		\hline                                
		2010 & 4.93 & +73.3\% & +75.8\% & +13.1\% \\
		\hline                                
		999 & 7.96 & +81.2\% & +65.8\% & -13.8\% \\ 
		\hline                                
		2000 & 8.08 & +63.0\% & +48.6\% & -26.7\% \\
		\hline
	\end{tabular}    
\end{table}

\section{Optimal Control}\label{sec:opt-ctrl}
We will here illustrate optimal cold-start control as one application of the presented model. Specifically, we generate an optimal state-feedback controller suitable for on-line operation that balances the combustion engine's fuel efficiency and tailpipe emissions. We model the combustion engine exhaust using a static mean-value engine model, and allow the controller to freely choose the engine's speed, BMEP, and spark angle \cite{chan_significance_1996}. We will then use the controller to generate simulated cold-start temperature and emission trajectories and compare the results for different weightings of fuel efficiency and tailpipe emissions.

\subsection{Problem formulation}
We introduce the optimal control problem as
\begin{subequations}
	\begin{align}
	J^* &= \min_u \lim_{K \to \infty} \sum_{k=0}^K \mathrm{BSFC}(k) + \Lambda^T \cdot
	\begin{bmatrix}
	\dot{m}_\mathrm{tp}^\mathrm{\ce{CO}}(k)\\ 
	\dot{m}_\mathrm{tp}^\mathrm{\ce{NO_x}}(k)\\ 
	\dot{m}_\mathrm{tp}^\mathrm{\ce{THC}}(k)
	\end{bmatrix} \label{eq:opt-cost}
	\shortintertext{subject to}
	x(k+1) &= f_d(x(k),u(k)) \\
	g(x(k),u(k)) &\leq 0\,.
	\end{align}
	\label{eq:opt}
\end{subequations}
Here, $x$ is a state vector corresponding to both TWC's (i.e. $[T_\mathrm{TWC1}; T_\mathrm{TWC2}]$), $u$ is a discrete control variable corresponding to the requested operating point of the combustion engine (i.e. an integer value that indexes the operating points in \cref{tab:steady-state-operation-points}), $\mathrm{BSFC}(k)$ is the mean BSFC associated with operating point $u(k)$, $\Lambda$ is a $3 \times 1$ tuning parameter that balances the relative weight given to fuel-efficient operation and minimizing emissions (where smaller $\Lambda$ prioritizes the BSFC and larger $\Lambda$ prioritizes the level of emissions), $f_d$ is the system dynamics given by solving \cref{eq:diff-eq} for a given sample time for each TWC, and $g$ is a constraint function that bounds $u$ to the integer values that index the tested operating points and bounds $x$ to safe TWC temperatures.

The cost function \cref{eq:opt-cost} is specifically formulated to be of form $a + \Lambda^T \cdot b$, as this is equivalent to minimizing $a$ while limiting $b \leq B$, i.e.~minimizing the average BSFC while limiting the vector of cumulative emissions to a given level. The same structure is also commonly seen in Equivalent Consumption Minimization Strategy (ECMS) controllers \cite{MUSARDO2005509, 1002989} for the equivalent purpose balancing fuel consumption and electric energy consumption. We here notationally use $\Lambda$ rather than $\lambda$ to avoid confusion with conventional notation where $\lambda$ is used to denote the air-fuel ratio.

Note that \cref{eq:opt-cost} is formulated as an undiscounted infinite-horizon problem, as this penalizes the BSFC and emissions without requiring time to heat the TWC to be explicitly specific or known beforehand. Furthermore, by permitting the engine to operate at any of the points in \cref{tab:steady-state-operation-points} we also allow the engine power to freely vary. The hybrid vehicle cold-start problem is one example of an application that is well-suited to this cost formulation, as the electric machine can typically either supply or consume the difference between the combustion engine power and traction power.

We have solved \cref{eq:opt} using a method developed by the authors \cite{lock2021undiscounted} based on approximate dynamic programming and similar to policy iteration methods. The method, \emph{Undiscounted Control Policy generation by Approximate Dynamic Programming} (UCPADP) extends on existing approximate dynamic programming policy iteration methods by allowing for undiscounted problem formulations, i.e.~infinite-horizon problems where the cost function does not decay with increasing $k$. In principle, we can solve \cref{eq:opt} without using the above method by setting $K$ to a sufficiently large value and using a conventional ADP method \cite{Bertsekas2017, Bertsekas2012} to generate a solution. However, it is difficult to manually determine a sufficiently but not excessively large value $K$. Conveniently, the UCPADP method also returns a sufficient horizon, which for the specific TWC and cost formulation studied here was found to be 145 seconds.

One major benefit with UCPADP and other policy iteration methods is that the optimal control signal can be represented as a control law, i.e.~the optimal control signal can be simply tabulated by the state values. This implies that a controller can be implemented by simply looking up the optimal control for the current state. However, this does require knowledge of the current system state, either by direct measurement or by a state observer that estimates the system state. Furthermore, note that the selection of the engine's operating point is not formulated as a dynamic problem, which in principle implies that there is no cost associated with rapidly changing the engine's operating point. Though our solutions did not exhibit very impractical operating point changes, we have in our presented results applied a 5-second rolling average filter to the engine's target speed, BMEP, and spark angle. Though this makes for solutions that are not optimal with respect to \cref{eq:opt}, we believe that this results in more suitable engine behavior with somewhat damped transients.

Solving \cref{eq:opt} using an ADP methods first requires the state and control variables to be discretized. As $u$ is inherently discrete (indexing the operating points in \cref{tab:steady-state-operation-points}) we only need to discretize the states $x$. A denser discretization will give a solution closer to the true optimal solution, but at cost of increased memory and computational demand. We have chosen to discretize the states in TWC1 as
\begin{subequations}
\begin{align}
	T_1 &= [0, 25, 50, 75, \dots, 900]\\
	T_2 &= [0, 100, 200, \dots, 900]\\
	T_3 &= [0, 100, 200, \dots, 900]\\
	\Delta_T &= [-200, 100]
\shortintertext{and for TWC2 as}
	T_1 &= [0, 100, 200, \dots, 900]\\
	\Delta_T &= [-200, 100]\,,
\end{align} \label{eq:discretization}
\end{subequations}
i.e.~we resolve the first axial slice in TWC1 with fairly high detail, while the remaining slices and $\Delta_T$ is more coarsely resolved.

\subsection{Optimal results}
We have solved \cref{eq:opt} for a range of different normalized weights $\Lambda_n$, defined element-wise as
\begin{align}
	\Lambda_n^s &= \frac{\Lambda^s}{\min \mathrm{BSFC} / \min \dot{m}_\mathrm{exh}^{s}}
\end{align}
with results listed in \cref{tab:opt-res}. We use the normalized $\Lambda_n$ for ease of reference as $\Lambda_n=[1,1,1]$ in some sense equally weighs the fuel consumption and engine-out emissions. As tailpipe emissions approach zero when the TWC heats up we can thus view $\Lambda_n=1$ as a lower bound of relevant values to consider.

\begin{table}
	\caption{Performance of optimal controller for varying $\Lambda_n$ during a cold-start (25\,\si{\degreeCelsius}, no radial distribution) and half-warm start (200\,\si{\degreeCelsius} some radial distribution). Green-colored cells indicate cases where $\Lambda_n$ penalizes only one emission species, while other are ignored (red), and can be compared with the equally-penalized case (blue). \label{tab:opt-res}}
	\begin{subtable}{\columnwidth}
		\caption{$T_{1-3,\mathrm{TWC1}} = T_{1,\mathrm{TWC2}} = 25$\,\si{\degreeCelsius}, $\Delta_{T,\mathrm{TWC1}} = \Delta_{T,\mathrm{TWC2}} = 0$\,\si{\degreeCelsius}.}
		\resizebox{\columnwidth}{!}{%
		\begin{tabular}{c|c|c|c|c|c|c|c}
			$\Lambda_n^{\ce{CO}}$ & $\Lambda_n^{\ce{THC}}$ & $\Lambda_n^{\ce{NO_x}}$ 
			& \ce{CO} & \ce{THC} & \ce{NO_x}
			& BSFC & $m_\mathrm{fuel}$ \\
			{[-]} & [-] & [-] & [\si{\milli \g}] & [\si{\milli \g}] & [\si{\milli \g}] & [\si{\g \per {\kilo \watt \hour}}] & [\si{\g}]\\
			\hline
			\hline
			0 & 0 & 0 & 521 & 19.5 & 262 & 250 & 165 \\
			\hline
			1 & 1 & 1 & 521 & 19.5 & 262 & 250 & 165 \\
			\hline
		    10 & 10 & 10 & \cellcolor{blue!20} 423 & \cellcolor{blue!20} 18.7 & \cellcolor{blue!20} 153 & 252 & 155 \\
			0 & 0 & 10 & \cellcolor{red!20} 423 & \cellcolor{red!20} 18.6 & \cellcolor{green!20} 152 & 252 & 155 \\
			0 & 10 & 0 & \cellcolor{red!20} 521 & \cellcolor{green!20} 19.5 & \cellcolor{red!20} 262 & 250 & 165 \\
			10 & 0 & 0 & \cellcolor{green!20} 388 & \cellcolor{red!20} 21.2 & \cellcolor{red!20} 239 & 250 & 162 \\
			\hline                                                              
			$10^2$ & $10^2$ & $10^2$ & \cellcolor{blue!20} 526 & \cellcolor{blue!20} 35.0 & \cellcolor{blue!20} 41 & 262 & 130 \\
			0 & 0 & $10^2$ & \cellcolor{red!20} 687 & \cellcolor{red!20} 42.6 & \cellcolor{green!20} 38 & 262 & 126 \\
			0 & $10^2$ & 0 & \cellcolor{red!20} 494 & \cellcolor{green!20} 19.7 & \cellcolor{red!20} 242 & 250 & 162 \\
			$10^2$ & 0 & 0 & \cellcolor{green!20} 340 & \cellcolor{red!20} 24.6 & \cellcolor{red!20} 229 & 251 & 157 \\
			\hline
			$10^3$ & $10^3$ & $10^3$ & \cellcolor{blue!20} 529 & \cellcolor{blue!20} 32.7 & \cellcolor{blue!20} 29 & 271 & 117 \\
			0 & 0 & $10^3$ & \cellcolor{red!20} 530 & \cellcolor{red!20} 32.9 & \cellcolor{green!20} 29 & 270 & 116 \\
			0 & $10^3$ & 0 & \cellcolor{red!20} 432 & \cellcolor{green!20} 11.5 & \cellcolor{red!20} 151 & 254 & 158 \\
			$10^3$ & 0 & 0 & \cellcolor{green!20} 338 & \cellcolor{red!20} 33.4 & \cellcolor{red!20} 245 & 253 & 143 \\
			\hline
			$10^4$ & $10^4$ & $10^4$ & \cellcolor{blue!20} 526 & \cellcolor{blue!20} 37.0 & \cellcolor{blue!20} 26 & 298 & 90 \\
			0 & 0 & $10^4$ & \cellcolor{red!20} 668 & \cellcolor{red!20} 85.1 & \cellcolor{green!20} 15 & 356 & 60 \\
			0 & $10^4$ & 0 & \cellcolor{red!20} 433 & \cellcolor{green!20} 11.4 & \cellcolor{red!20} 145 & 255 & 155 \\
			$10^4$ & 0 & 0 & \cellcolor{green!20} 340 & \cellcolor{red!20} 33.0 & \cellcolor{red!20} 232 & 255 & 140 \\
		\end{tabular} 
		}
	\end{subtable}
	\\
	\\
	
	\begin{subtable}{\columnwidth}
		\caption{$T_{1-3,\mathrm{TWC1}} = T_{1,\mathrm{TWC2}} = 200$\,\si{\degreeCelsius}, $\Delta_{T,\mathrm{TWC1}} = \Delta_{T,\mathrm{TWC2}} = -50$\,\si{\degreeCelsius}.}
		\resizebox{\columnwidth}{!}{%
		\begin{tabular}{c|c|c|c|c|c|c|c}
			$\Lambda_n^{\ce{CO}}$ & $\Lambda_n^{\ce{THC}}$ & $\Lambda_n^{\ce{NO_x}}$ 
			& \ce{CO} & \ce{THC} & \ce{NO_x}
			& BSFC & $m_\mathrm{fuel}$ \\
			{[-]} & [-] & [-] & [\si{\milli \g}] & [\si{\milli \g}] & [\si{\milli \g}] & [\si{\g \per {\kilo \watt \hour}}] & [\si{\g}]\\
			\hline
			\hline
			0 & 0 & 0 & 265 & 11.9 & 153 & 250 & 165 \\
			\hline
			1 & 1 & 1 & 265 & 11.9 & 153 & 250 & 165 \\
			\hline
			10 & 10 & 10 & \cellcolor{blue!20} 188 & \cellcolor{blue!20} 14.2 & \cellcolor{blue!20} 100 & 251 & 155 \\                                                            
			0 & 0 & 10 & \cellcolor{red!20}187 & \cellcolor{red!20}14.1 & \cellcolor{green!20}99 & 251 & 155 \\                                                               
			0 & 10 & 0 & \cellcolor{red!20}266 & \cellcolor{green!20}11.9 & \cellcolor{red!20}153 & 250 & 165 \\                                                              
			10 & 0 & 0 & \cellcolor{green!20}207 & \cellcolor{red!20}14.3 & \cellcolor{red!20}129 & 250 & 162 \\
			\hline                                                           
			$10^2$ & $10^2$ & $10^2$ & \cellcolor{blue!20} 219 & \cellcolor{blue!20} 19.5 & \cellcolor{blue!20} 20 & 258 & 138 \\                                                          
			0 & 0 & $10^2$ & \cellcolor{red!20}271 & \cellcolor{red!20}22.8 & \cellcolor{green!20}21 & 257 & 137 \\                                                              
			0 & $10^2$ & 0 & \cellcolor{red!20}219 & \cellcolor{green!20}12.1 & \cellcolor{red!20}126 & 250 & 161 \\                                                             
			$10^2$ & 0 & 0 & \cellcolor{green!20}167 & \cellcolor{red!20}23.0 & \cellcolor{red!20} 120 & 251 & 149 \\
			\hline
			$10^3$ & $10^3$ & $10^3$ & \cellcolor{blue!20} 212 & \cellcolor{blue!20} 17.4 & \cellcolor{blue!20} 16 & 267 & 125 \\                                                       
			0 & 0 & $10^3$ & \cellcolor{red!20}235 & \cellcolor{red!20}19.4 & \cellcolor{green!20}13 & 264 & 126 \\                                                             
			0 & $10^3$ & 0 & \cellcolor{red!20}201 & \cellcolor{green!20}9.0 & \cellcolor{red!20}80 & 253 & 157 \\                                                              
			$10^3$ & 0 & 0 & \cellcolor{green!20}131 & \cellcolor{red!20}32.0 & \cellcolor{red!20}122 & 252 & 132 \\
			\hline
			$10^4$ & $10^4$ & $10^4$ & \cellcolor{blue!20} 203 & \cellcolor{blue!20} 26.9 & \cellcolor{blue!20} 12 & 303 & 86 \\                                                     
			0 & 0 & $10^4$ & \cellcolor{red!20}208 & \cellcolor{red!20} 35.3 & \cellcolor{green!20}6 & 322 & 75 \\                                                              
			0 & $10^4$ & 0 & \cellcolor{red!20}212 & \cellcolor{green!20}9.9 & \cellcolor{red!20}76 & 255 & 152 \\                                                             
			$10^4$ & 0 & 0 & \cellcolor{green!20}168 & \cellcolor{red!20}36.0 & \cellcolor{red!20}81 & 265 & 117 \\
			\hline
		\end{tabular}  
		}
	\end{subtable}
\end{table} 

We have simulated the performance of the optimal controller and list the cumulative emissions, fuel efficiency, and consumed fuel in \cref{tab:opt-res} for several different $\Lambda_n$ and two initial conditions. As expected, with increasing $\Lambda_n$ the sum of penalized emissions decrease, while the mean BSFC increases. This data also indicates that the potential for reducing \ce{NO_x} emissions is significantly larger than \ce{CO} and \ce{THC} emissions, as shown in the last 3 rows where \ce{NO_x} emissions are reduced by 94\% compared to the unpenalized case, while \ce{CO} and \ce{THC} emissions are reduced by 35\% and 41\% respectively. Furthermore, there seems to be some degree of conflict with respect to the individual emissions, as penalizing one species tends to increase the production of others. This indicates that the solutions shown in \cref{tab:opt-res} are Pareto optimal, i.e.~a given emission species mass cannot be reduced without either increasing another species' or the BSFC.

\begin{figure}
	\centering
	\includegraphics[width=\columnwidth]{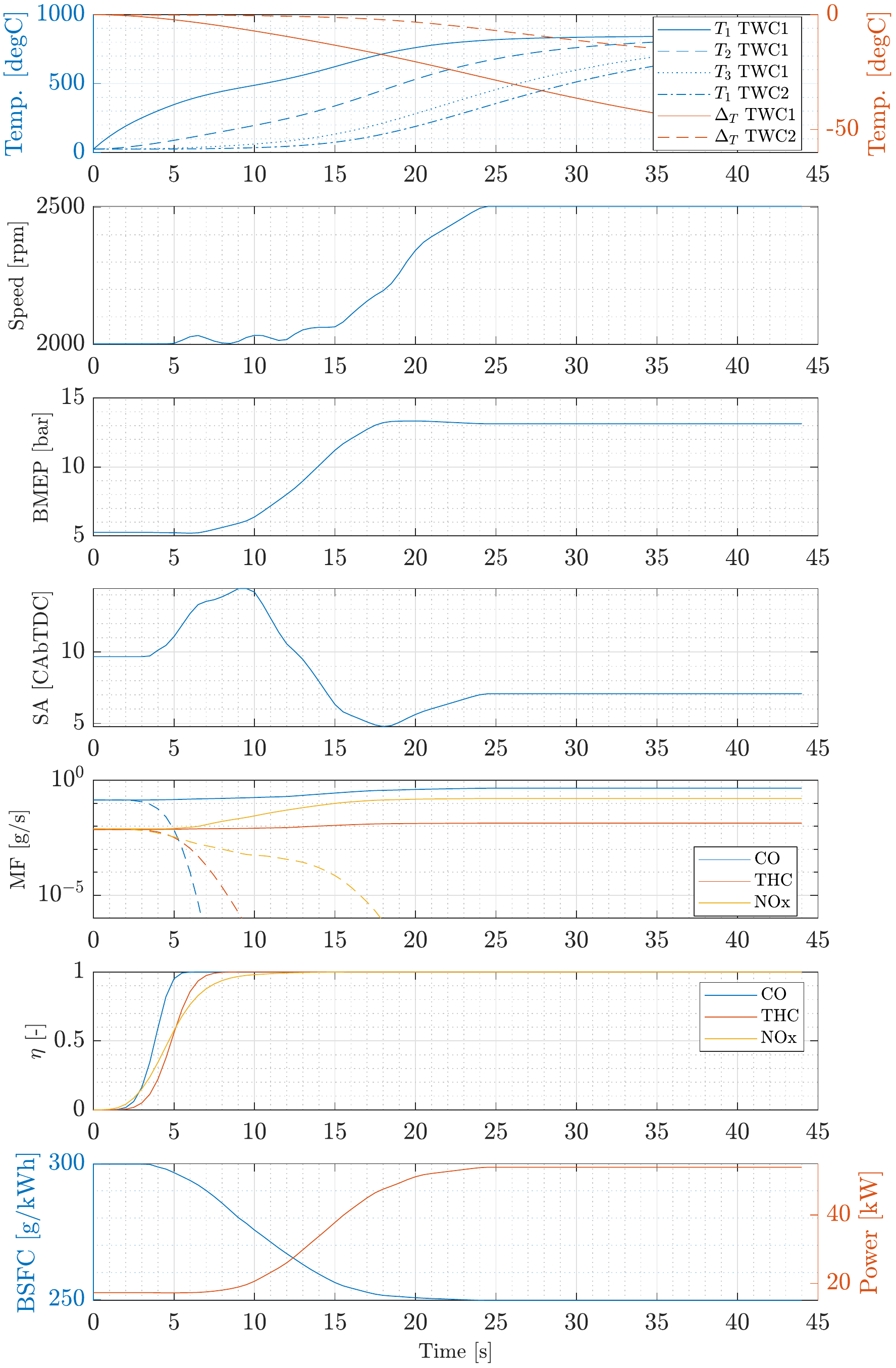}
	\caption{Simulated trajectories for $\Lambda_n = [10^2, 10^2, 10^2]$. SA indicates the spark angle, MF the engine-out emissions (solid) and tailpipe emissions (dashed), and $\eta$ the net conversion efficiency.}\label{fig:opt-traj}
\end{figure}

An illustration of the controller's time-evolution is shown in \cref{fig:opt-traj} for $\Lambda_n = [10^2, 10^2, 10^2]$. The cold-start trajectory can be divided into three sections;

\emph{$t < 5$: Initial heating phase.} The engine-out species massflow is kept low (reducing tailpipe emissions) and the BSFC is not prioritized. At the end of this phase the first axial slice is hot enough to convert emissions at low mass-flows.

\emph{$5 < t < 25$: Intermediary phase.} With increasing conversion efficiency the engine-out emissions are gradually allowed to increase, allowing for the BSFC to be increasingly prioritized. At the end of this phase the first two axial slices are hot enough to convert emissions at the massflow associated with the minimum-BSFC operating point.

\emph{$t > 25$: Sufficiently-heated phase.} Here the TWC is sufficiently hot for operation at the minimum-BSFC operating point, which the engine is statically operated at while the TWC converts virtually all emissions. Note that the entire TWC is well above light-off after 40\,\si{\second}, i.e.~a relatively short heating interval \cite{hedinger_optimal_2017}.

These sections can largely be seen for other values of $\Lambda_n$, with shorter times allocated to the initial- and intermediary heating phases as $\Lambda_n$ decreases and longer times with larger $\Lambda_n$.

\begin{figure}
	\centering
	\includegraphics[width=\columnwidth]{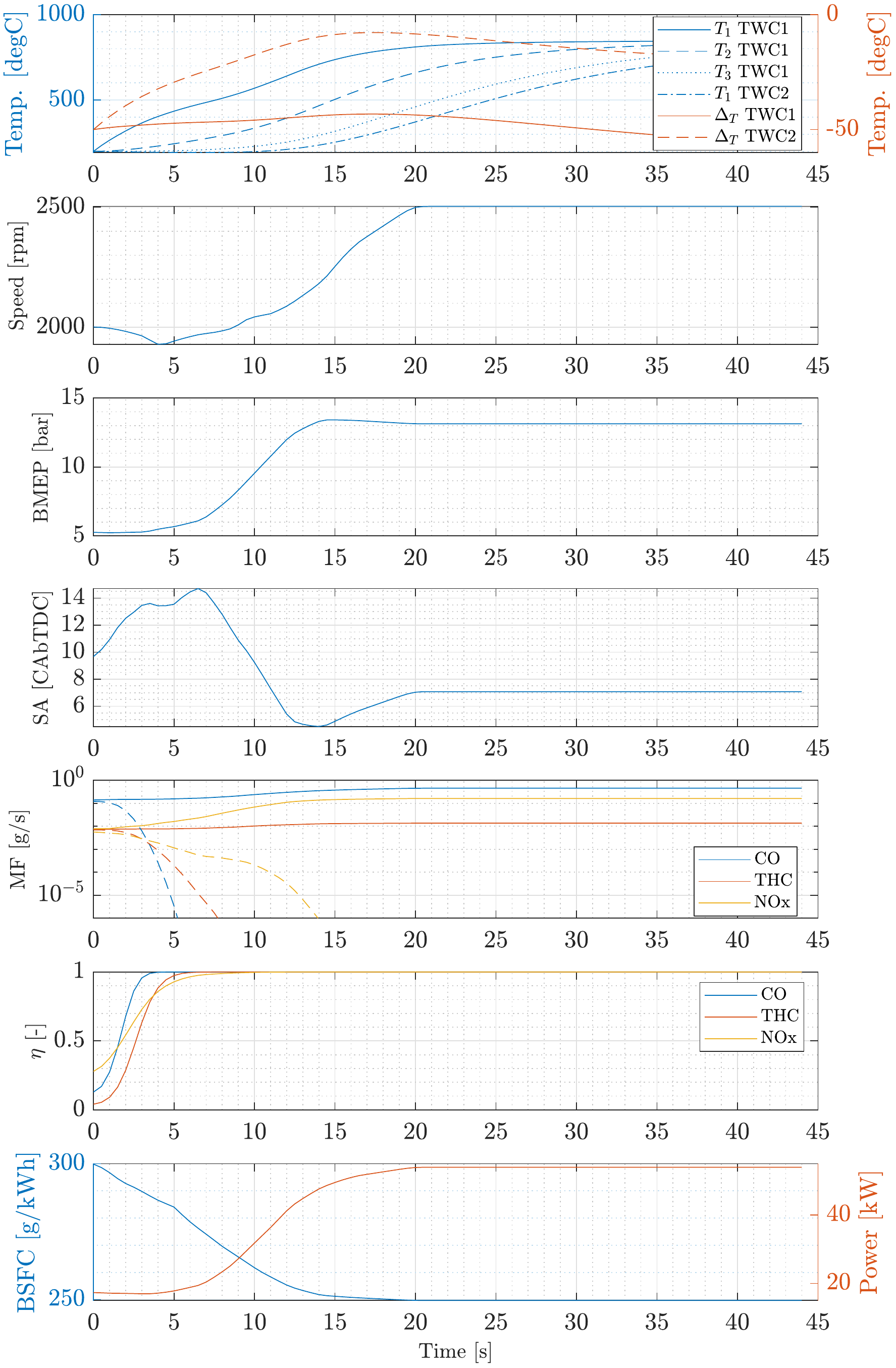}
	\caption{Simulated trajectories for $\Lambda_n = [10^2, 10^2, 10^2]$ with an initial condition $T_{1-3,\mathrm{TWC1}} = T_{1,\mathrm{TWC2}} = 200$\,\si{\degreeCelsius}, $\Delta_{T,\mathrm{TWC1}} = \Delta_{T,\mathrm{TWC2}} = -50$\,\si{\degreeCelsius}.}\label{fig:opt-traj-warmstart}
\end{figure}

Though an open-loop control scheme could be easily implemented, i.e.~by ``playing back'' the speed, BMEP, and spark angle trajectory shown in \cref{fig:opt-traj} without any temperature feedback, this type of controller is potentially sensitive to system variations. This includes both the initial temperature of the TWC (where a hotter initial condition will reach light-off more quickly) as well as variations in the exhaust gas temperature due to the fuel's composition, combustion variability, and so on. We found that the optimal control trajectory for half-warm starts is for some values of $\Lambda_n$ nearly identical to a time-shifted version of \cref{fig:opt-traj}, as exemplified in \cref{fig:opt-traj-warmstart}, while others displayed significant differences. This indicates the potential for implementing a quasi-optimal open-loop heating strategy for some $\Lambda_n$.

\subsection{Comparison to suboptimal control}\label{subsec:subopt-control}
\begin{figure}
	\centering
	\includegraphics[width=\columnwidth]{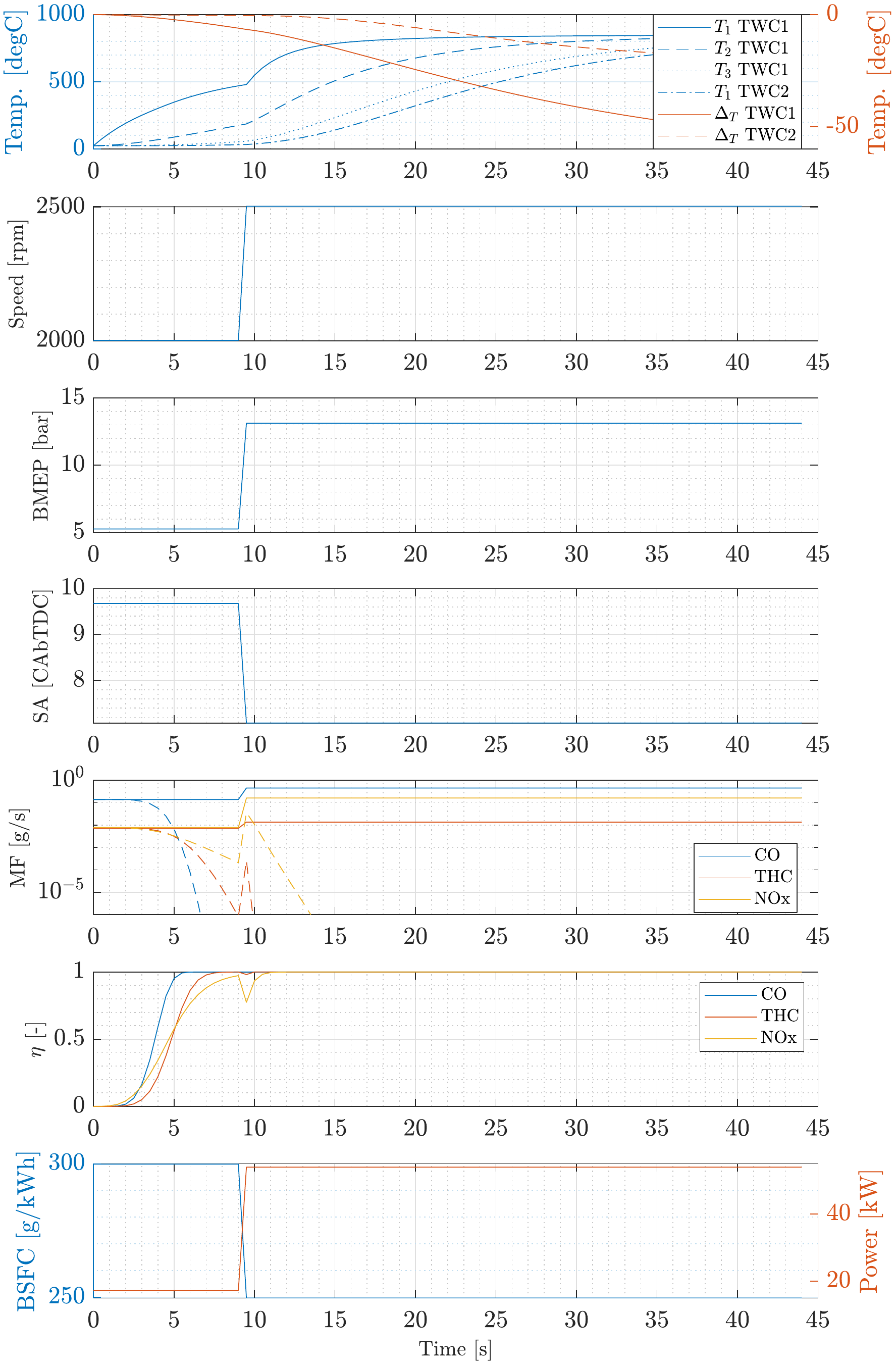}
	\caption{Simulated trajectory for the suboptimal controller. Note the significant \ce{NO_x} slip at $t=10$\,\si{\second}, which is not present in \cref{fig:opt-traj}.}\label{fig:subopt-traj}
\end{figure}

We have compared the optimal controller with a traditional suboptimal heating strategy. The suboptimal controller was defined such that combustion engine is run at a constant operating point for $t'$ seconds, and then switches to the minimum-BSFC operating point. The initial operating point was chosen to be the same as the point chosen by the optimal controller at $t=0$, and $t'$ selected to give the same average BSFC as the case for $\Lambda_n = [10^2, 10^2, 10^2]$. A comparison of the optimal and suboptimal controllers is listed in \cref{tab:subopt-ctrl}, and the time-evolution is shown in \cref{fig:subopt-traj}. The \ce{CO} and \ce{THC} emissions are virtually identical, but the \ce{NO_x} emissions are reduced by 35\% in the optimal controller. We can see the source of this in \cref{fig:subopt-traj}, where there is significant \ce{NO_x} slip at $t \in [9,13]$ when the engine transitions from the heating phase to the minimum-BSFC operation phase. Though the heating phase could be extended, this would be at cost of reduced average BSFC.

We can also compare the optimal and suboptimal controllers for the half-warm start case. It we consider the same optimal and suboptimal controllers, \cref{tab:opt-res} indicates that the optimal controller attains a mean BSFC of 258\,\si{\g \per {\kilo \watt \hour}}, in comparison to the suboptimal controller's 262\,\si{\g \per {\kilo \watt \hour}}. This difference corresponds to a 33\% reduction relative to the minimum BSFC of 250\,\si{\g \per {\kilo \watt \hour}}, indicating the potential for fuel savings by using a closed-loop cold-start strategy.

\begin{table}
	\centering
	\caption{Comparison of the optimal (here shown for $\Lambda_n = [10^2 10^2 10^2]$) and sub-optimal controllers and the relative reduction in emissions for the optimal controller. \label{tab:subopt-ctrl}}
	\begin{tabular}{c|c|c|c|c|c}
		Cont & \ce{CO} & \ce{THC} & \ce{NO_x} & BSFC & $m_\mathrm{fuel}$\\
		{[-]} & [\si{\milli \g}] & [\si{\milli \g}] & [\si{\milli \g}] & [\si{\g \per {\kilo \watt \hour}}] & [\si{\g}]\\
		\hline
		Optimal & 526 & 35.0 & 41.0 & 262 & 130 \\
		\hline
		Suboptimal & 525 & 34.9 & 62.9 & 262 & 143 \\
		\hline
		Difference & -0.2\% & -0.3\% & 34.8\% & - & - \\
		\hline
	\end{tabular}
\end{table}

\subsection{Memory footprint}
Though manageable in a PC, the memory demand associated with the discretization in \cref{eq:discretization} (with 148,000 permutations) can be problematic in an ECU that has a wide range of other tasks to perform. However, we can apply a simple space reduction scheme to significantly reduce the used memory. Rather than store the full discretization, we can reduce the number of stored elements by letting successive axial slices only be resolved for temperatures equal to or below the preceding slices (rounded up to the nearest discretized value), i.e.~storing a table of form similar to that in \cref{tab:state-lut}.

\begin{table}
	\centering
	\caption{Representative table of stored states. Here sorted from first to last column in ascending order.}
	\label{tab:state-lut}
	\begin{tabular}{c|c|c|c|c|c}
		\multicolumn{4}{c|}{TWC1} & \multicolumn{2}{c}{TWC2}\\
		$T_1$ & $T_2$ & $T_3$ & $\Delta_T$ & $T_1$ & $\Delta_T$\\
		\hline
		0 & 0 & 0 & -200 & 0 & -200 \\
		\hline
		0 & 0 & 0 & -200 & 0 & 100 \\
		\hline
		0 & 0 & 0 & -200 & 100 & -200 \\
		\hline
		$\vdots$ & 	$\vdots$ & 	$\vdots$ & 	$\vdots$ & 	$\vdots$ & 	$\vdots$\\
		\hline
		250 & 300 & 300 & 100 & 200 & -200\\
		\hline
		250 & 300 & 300 & 100 & 200 & 100\\
		\hline
		250 & 300 & 300 & 100 & 300 & -200\\
		\hline
		250 & 300 & 300 & 100 & 300 & 100\\
		\hline
		275 & 0 & 0 & -200 & 0 & -200\\
		\hline
		275 & 0 & 0 & -200 & 0 & 100\\
		\hline
		275 & 0 & 0 & -200 & 100 & -200\\
		\hline
		$\vdots$ & 	$\vdots$ & 	$\vdots$ & 	$\vdots$ & 	$\vdots$ & 	$\vdots$\\
		\hline
		900 & 900 & 900 & 100 & 900 & 100\\
		\hline
	\end{tabular}
\end{table}

Furthermore, we can reduce the resolution of $T_1$ for temperatures significantly below and above light-off by resolving $T_1$ with 100\,\si{\degreeCelsius} increments for temperatures below 100\,\si{\degreeCelsius} and above 400\,\si{\degreeCelsius}, i.e.
\begin{align}
T_1 &= [0, 100, 125, 150, \dots, 350, 375, 400, 500, \dots, 900]\,.
\end{align}
Reducing the range of considered values in this manner reduces the number of stored states to 9500 permutations. Each state permutation is associated with an optimal engine speed, BMEP, and spark angle. Assuming 4 bits of information are allocated for each parameter (allowing resolving the speed, BMEP and spark angle to 16 different and independent values) gives a total storage requirement of 12 bits per state permutation, for a total non-volatile memory requirement of $12/8\cdot9500 \approx 13.9$\,KiB, which is feasible with existing ECU hardware. More sophisticated compression schemes have the potential to further reduce the required memory, for instance by using decision trees to avoid the need to exhaustively storing every state permutation in regions where the optimal control is constant.

\section{Conclusions}
In this paper we have extended a physics-based TWC model previously presented by the authors \cite{lock2021} suited for on-line optimal control. The previously presented model resolves both axial and radial temperature variations while limiting the number of state variables, allowing for use with optimal control methods that construct an optimal control policy (e.g.~nonlinear state-feedback and explicit MPC). In this paper we extended the model to support varying axial discretization lengths, use tuning parameters expressed in well-known SI units, model heat generation by the oxidation of hydrogen, consider a TWC consisting of two separate monoliths of different construction, and use a more rigorous evaluation method with separate tuning and validation datasets. Finally, we have used the model to generate a near-optimal controller \cite{lock2021undiscounted} that can easily be implemented in existing ECU hardware, requiring no more than 13.9\,KiB (14250 bytes) of nonvolatile memory and at virtually no computational cost (as the optimal control is given by a simple linear interpolation operation and a linear rolling average filter). The specific construction of the cost function allows for systematically trading off fuel consumption and each individual emission species, giving the ability to tune the cold-start controller to minimize fuel consumption while individually limiting the specific level \ce{CO}, \ce{THC}, and \ce{NO_x} emissions.

Our experimental study, though limited by the measurement equipment, shows the potential for use both for off-line simulation as well as for generating a near-optimal cold-start controller. Though we experimentally studied the case of a warm engine and a cold TWC for improved experimental repeatability, we hypothesize that the controller can be extended to the cold engine case by a suitable update to the combustion engine exhaust model. Though the measured predictive accuracy is fairly low (with cumulative cold-start emissions typically estimated at -20\% to +80\% of the measured emissions), it is likely that the experimental setup significantly contributes to this error. The second monolith is situated after a sharp bend, giving a temperature distribution that is not particularly well-captured with an axi-radial model. However, as the majority of the emissions can be converted in the first monolith for low to moderate load-points this inaccuracy might not be of great importance. It may therefore be a prudent design decision to solely model and control the first monolith dynamics in an effort to further reduce the memory requirements of the controller.

We have simulated the performance of the Pareto-optimal cold-start controller for several different relative weightings of fuel efficiency (BSFC) and cumulative emissions for each emissions species, i.e.~different points on the Pareto front. For one representative weighting the optimal controller gives \ce{NO_x} emissions that are 35\% lower than a traditional cold-start controller with otherwise identical BSFC and \ce{CO} and \ce{THC} emissions. This indicates that an optimal controller that is generated using the presented model has the potential to reduce the cold-start emissions, as well as allowing for systematically adjusting the trade-off between each emission species and fuel consumption. Furthermore, for some regions on the Pareto-front the optimal controllers display similar speed, load, and spark angle trajectories for varying initial TWC temperatures (up to a shift in time). It is therefore plausible that a close-to optimal controller could be implemented with only a single temperature sensor by ``playing back'' a section of the temporally-resolved optimal control trajectory based on the measured temperature.

Relevant future work includes performing an experimental study using measurement equipment more suited to transient conditions and that is capable of measuring emissions at two locations simultaneously. Furthermore, it would be prudent to experimentally validate the performance of the presented controller, which in turn requires a method for measuring or estimating the temperatures in the TWC. Additionally, we hypothesize that characterizing the engine emissions for different average air-fuel ratios near stoichiometry could allow for the optimal controller to gain an additional degree of freedom in balancing the ratio of \ce{CO} and \ce{THC} emissions to the \ce{NO_x} emissions during a cold-start.

\section*{Funding}
This work was performed within the Combustion Engine Research Center at Chalmers (CERC) with financial support from the Swedish Energy Agency.

\printbibliography

@article{lock2021,
	title = {A Control-Oriented Spatially Resolved Thermal Model of the Three-Way-Catalyst},
	url = {https://www.sae.org/publications/technical-papers/content/2021-01-0597/},
	journal = {{SAE} International},
	author = {Jonathan Lock and Kristoffer Clas\'{e}n and Jonas Sj\"{o}blom and Tomas McKelvey},
	date = {2021-04-06},
	year = "2021",
}

@misc{lock2021undiscounted,
      title={Undiscounted Control Policy Generation for Continuous-Valued Optimal Control by Approximate Dynamic Programming}, 
      author={Jonathan Lock and Tomas McKelvey},
      year={2021},
      eprint={2104.11093},
      archivePrefix={arXiv},
      primaryClass={math.OC}
}

@article{ramanathan_optimal_2004,
	title = {Optimal design of catalytic converters for minimizing cold-start emissions},
	volume = {98},
	issn = {09205861},
	url = {https://linkinghub.elsevier.com/retrieve/pii/S0920586104005139},
	doi = {10.1016/j.cattod.2004.08.003},
	pages = {357--373},
	number = {3},
	journal = {Catalysis Today},
	shortjournal = {Catalysis Today},
	author = {Ramanathan, Karthik and West, David H. and Balakotaiah, Vemuri},
	date = {2004-12},
	year = "2004",
}

@article{gao_review_2019,
	title = {Review of thermal management of catalytic converters to decrease engine emissions during cold start and warm up},
	volume = {147},
	issn = {13594311},
	url = {https://linkinghub.elsevier.com/retrieve/pii/S1359431118336081},
	doi = {10.1016/j.applthermaleng.2018.10.037},
	pages = {177--187},
	journal = {Applied Thermal Engineering},
	author = {Gao, Jianbing and Tian, Guohong and Sorniotti, Aldo and Karci, Ahu Ece and Di Palo, Raffaele},
	date = {2019-01},
	year = "2019",
	langid = {english},
}

@article{schori_optimal_2014,
	title = {Optimal catalytic converter heating in hybrid vehicles},
	journal = {{SAE} Technical Paper},
	author = {Schori, Markus and Boehme, Thomas and Jeinsch, Torsten and Schultalbers, Matthias},
	year = {2014},
}

@misc{pannag_coldstart_2009,
	title = {Coldstart {Modeling} and {Optimal} {Control} {Design} for {Automotive} {SI} engine},
	author = {Pannag, R. Sanketi},
	year = {2009},
}

@article{fiengo_control_2002,
	title = "Control of the exhaust gas emissions during the warm-up process of a {TWC}-equipped {SI} engine",
	volume = {35},
	issn = {14746670},
	url = {https://linkinghub.elsevier.com/retrieve/pii/S1474667015399390},
	doi = {10.3182/20020721-6-ES-1901.01518},
	language = {en},
	number = {1},
	journal = {IFAC Proceedings Volumes},
	author = {Fiengo, Giovanni and Glielmo, Luigi and Santini, Stefania and Serra, Gabriele},
	year = {2002},
	pages = {301--306},
}

@article{keynejad_suboptimal_2013,
	title = {Suboptimal {Cold} {Start} {Strategies} for {Spark} {Ignition} {Engines}},
	volume = {21},
	url = {http://ieeexplore.ieee.org/document/6236300/},
	doi = {10.1109/TCST.2012.2203821},
	number = {4},
	journal = {IEEE Transactions on Control Systems Technology},
	author = {Keynejad, Farzad and Manzie, Chris},
	month = jul,
	year = {2013},
	pages = {1295--1308},
}

@article{brandt_dynamic_2000,
	title = {Dynamic modeling of a three-way catalyst for {SI} engine exhaust emission control},
	volume = {8},
	issn = {10636536},
	url = {http://ieeexplore.ieee.org/document/865850/},
	doi = {10.1109/87.865850},
	language = {en},
	number = {5},
	journal = {IEEE Transactions on Control Systems Technology},
	author = {Brandt, E.P. and {Yanying Wang} and Grizzle, J.W.},
	month = sep,
	year = {2000},
	pages = {767--776},
}

@article{pattas_transient_1994,
	title = {Transient Modeling of 3-Way Catalytic Converters},
	url = {https://www.sae.org/content/940934/},
	doi = {10.4271/940934},
	journal = {International Congress \&  Exposition},
	author = {Pattas, K. N. and Stamatelos, A. M. and Pistikopoulos, P. K. and Koltsakis, G. C. and Konstandinidis, P. A. and Volpi, E. and Leveroni, E.},
	date = {1994-03-01},
	year = "1994",
}

@book{heywood1988internal,
	title = {Internal Combustion Engine Fundamentals},
	author = "John B. Heywood",
	year = {1988},
	publisher = {McGraw Hill}
}

@article{tischer_three-way-catalyst_2007,
	title = {Three-Way-Catalyst Modeling - A Comparison of 1D and 2D Simulations},
	url = {https://www.sae.org/content/2007-01-1071/},
	doi = {10.4271/2007-01-1071},
	journal = {{SAE} World Congress \&  Exhibition},
	pages = {2007--01--1071},
	author = {Tischer, Steffen and Jiang, Yi and Hughes, Katherine W. and Patil, M. D. and Murtagh, Michael},
	date = {2007-04-16},
	year = "2007",
}

@mastersthesis{svraka_model_2019,
	title = {Model Based Catalyst Control},
	author = {Svraka, Irman and Linus, {\"O}sterdahl Wetterhag},
	year = {2019},
	school = "Link{\"o}ping University"
}

@book{logan1998applied,
  title={Applied Partial Differential Equations},
  author={Logan, J.D. and Gehring, F.W. and Halmos, P.R.},
  series={Springer Undergraduate Mathematics Series},
  year={1998},
  publisher={Springer}
}

@article{hedinger_optimal_2017,
	title = {Optimal {Cold}-{Start} {Control} of a {Gasoline} {Engine}},
	volume = {10},
	issn = {1996-1073},
	url = {http://www.mdpi.com/1996-1073/10/10/1548},
	doi = {10.3390/en10101548},
	language = {en},
	number = {10},
	journal = {Energies},
	author = {Hedinger, Raffael and Elbert, Philipp and Onder, Christopher},
	month = oct,
	year = {2017},
	pages = {1548},
}

@article{zygourakis1989transient,
  title={Transient operation of monolith catalytic converters: a two-dimensional reactor model and the effects of radially nonuniform flow distributions},
  author={Zygourakis, Kyriacos},
  journal={Chemical Engineering Science},
  volume={44},
  number={9},
  pages={2075--2086},
  year={1989},
}

@techreport{braun2002three,
  title={Three-dimensional simulation of the transient behavior of a three-way catalytic converter},
  author={Braun, Joachim and Hauber, Thomas and T{\"o}bben, Heike and Windmann, Julia and Zacke, Peter and Chatterjee, Daniel and Correa, Chrys and Deutschmann, Olaf and Maier, Lubow and Tischer, Steffen and others},
  year={2002},
  institution={SAE Technical Paper}
}

@article{chen1988three,
  title={A three-dimensional model for the analysis of transient thermal and conversion characteristics of monolithic catalytic converters},
  author={Chen, David KS and Bissett, Edward J and Oh, Se H and Van Ostrom, David L},
  journal={SAE transactions},
  pages={177--189},
  year={1988},
}

@article{baba1996numerical,
  title={Numerical approach for improving the conversion characteristics of exhaust catalysts under warming-up condition},
  author={Baba, Naoki and Ohsawa, Katsuyuki and Sugiura, Shigeki},
  journal={SAE transactions},
  pages={2064--2079},
  year={1996},
}

@article{oh1982transients,
  title={Transients of monolithic catalytic converters. Response to step changes in feedstream temperature as related to controlling automobile emissions},
  author={Oh, Se H and Cavendish, James C},
  journal={Industrial \& Engineering Chemistry Product Research and Development},
  volume={21},
  number={1},
  pages={29--37},
  year={1982},
}

@article{pontikakis2004three,
  title={Three-way catalytic converter modeling as a modern engineering design tool},
  author={Pontikakis, GN and Konstantas, GS and Stamatelos, AM},
  journal={J. Eng. Gas Turbines Power},
  volume={126},
  number={4},
  pages={906--923},
  year={2004}
}

@article{real_modelling_2021,
	title = {Modelling three-way catalytic converter oriented to engine cold-start conditions},
	volume = {22},
	url = {http://journals.sagepub.com/doi/10.1177/1468087419853145},
	doi = {10.1177/1468087419853145},
	number = {2},
	journal = {International Journal of Engine Research},
	author = {Real, Marcelo and Hedinger, Raffael and Pla, Benjamín and Onder, Christopher},
	month = feb,
	year = {2021},
	pages = {640--651},
}

@article{yan_modeling_2019,
	title = {Modeling {Three}-{Way} {Catalyst} {Converters} {During} {Cold} {Starts} {And} {Potential} {Improvements}},
	url = {https://www.sae.org/content/2019-01-2326/},
	doi = {10.4271/2019-01-2326},
	author = {Yan, Xieyang and Sone, Ryota and Inoue, Ryoya and Kusaka, Jin and Umezawa, Katunori and Kondo, Yasuhiro},
	month = dec,
	year = {2019},
	pages = {2019--01--2326},
	journal = {SAE International},
}

@article{ramanathan_kinetic_2011,
	title = {Kinetic {Parameters} {Estimation} for {Three} {Way} {Catalyst} {Modeling}},
	volume = {50},
	url = {https://pubs.acs.org/doi/10.1021/ie200726j},
	doi = {10.1021/ie200726j},
	language = {en},
	number = {17},
	journal = {Industrial \& Engineering Chemistry Research},
	author = {Ramanathan, Karthik and Sharma, Chander Shekhar},
	month = sep,
	year = {2011},
	pages = {9960--9979},
}

@article{shaw_simplified_2002,
	title = {A {SIMPLIFIED} {COLDSTART} {CATALYST} {THERMAL} {MODEL} {TO} {REDUCE} {HYDROCARBON} {EMISSIONS}},
	volume = {35},
	issn = {14746670},
	url = {https://linkinghub.elsevier.com/retrieve/pii/S1474667015399407},
	doi = {10.3182/20020721-6-ES-1901.01519},
	language = {en},
	number = {1},
	journal = {IFAC Proceedings Volumes},
	author = {Shaw, Byron T. and Fischer, Gerald D. and Hedrick, J. Karl},
	year = {2002},
	pages = {307--312},
}

@article{chan_significance_1996,
	title = {The {Significance} of {High} {Value} of {Ignition} {Retard} {Control} on the {Catalyst} {Lightoff}},
	url = {https://www.sae.org/content/962077/},
	doi = {10.4271/962077},
	author = {Chan, S. H. and Zhu, J.},
	month = oct,
	year = {1996},
	pages = {962077},
	journal = {SAE International},
}

@book{Bertsekas2017,
	Author = "Dimitri P. Bertsekas",
	Publisher = "Athena Scientific",
	Title = "Dynamic Programming and Optimal Control",
	Year = 2017,
	Volume = 1,
	edition = 4
}

@book{Bertsekas2012,
	Author = "Dimitri P. Bertsekas",
	Publisher = "Athena Scientific",
	Title = "Dynamic Programming and Optimal Control",
	Year = 2012,
	Volume = 2,
	edition = 4
}

@article{michel_optimizing_2017,
	title = {Optimizing fuel consumption and pollutant emissions of gasoline-{HEV} with catalytic converter},
	language = {en},
	journal = {Control Engineering Practice},
	author = {Michel, Pierre},
	year = {2017},
	pages = {8}
}

@article{azad_determining_2012,
	title = {Determining {Model} {Accuracy} {Requirements} for {Automotive} {Engine} {Coldstart} {Hydrocarbon} {Emissions} {Control}},
	volume = {134},
	url = {https://asmedigitalcollection.asme.org/dynamicsystems/article/doi/10.1115/1.4006217/394995/Determining-Model-Accuracy-Requirements-for},
	doi = {10.1115/1.4006217},
	language = {en},
	number = {5},
	urldate = {2019-10-22},
	journal = {Journal of Dynamic Systems, Measurement, and Control},
	author = {Azad, Nasser L. and Sanketi, Pannag R. and Hedrick, J. Karl},
	month = sep,
	year = {2012},
	pages = {051002},
}

@article{zhu_development_2019,
	title = {Development of physics-based three-way catalytic converter model for real-time distributed temperature prediction using proper orthogonal decomposition and collocation},
	url = {http://journals.sagepub.com/doi/10.1177/1468087419876127},
	doi = {10.1177/1468087419876127},
	language = {en},
	urldate = {2020-02-18},
	journal = {International Journal of Engine Research},
	author = {Zhu, Zhaoxuan and Midlam-Mohler, Shawn and Canova, Marcello},
	month = sep,
	year = {2019},
	pages = {146808741987612},
}

@book{bryson1975applied,
  title={Applied optimal control: optimization, estimation and control},
  author={Bryson, Arthur Earl},
  year={1975},
  publisher={CRC Press}
}

@article{kang_detailed_2014,
	title = {Detailed reaction kinetics for double-layered Pd/Rh bimetallic {TWC} monolith catalyst},
	volume = {241},
	issn = {13858947},
	url = {https://linkinghub.elsevier.com/retrieve/pii/S1385894713016203},
	doi = {10.1016/j.cej.2013.12.039},
	pages = {273--287},
	journal = {Chemical Engineering Journal},
	shortjournal = {Chemical Engineering Journal},
	author = {Kang, Sung Bong and Han, Seok Jun and Nam, In-Sik and Cho, Byong K. and Kim, Chang Hwan and Oh, Se H.},
	urldate = {2020-08-27},
	date = {2014-04},
	year = "2014",
	langid = {english},
}

@article{MUSARDO2005509,
title = {A-ECMS: An Adaptive Algorithm for Hybrid Electric Vehicle Energy Management},
journal = {European Journal of Control},
volume = {11},
number = {4},
pages = {509-524},
year = {2005},
issn = {0947-3580},
doi = {https://doi.org/10.3166/ejc.11.509-524},
url = {https://www.sciencedirect.com/science/article/pii/S0947358005710487},
author = {Cristian Musardo and Giorgio Rizzoni and Yann Guezennec and Benedetto Staccia},
}

@article{1002989,
  author={G. {Paganelli} and S. {Delprat} and T. M. {Guerra} and J. {Rimaux} and J. J. {Santin}},
  journal={Vehicular Technology Conference. IEEE 55th Vehicular Technology Conference. VTC Spring 2002 (Cat. No.02CH37367)}, 
  title={Equivalent consumption minimization strategy for parallel hybrid powertrains}, 
  year={2002},
  volume={4},
  number={},
  pages={2076-2081 vol.4},
  doi={10.1109/VTC.2002.1002989}
}

\end{document}